\documentclass[12pt]{amsart}
\usepackage{amsfonts,amsmath,amsthm,amssymb,tikz,tikz-cd}
\usetikzlibrary{matrix}
\usepackage[all]{xy}
\usepackage[vcentermath]{youngtab}
\usepackage{hyperref}
\usepackage{epsfig}

\usepackage[margin=1.1in]{geometry}

\definecolor{darkred}{rgb}{1,0,0}

\newtheorem{theorem}{Theorem}[section]
\newtheorem{thm}[theorem]{Theorem}
\newtheorem{lem}[theorem]{Lemma}

\newtheorem{prop}[theorem]{Proposition}

\newtheorem{conj}[theorem]{Conjecture}

\theoremstyle{remark}

\theoremstyle{definition}

\newtheorem{rmk}[theorem]{Remark}
\newtheorem{example}[theorem]{Example}
\newtheorem{defn}[theorem]{Definition}

\newcommand{\bbc}{\mathbb{C}}

\newcommand{\bbp}{\mathbb{P}}

\newcommand{\bbr}{\mathbb{R}}

\newcommand{\bbt}{\mathbb{T}}

\newcommand{\bbz}{\mathbb{Z}}

\newcommand{\mca}{\mathcal{A}}

\newcommand{\mcc}{\mathcal{C}}

\newcommand{\mce}{\mathcal{E}}
\newcommand{\mcf}{\mathcal{F}}

\newcommand{\mcn}{\mathcal{N}}

\newcommand{\mcw}{\mathcal{W}}
\newcommand{\mcx}{\mathcal{X}}

\newcommand{\mfs}{\mathfrak{S}}

\renewcommand{\aa}{\alpha}
\newcommand{\bb}{\beta}

\newcommand{\ov}{\overline}
\newcommand{\eps}{\epsilon}

\newcommand{\CMG}{\mathcal{G}}

\newcommand{\la}{\langle}
\newcommand{\ra}{\rangle}
\newcommand{\injects}{\hookrightarrow}
\newcommand{\surjects}{\twoheadrightarrow}

\DeclareMathOperator{\MCG}{\textnormal{MCG}}
\DeclareMathOperator{\Gr}{\textnormal{Gr}}
\DeclareMathOperator{\tGr}{\widetilde{\textnormal{Gr}}}
\newcommand{\tQ}{\tilde{Q}}

\DeclareMathOperator{\SL}{\text{SL}}
\DeclareMathOperator{\PSL}{\text{PSL}}
\DeclareMathOperator{\Conf}{\textnormal{FG}}
\DeclareMathOperator{\FG}{\textnormal{FG}}

\DeclareMathOperator{\Aut}{Aut}

\DeclareMathOperator{\Arb}{\textnormal{Arb}}
\DeclareMathOperator{\Span}{\textnormal{span}}
\newcommand{\Le}{\textup{\protect\scalebox{-1}[1]{L}}}
\newcommand{\Imod}{I_{\textnormal{mod}}}
\newcommand{\Jmod}{J_{\textnormal{mod}}}
\newcommand{\mccmod}{\mcc_{\textnormal{mod}}}

\title[Braid group symmetries of Grassmannian cluster algebras]{Braid group symmmetries of Grassmannian cluster algebras}
\author{Chris Fraser}
\date{}
\keywords{Cluster algebra, Grassmannian, braid group, quasi-homomorphism, web.}
\subjclass[2010]{13F60}
\address{Department of Mathematical Sciences, Indiana University - Purdue University, Indianapolis.}
\email{cmfra@umich.edu}
\thanks{Portions of this work were supported by a graduate fellowship from the National Physical Science Consortium and NSF grant DMS-1361789.}

\begin{document}

\begin{abstract}
Let $\Gr^\circ(k,n) \subset \Gr(k,n)$ denote the open positroid stratum in the Grassmannian. We define an action of the extended affine $d$-strand braid group on $\Gr^\circ(k,n)$ by regular automorphisms, for $d$ the greatest common divisor of $k$ and $n$. The action is by quasi-automorphisms of the cluster structure on $\Gr^\circ(k,n)$, determining a homomorphism from the extended affine braid group to the cluster modular group for $\Gr(k,n)$. We also define a quasi-isomorphism between the Grassmannian $\Gr(k,rk)$ and the Fock-Goncharov configuration space of $2r$-tuples of affine flags for $\SL_k$. This identifies the cluster variables, clusters, and cluster modular groups, in these two cluster structures. 

Fomin and Pylyavskyy proposed a description of the cluster combinatorics for $\Gr(3,n)$ in terms of Kuperberg's basis of non-elliptic webs. As our main application, we prove many of their conjectures for $\Gr(3,9)$ and give a presentation for its cluster modular group. We establish similar results 
for $\Gr(4,8)$. These results rely on the fact that both of these Grassmannians have finite mutation type.
\end{abstract}

\maketitle
\vspace{-.3in}

\section*{Introduction}\label{secn:Intro}
In its most combinatorial formulation, the theory of cluster algebras concerns itself with identifying \emph{which} elements of a cluster algebra are cluster variables, and how these elements are grouped into clusters. Some of the most central examples of cluster algebras occur as (homogeneous) coordinate rings of important algebraic varieties in Lie theory or geometry -- Grassmannians and other partial flag varieties, double Bruhat and positroid cells, and (complexifications of) decorated Teichm\"uller spaces. A motivation for computing \emph{all of} the clusters in these important examples is that each cluster variable, moreover each \emph{monomial} in any cluster, is expected to lie in a version of a \emph{canonical basis} for the coordinate ring.  

The definition of cluster variables is recursive and somewhat technical.  It begins with an initial choice of a \emph{cluster} -- a distinguished collection of elements in the cluster algebra  -- and produces new cluster variables one at a time, by exchanging a current cluster variable for a neighboring one. The new cluster variable is defined algebraically in terms of the current cluster, in what is known as an \emph{exchange relation}. An extra layer of subtlety is provided by designating certain initial variables as \emph{frozen}: these variables are never themselves exchanged, but they appear in the exchange relations defining new cluster variables. In any given cluster, the list of exchange relations defining the neighboring cluster variables is encoded by a \emph{quiver}. Each time a cluster is exchanged for a new one, this list updates according to its own dynamical rules (in a process called \emph{quiver mutation}). 

This paper presents a non-recursive way of calculating cluster variables and clusters for Grassmannians, by an action of an appropriate group of symmetries. Let $\tGr(k,n)$ denote (the affine cone over) the Grassmannian of $k$-dimensional subspaces in~$\bbc^n$ and $\tGr^\circ(k,n) \subset \tGr(k,n)$ the Zariski-open subset cut out by the non-vanishing of the frozen variables. Let $d = \gcd(k,n)$ and $B_d$ denote the braid group on $d$ strands. We define (cf~Definition~\ref{defn:sigmai}) regular automorphisms $\sigma_i \colon \tGr^\circ(k,n) \to \tGr^\circ(k,n)$ which satisfy the braid relations. 

The pullback $\sigma_i^*$ to the coordinate ring $\bbc[\tGr(k,n)]$ is a \emph{quasi-automorphism}~\cite{Fraser} of the cluster structure, hence  $\sigma_i^*$ induces a permutation of the cluster variables, clusters, and cluster monomials, in $\bbc[\tGr(k,n)]$. The resulting $B_d$ action  
preserves (the non-frozen part of) each quiver, defining a homomorphism from $B_d$ to the \emph{cluster modular group} $\CMG = \CMG(\Gr(k,n))$, a certain symmetry group of a cluster algebra introduced by Fock and Goncharov. Together with the well-known twisted cyclic shift automorphism $\rho$ of $\Gr(k,n)$, this homomorphism can be enriched to a homomorphism $\hat{B}_{\hat{A}_{d-1}} \to \CMG$, where $\hat{B}_{\hat{A}_{d-1}}$ is the $d$-strand \emph{extended affine braid group} \cite{ExtendedBraid}. We note two extreme cases: when $k$ divides $n$, the extended affine braid group action reduces to an action by the ordinary braid group~$B_k$ (cf. Lemma~\ref{lem:consecutivecomp}). On the other hand, when $k$ and $n$ are coprime, the braid group $B_d=B_1$ is trivial and our constructions do not give rise to new symmetries. 

The homomorphism $\hat{B}_{\hat{A}_{d-1}} \to \CMG$ is not faithful: $\hat{B}_{\hat{A}_{d-1}}$ has no elements of finite order, but the element $\rho \in \hat{B}_{\hat{A}_{d-1}}$ that maps to the cyclic shift automorphism has finite order~$n$ inside $\CMG$. The element $\rho^n$ is central in $\hat{B}_{\hat{A}_{d-1}}$, and we expect that our action (cf.~Conjecture \ref{conj:groups}) is faithful modulo the center for $n>>k$.  

\vspace{.3cm}

Another important family of cluster algebras were introduced by Fock and Goncharov~\cite{ModuliSpaces} in a pioneering series of papers on higher Teichm\"uller theory. For a simple Lie group~$G$, they considered spaces of (twisted, decorated) $G$-local systems on a (bordered, marked, oriented) surface $S$. These spaces have cluster structures \cite{ModuliSpaces,LeClassicalGroups}. This paper focuses on the space $\Conf(k,r)$ arising when $G = \SL_k$ and $S$ is a disk with $r$ points on its boundary.  

As a second main result, we show that the space $\Conf(k,2r)$ is quasi-isomorphic to a Grassmannian $\tGr(k,rk)$. That is, the cluster variables and clusters in these two spaces can be identified by a pair of rational maps. The maps are not inverse, but each composite is the identity up to monomials in the frozen variables. In particular, the cluster structure on $\Conf(k,2r)$ inherits a braid group action. Note that $\tGr(k,rk)$ and $\Conf(k,2r)$ are not birationally isomorphic -- they have different dimensions -- but they become isomorphic after taking products with complex tori of appropriate size (this follows from our quasi-isomorphism and \cite[Proposition 5.11]{LamSpeyer}). 

\vspace{.3cm}

When $k=2$, the cluster combinatorics for $\Gr(2,n)$ -- more generally, for the space of twisted decorated $\SL_2$-local systems on any surface \cite{CATSI} --  has an elegant description (cf.~Example~\ref{eg:Pluckers}). Fomin and Pylyavskyy have conjectured a combinatorial description of the cluster combinatorics for $\Gr(3,n)$ in terms of Kuperberg's \emph{non-elliptic web basis}~\cite{Kuperberg}. The combinatorics involved is much more complicated than in the $\SL_2$ case. When $n \leq 8$, there are only finitely many cluster variables for $\Gr(3,n)$, and verifying the Fomin-Pylyavskyy description is a finite check. The first nontrivial case is $n=9$, which is of \emph{infinite type, but finite mutation type}:  it has infinitely many cluster variables, but only finitely many quivers. 

As a third main result, we prove (most of) the Fomin-Pylyavskyy conjectures in the case of~$\Gr(3,9)$, and give a presentation for the cluster modular group. The key point is that the braid group action preserves all relevant notions from web combinatorics, after which we verify the conjecture with a SAGE program \cite{SAGEPackage,SAGE} available as an ancillary file \cite{ArxivFraser}. We obtain similar results for $\Gr(4,8)$ (the other $\Gr(k,n)$ of finite mutation type). 

\vspace{.2cm}

Organization: The following sections contain standard background material. Section~\ref{secn:clusteralgebras} introduces cluster algebras; Section~\ref{secn:QHs} reviews quasi-homomorphisms and the cluster modular group;  Section~\ref{secn:Grassmannians} reviews Grassmannian cluster algebras; Section~\ref{secn:BraidGroup} reviews some braid group facts and definitions; Section~\ref{secn:FockGoncharovSpaces} reviews the cluster structure on ~$\Conf(k,r)$. 

The following sections contain results: Section~\ref{secn:BraidGroupAction} defines the maps~$\sigma_i$ and states and proves our main theorem (Theorem~\ref{thm:braidgroupacts}). Section~\ref{secn:themaps} describes a quasi-isomorphism of $\Gr(k,rk)$ with $\Conf(k,2r)$. Section~\ref{secn:groups} summarizes what is known about the cluster modular group of $\Gr(k,n)$ and makes a conjecture describing them. Section \ref{secn:Webs} narrows our focus to $\Gr(3,9)$ and $\Gr(4,8)$. We review web combinatorics and state the Fomin-Pylyavskyy conjectures for $\Gr(3,n)$. Theorem~\ref{thm:GrThreeNine} proves most of these conjectures for $\Gr(3,9)$, and \ref{thm:GrFourEight} proves an analogue for $\Gr(4,8)$. Theorems \ref{thm:CMG3presentation} and \ref{thm:CMG4presentation} give a presentation of their cluster modular groups. The proofs for Section~\ref{secn:Webs} are in Section~\ref{secn:proofs}.

\vspace{.2cm}
       
\begin{center}
\textsc{Prior work}
\end{center}
A connection between braid groups and cluster modular groups is inspired by \cite{FGXVarieties,ModuliSpaces}. 
For \emph{any} $G$, and for $S$ a \emph{punctured} disk with $2r$ marked points on its boundary, they stated the existence \cite{ModuliSpaces} of a homomorphism from the $G$-braid group to the cluster modular group, but the details have not been published. This is similar to our Theorem \ref{thm:braidgroupacts}, but concerns a different class of cluster algebras because our disk is unpunctured. 

Our results in $\Gr(3,9)$ and $\Gr(4,8)$ also have antecedents. Barot, Geiss, and Jasso \cite{TubularII}, as well as Felikson, Shapiro, Thomas, and Tumarkin \cite{GrowthRate}, gave Ping-Pong lemma arguments for a relation between these cluster modular groups and $\PSL_2(\bbz)$. We also mention that a certain map $(\Psi \circ P \circ \Phi) \colon \Gr(3,9) \to \Gr(3,9)$ is a generator in our description of the cluster modular group. This map was also discovered by Morier-Genoud, Ovsienko, and Tabachnikov \cite[Section 4.6]{MgOT}, who thought of it as a symmetry of the space of convex 9-gons in $\bbr \bbp^2$.

\vspace{.2cm}

\begin{center}
\textsc{Acknowledgements}
\end{center}
The early stages of this work appeared in the extended abstract \cite{FPSAC}. 
Thanks to: Andrew Nietzke, Pavel Tumarkin, Pavlo Pylyavskyy for suggesting the Grassmann-Cayley algebra, Kontstanze Rietsch for suggesting the braid group, Dylan Thurston for suggesting the correct construction when $k$ does not divide $n$, and Ian Le for many ideas and conversations. The SAGE program was begun with David Speyer during SAGE days 64.5, and I thank him and the organizers. Most of all, I thank my Ph.D. advisor Sergey Fomin for his wisdom and encouragement, and for suggesting this line of inquiry.

\section{Cluster algebras}\label{secn:clusteralgebras}
\begin{defn}
A \emph{quiver} is a finite directed graph $Q$ on the vertex set $[1,n]$, without loops or directed $2$-cyles. An \emph{extended quiver} is a finite directed graph on the vertex set $[1,n+m]$ without loops or $2$-cycles. The last $m$ vertices are called \emph{frozen vertices}, and the first $n$ vertices are \emph{mutable vertices}. We disallow arrows between frozen vertices. The integer $n$ is called the \emph{rank} of the extended quiver. 
\end{defn}  

We denote extended quivers by $\tQ$ and denote by $Q$ their underlying \emph{mutable subquivers} obtained by restricting to the mutable vertices.  

\begin{defn}
Let $\tilde{Q}$ be an extended quiver and $k \in [1,n]$  a mutable vertex. The operation of \emph{quiver mutation in direction $k$} replaces $\tilde{Q}$ by a new extended quiver $\tilde{Q}' = \mu_k(\tilde{Q})$. The quiver $\mu_k(\tilde{Q})$ is obtained from $\tilde{Q}$ in three steps: 
\begin{enumerate}
\item For each directed path $i \to k \to j$ of length two through $k$ in $\tQ$, add an arrow $i \to j$ (do not perform this step if both of $i$ and $j$ are frozen).
\item Reverse the direction of all arrows incident to vertex $k$.
\item Remove any oriented $2$-cycles created in performing steps~1 and~2. 
\end{enumerate}
\end{defn}

Mutation commutes with the operation of restricting to mutable subquivers. Furthermore, $\mu_k^2(\tilde{Q}) = \tilde{Q}$ for any mutable vertex~$k$.

\begin{defn}[Seed]
Let $\mcf$ be a field isomorphic to a field of rational functions in $n+m$ variables. A \emph{seed in $\mcf$} is a pair $(\tQ,\mathbf{x})$ where $\tQ$ is an extended quiver on $n+m$ variables, and $\mathbf{x} = (x_1,\dots,x_n; x_{n+1},\dots,x_{n+m})$ is a transcendence basis for $\mcf$. The elements $x_{n+1},\dots,x_{n+m}$ are called \emph{frozen variables}. The set $\{x_1,\dots,x_n\}$ is a \emph{cluster}, and the set $\{x_1,\dots,x_{n+m}\}$ is an \emph{extended cluster}. 
\end{defn}

\begin{defn}[Seed mutation, exchange ratio]
Let $\Sigma = (\tQ,\mathbf{x})$ be a seed and $k$ be a mutable vertex. The operation of \emph{seed mutation in direction $k$} replaces $\Sigma$ by a seed $\Sigma' = \mu_k(\Sigma) = (\mu_k(\tQ),\mathbf{x}')$. The new extended cluster $\mathbf{x}'$ satisfies $\mathbf{x}' = \mathbf{x} \backslash \{x_k\} \cup \{x_k'\}$.  The new cluster variable $x_k'$ is defined by an \emph{exchange relation}: 
\begin{equation}\label{eq:exchangereln}
x'_k x_k = \prod_{i \in [1,n+m]} x_i^{\text{number of edges $i \to k$}}+ \prod_{i \in [1,n+m]} x_i^{\text{number of edges $k \to i$}},
\end{equation}
where the numbers in the right hand side of \eqref{eq:exchangereln} refer to edges in $\tQ$. 
We define also the \emph{exchange ratio} to be the Laurent monomial
\begin{equation}\label{}
\hat{y}_\Sigma(x_k) = \frac{\prod_{i \in [1,n+m]} x_i^{\text{number of edges $i \to k$}}}{\prod_{i \in [1,n+m]} x_i^{\text{number of edges $k \to i$}}},
\end{equation}
i.e. the ratio of the terms in the right hand side of the exchange relation \eqref{eq:exchangereln}. 
\end{defn}

\begin{defn}[Cluster algebra]\label{defn:CA}
Let $\Sigma$ be a seed in $\mcf$. The \emph{seed pattern} $\mce$ determined by $\Sigma$ is the collection of seeds which can be obtained from $\Sigma$ by performing an arbitrary sequence of mutations from $\Sigma$. The \emph{cluster algebra} associated to $\mce$ is the $\bbc$-algebra generated by the frozen variables, the inverses of the frozen variables, and all of the cluster variables arising in the seeds of $\mce$. 
\end{defn}

Sometimes, rather than making Definition~\ref{defn:CA}, one instead defines the cluster algebra as the algebra generated by the frozen and cluster variables only (i.e., inverses of frozen variables are not taken as generators). For combinatorial purposes, either of these conventions is equally good. A \emph{cluster monomial} is an element that can be expressed as a monomial in any extended cluster (sometimes one allows frozen variables in the denominator, but we do not take this convention). 

Seed mutation $\mu_k$ is an involution. As such, the seed pattern~$\mce$ is determined by any of its seeds, and thus by a choice of an extended quiver $\tQ$. By a deep result \cite{IKLP}, the combinatorics of the seed pattern is in fact determined by any mutable subquiver $Q \subset \tQ$. Precisely, if $Q \subset \tQ$ and $Q \subset \tQ'$ are two different extensions of a mutable quiver to an extended quiver, determining seeds $\Sigma(\tQ)$ and $\Sigma(\tQ')$, then a sequence of mutations $\vec{\mu}$ satisfies $\vec{\mu}(\Sigma(\tQ)) = \Sigma(\tQ)$ if and only if $\vec{\mu}(\Sigma(\tQ')) = \Sigma(\tQ')$. A seed pattern (and its cluster algebra) is \emph{finite type} if the seed pattern consists of only finitely many seeds. Less restrictively,it is of \emph{finite mutation type} if the seed pattern contains only finitely many (isomorphism classes of) mutable subquivers.

\section{Symmetries of cluster algebras}\label{secn:QHs}
A quasi-homomorphism is a notion of map between cluster algebras, introduced and systematically studied in \cite{Fraser}. We give a streamlined account of the definition here. 

For a cluster algebra $\mca$ we denote by $\bbp$ the group of Laurent monomials in the frozen variables for $\mca$. For elements $x,y \in \mca$, we say that $x$ \emph{is proportional to} $y$, writing $x \propto y$, if $x = M y$ for some Laurent monomial $M \in \bbp$. Likewise, let $\mca$ and $\ov{\mca}$ be a pair of cluster algebras with respective groups $\bbp$ and $\ov{\bbp}$. If $f_1 \colon \mca \to \ov{\mca}$ and $f_2 \colon \mca \to \ov{\mca} $ are algebra homomorphisms satisfying  $f_1(\bbp) \subset \ov{\bbp}$ and $f_2(\bbp) \subset \ov{\bbp}$, then we say that $f_1$ \emph{is proportional to} $f_2$, and write $f_1 \propto f_2$, if $f_1(x) \propto f_2(x)$ for every cluster variable $x \in \mca$. 

\begin{defn}[Quasi-homomorphism]\label{defn:QH}
Consider a pair of cluster algebras $\mca$ and $\ov{\mca}$, both of the same rank $n$, and with respective groups $\bbp$ and $\ov{\bbp}$. Then an algebra map $f \colon \mca \to \ov{\mca}$ that satisfies $f(\bbp) \subset \ov{\bbp}$ is a \emph{quasi-homomorphism} from $\mca$ to $\ov{\mca}$ if there are seeds $\Sigma = (\tilde{Q},\mathbf{x})$ for $\mca$ and $\ov{\Sigma} = (\ov{\tilde{Q}},\ov{\mathbf{x}})$ for $\ov{\mca}$, and a sign $\eps \in \{1,-1\}$, for which
\begin{equation}\label{eq:QH}f(x_i) \propto \ov{x}_i \text{ and } f(\hat{y}_\Sigma(x_i)) = \hat{y}_{\ov{\Sigma}}(\ov{x}_i)^{\eps}, \text{ for } i=1,\dots,n.\end{equation}
\end{defn}

If the conditions in this definition hold for \emph{some} pair of seeds $\Sigma,\ov{\Sigma}$, then it holds for all seeds in the respective seed patterns \cite{Fraser}. That is, for every seed $\Sigma$ in $\mca$, one can define a seed $\ov{\Sigma}$ in $\ov{\mca}$ satisfying \eqref{eq:QH}, so that the map $\Sigma \mapsto \ov{\Sigma}$ commutes with mutation. 

A \emph{quasi-isomorphism} of two cluster algebras $\mca$ and $\ov{\mca}$ is a pair of quasi-homomorphisms $f \colon \mca \to \ov{\mca}$ and $g \colon \ov{\mca} \to \mca$ such that the composite $g \circ f $ is proportional to the identity map on $\mca$. A \emph{quasi-automorphism} is a quasi-isomorphism of a cluster algebra with itself.

In general, a quasi-homomorphism $f$ does not necessarily have a ``quasi-inverse'' $g$. Furthermore, if such a $g$ \emph{does} exist, it is not uniquely defined. Rather, there is a family of possible~$g$, all proportional to each other, and related to each other by certain ``rescalings'' \cite[Remark 4.7]{Fraser} of the cluster variables in $\mca$ by elements of~$\bbp$. On the other hand, a quasi-homomorphism from a cluster algebra to itself always has a quasi-inverse \cite[Lemma 6.4]{Fraser} (and in fact, a family of quasi-inverses related to each other by rescalings). Thus, we typically suppress the choice of quasi-inverse when thinking about quasi-automorphisms.

Quasi-automorphisms are closely related to other notions of symmetry in cluster theory (cf.~\cite{ADS,ASS,TubularII, FG}). Of these, we will focus on the \emph{cluster modular group}. To define it, we use a result \cite[Theorem 4]{GSV} of Gekhtman, Shapiro, and Vainshtein: in any cluster algebra, the mutable subquiver $Q(\Sigma)$ in a seed $\Sigma$ is determined by its cluster $\textbf{x}$. Thus, we can write $Q = Q(\textbf{x})$. For a quiver $Q$, we let $Q^{\text{opp}}$ denote the \emph{opposite quiver} in which the orientations of edges is reversed. 

\begin{defn}[Cluster modular group]\label{defn:CMG}
Let $\mca$ be a cluster algebra. Let $\pi$ be a permutation of the cluster variables in $\mca$ that preserves clusters and commutes with mutation. The \emph{cluster modular group} $\CMG = \CMG(\mca)$ consists of such permutations $\pi$ for which one of the following holds: for every cluster ~$\mathbf{x}$ in~$\mca$, the induced map on mutable subquivers $Q(\mathbf{x}) \to Q(\pi(\mathbf{x}))$ is a quiver isomorphism, or; this is true for the induced map $Q(\mathbf{x})^{\text{opp}} \to Q(\pi(\mathbf{x}))$. The element $\pi$ is called \emph{orientation-preserving} or \emph{orientation-reversing} respectively. We denote by $\CMG^+$ the subgroup of orientation-preserving elements. 
\end{defn}

The group $\CMG$ only depends on the mutation pattern of mutable subquivers $Q$ (not on the pattern of extended quivers $\tQ$). Every quasi-automorphism determines an element of the cluster modular group via its underlying map on cluster variables $x \mapsto \ov{x}$; this element is in $\CMG^+$ or in $\CMG \setminus \CMG^+$ according to whether $\eps = 1$ or $=-1$ in Definition~\ref{defn:QH}. Provided a technical ``row span'' condition is satisfied \cite[Corollary 4.5]{Fraser}, \emph{every} element $g \in \CMG$ comes from some quasi-automorphism, and in fact each $g$ corresponds to a family of quasi-automorphisms, all of which are proportional to each other. This technical condition is satisfied for the cluster structures on $\Gr(k,n)$ and $\Conf(k,r)$. Thus, we can identify $\CMG$ with the group of proportionality classes of quasi-automorphisms. We summarize cluster algebras for which $\CMG$ has been computed in Remarks~\ref{rmk:finitetype} and \ref{rmk:fmtype}.

The groups $\CMG^+$ and $\CMG$ are related as follows \cite[Theorem 2.11]{ASS}. One has $\CMG^+ = \CMG$ provided for some (equivalently, any) mutable subquiver $Q$ in the seed pattern, there is a sequence of quiver mutations $Q \to Q^{\text{opp}}$. Otherwise, $\CMG^+$ is an index two subgroup of $\CMG$. For the cluster algebras $\Gr(k,n)$ and $\FG(k,r)$ that we consider in this paper, $Q$ is mutation equivalent to $Q^{\text{opp}}$ (this follows from the theory of weak separation, cf.~Section~\ref{secn:Grassmannians}), and therefore $\CMG^+$ has index two in $\CMG$.

Finally we record here a criterion for verifying that a map between cluster algebras is a quasi-isomorphism. By a \emph{nerve} $\mcn$ for $\mca$ we will mean a finite subset of the clusters in $\mca$, such that the clusters in $\mcn$ are pairwise connected to each other by sequences of mutations that stay in $\mcn$, and such that the intersection over all clusters in $\mcn$ is empty. The empty intersection hypothesis is the same as requiring that every cluster variable that shows up in a cluster in $\mcn$ is mutated at least once on the nerve. The simplest example of a nerve (but not the one we will use in the sequel) is a cluster together with its $n$ neighboring clusters. 

\begin{lem}[Constructing a quasi-isomorphism {\cite[Lemma 3.7 and Proposition 5.2]{Fraser}} 
]\label{lem:QHnerves} Let $\mcn$ be a nerve for $\mca$. Let $f \colon \mca \to \ov{\mca}$ be an algebra map satisfying $f(\bbp) \subset \ov{\bbp}$. Suppose for each cluster variable $x$ on $\mcn$, there is a cluster variable $\ov{x} \in \ov{\mca}$ such that $f(x) \propto \ov{x}$. Suppose this map sends clusters on $\mcn$ to clusters in $\ov{\mca}$ in a way that is compatible with mutation. 
Then~$f$ is a quasi-homomorphism. If an algebra map $g \colon \ov{\mca} \to \mca$ satisfies $g(\ov{\bbp}) \subset \bbp$ and $g(f(x)) \propto x$ for all cluster variables $x \in \mcn$, then $g$ is a quasi-inverse quasi-homomorphism to~$f$.  
\end{lem}

\section{Grassmannian cluster algebras}\label{secn:Grassmannians}
We introduce four spaces $(V^n)^\circ, \tGr(k,n), \tGr^\circ(k,n)$, and $\tGr^\circ(k,n) / (\bbc^*)^n$, each of which is a frequently encountered variant of the Grassmannian $\Gr(k,n)$ of $k$-planes in $\bbc^n$. We assume throughout that $k \geq 2$. 

Throughout this paper, we fix a $k$-dimensional complex vector space $V$, and a volume form $\omega^* \colon \bigwedge^k(V) \to \bbc$ with dual volume form $\omega \in \bigwedge^k(V)$. We denote by $\bigwedge(V)$ the exterior algebra for $V$, and we always denote the exterior product map $\bigwedge^a(V) \otimes \bigwedge^b(V) \to \bigwedge^{a+b}(V)$ by multiplication, writing $v_1 \cdots v_a$ rather than $v_1 \wedge \cdots \wedge v_a$. Every simple tensor $v = v_1 \cdots v_a \in \bigwedge^a(V)$ defines an $a$-dimensional subspace 
\begin{equation}\label{eq:tensortosubspace}
\ov{v} = \{w \in V \colon w  v = 0\},
\end{equation}
and the subspace $\ov{v}$ characterizes the tensor $v$ up to a scalar multiple. Geometrically, the tensor $v$ is the data of the subspace $\ov{v}$ together with a choice of volume form on $\ov{v}$.

We denote by $\tGr(k,n) \subset \bbc^{\binom n k}$ the affine cone over $\Gr(k,n) \subset \bbp^{\binom n k -1}$ in its Pl\"ucker embedding. Thus $\tGr(k,n)$ is the affine subvariety of $\bbc^{\binom n k}$ consisting of points whose coordinates satisfy the \emph{Pl\"ucker relations}. 

A \emph{configuration of $n$ vectors in~$V$} is a point in the space $\SL(V) \backslash V^n$, i.e. an $n$-tuple of vectors in $V$, with these vectors considered up to simultaneous $\SL(V)$ action. Such a configuration is \emph{weakly generic} if $v_1,\dots,v_n$ span $V$. The space of weakly generic configurations can be identified with $\tGr(k,n)\backslash \{0\} \subset \bbc^{\binom n k}$. The coordinate ring $\bbc[\tGr(k,n)]$ is generated by \emph{Pl\"ucker coordinates} $\Delta_I$, as $I$ ranges over subsets  $\{i_1 < i_2 < \dots < i_k\} \subset \{1,\dots,n\}$. The function $\Delta_I$ evaluates on a configuration
$(v_1,\dots,v_n)$ by 
\begin{equation}\label{eq:PluckeronConfig}
\Delta_I ((v_1,\dots,v_n))  =  \omega^* (v_{i_1} \cdots v_{i_k}) = \det (v_{i_1} | \cdots | v_{i_k}). 
\end{equation}

An $n$-tuple of vectors $(v_1,\dots,v_n)$ is \emph{consecutively generic} if every cyclically consecutive $k$-tuple of vectors is linearly independent, i.e. $\omega^*(v_{i+1} \cdots v_{i+k}) \neq 0$ for $i=1,\dots,n$ where we treat indices modulo $n$. We denote by $(V^n)^\circ \subset V^n$ the quasi-affine variety consisting of consecutively generic $n$-tuples. 

We also consider the quasi-affine variety $\tGr(k,n)^\circ \subset \tGr(k,n)$ defined by the non-vanishing of the cyclically consecutive Pl\"ucker coordinates $\Delta_{i+1,\dot,i+k}$. This space is the affine cone over the \emph{open positroid stratum} $ \Gr(k,n)^\circ \subset \Gr(k,n)$. Points in $\tGr(k,n)^\circ$ are identified with points in $\SL(V) \backslash (V^n)^\circ$.  The coordinate ring $\bbc[\tGr^\circ(k,n)]$ is the localization of $\bbc[\tGr(k,n)]$ at the cyclically consecutive Pl\"ucker coordinates.  

The space of $n$-tuples $V^n$ is endowed with a right action by an algebraic torus $\bbt = (\bbc^*)^n$ rescaling each of the vectors. This action commutes with the $\SL(V)$ action and preserves consecutive genericity, yielding a $\bbt$-action on $(V^n)^\circ,\tGr(k,n)$, and $\tGr^\circ(k,n)$. The space $\tGr^\circ(k,n) / \bbt$ bears a ``cluster-$\mcx$ structure'' (without any frozen variables). 

J.~Scott introduced a cluster algebra structure on $\bbc[\tGr(k,n)]$ \cite{Scott}. The frozen variables are the cyclically solid Pl\"ucker coordinates. Thus, Scott's recipe can also be thought of as a cluster algebra structure on $\bbc[\tGr^\circ(k,n)]$ in which the frozen variables are inverted. From a combinatorial perspective, the distinction between these two cluster algebras is unimportant, and we freely translate between them. We introduce certain seeds in the cluster structure on $\tGr(k,n)$ in the next two sections. The combinatorial details of these seeds are used in our proof that the braid group acts by mutations.

The torus action $\Gr(k,n) \curvearrowleft \bbt$ induces a $\bbz^n$-grading on the algebras $\bbc[\tGr(k,n)]$ and $\bbc[\tGr^\circ(k,n)]$. Our proofs make use of the well-known fact that every cluster variable $x \in \bbc[\tGr(k,n)]$ is homogeneous with respect to this grading, and the closely related fact that every exchange ratio $\hat{y}_\Sigma(x)$ is invariant under the $\bbt$ action. 

\begin{example}\label{eg:Pluckers}
When $k=2$, the cluster algebra $\bbc[\tGr(2,n)]$ has only finitely many clusters and cluster variables. The cluster variables are exactly the (non-frozen) Pl\"ucker coordinates. Let $D$ be an $n$-gon with its vertices numbered $1,\dots,n$ in clockwise order. We can index the Pl\"ucker coordinate $\Delta_{ij}$ by the straight line $(ij)$ connecting vertices $i$ and $j$ in~$D$.  The frozen variables are the sides of~$D$. A pair of Pl\"ucker coordinates are in a cluster if any only if the corresponding arcs $(ij)$ are pairwise noncrossing. The clusters are in bijection with triangulations of the $n$-gon $D$. 

For $k \geq 3$, the combinatorial description of clusters and cluster variables in $\bbc[\tGr(k,n)]$ is much more complicated. Besides a few small examples, these cluster algebras have infinitely many clusters and cluster variables.  We record here the Grassmannians of finite type, and of finite mutation type: 
\begin{itemize}
\item $\bbc[\tGr(2,n+3)]$ has finite Dynkin type $A_{n}$.
\item  $\bbc[\tGr(3,6)], \bbc[\tGr(3,7)]$ and $\bbc[\tGr(3,8)]$ have finite Dynkin types $D_4, E_6, E_8$ respectively.
\item $\bbc[\tGr(3,9)]$ and $\bbc[\tGr(4,8)]$ are of infinite type but finite mutation type. 
\item The cluster algebras $\bbc[\tGr(k,n)]$ and $\bbc[\tGr(n-k,n)]$ can be identified with each other 
by the complementation map on Pl\"ucker coordinates. Thus, $\bbc[\tGr(n+1,n+3)]$ has cluster type $A_n$ and so on.
\end{itemize}

When $k=3$, Fomin and Pylyavskyy \cite{tensors} have proposed a description of the cluster combinatorics in terms of Kuperberg's basis of \emph{non-elliptic webs}. We review this description in Section~\ref{secn:Webs}, and give strong evidence for its correctness for $\bbc[\tGr(3,9)]$ in Theorem~\ref{thm:GrThreeNine}. 
\end{example}

\subsection{Seeds from weakly separated collections} We let $\binom {[n]} k$ denote the set of subsets of $\{1,\dots,n\}$ of size $k$. We review now the construction of certain combinatorially defined seeds for $\tGr(k,n)$. A pair of subsets $I,J \in \binom {[n]} k$ is \emph{weakly separated} if the sets $I\setminus J$ and $J \setminus I$ are cyclically disjoint. In other words, there is no cyclic interval $a < b < c < d$ where $a,c \in I\setminus J$ and $b,d \in J\setminus I$. A {\it weakly separated collection} $\mcc \subset \binom{[n]} k$ is a collection whose members are pairwise weakly separated. Such a collection is \emph{maximal} if there is no larger weakly separated collection containing it; it is known \cite{OPS} that all maximal weakly separated collections have $\dim \tGr^\circ(k,n) = k(n-k)+1$ elements.

By a recipe initially due to Postnikov \cite{Postnikov}, every maximal weakly separated collection $\mcc$ determines a seed $\Sigma(\mcc) = (\tQ(\mcc),\textbf{x}(\mcc))$ in the cluster structure on $\tGr(k,n)$. Its extended cluster $\textbf{x}(\mcc)$ is the set of Pl\"ucker coordinates $\{\Delta_I \colon I \in \mcc\}$. The cyclically solid Pl\"ucker coordinates are the frozen variables. We describe $\tQ$ after the following definition. 

\begin{defn}
Given $S \in \binom {[n]}{k+1}$, the \emph{clique} determined by $S$ is the collection $\mcw(S) = \{I \in \mcc \colon \, I \subset S\}$. This clique is \emph{nontrivial} if its cardinality is at least three.  
\end{defn}

Clearly, if $k$-subsets $I$ and $J$ are in the same clique, then $|I \cap J| = k-1$. 
 
Now we describe the extended quiver $\tQ(\mcc)$ \cite{OPS}. First, $\tQ(\mcc)$ has not multiple arrows. Second, if $S$ is a nontrivial clique, and if $I$ is in $\mcw(S)$, then $S \setminus I$ is a singleton, and we can use these singleton sets to cyclically order the elements of $\mcw(S)$. Then if $I,J \in \mcc$ are two subsets not both of which are frozen, then there is an arrow $\Delta_I \to \Delta_J$ in $Q(\mcc)$ exactly when the clique $\mcw(I \cup J)$ is nontrivial and furthermore $I$ and $J$ are cyclically adjacent in $\mcw(I \cup J)$ (and with $I$ preceding $J$).  

For a maximal weakly separated collection $\mcc$, let $I \in \mcc$ be any mutable Pl\"ucker coordinate that has exactly 4 neighbors in $\tQ(\mcc)$. Then mutation at $\Delta_I$ produces a seed $\Sigma(\mcc')$ where $\mcc'$ is again a maximal weakly separated collection. The resulting operation on weakly separated collections is known as a \emph{square move}. By \cite{Postnikov}, any two maximal weakly separated collections are connected by a sequence of square moves. Thus the seeds $\Sigma(\mcc)$ we have just described lie in a common seed pattern.

\subsection{Convenient choice of initial seed} One has the following standard choice of maximal weakly separated collection: 
\begin{equation}\label{eq:stdseed}
\mcc_{\Le} = \{[1,a] \cup [b+1,b+k-a] \colon \, 0 \leq a \leq k , \, a+1 \leq b \leq n-k+a \} \subset \binom{[n]}k. 
\end{equation}
The collection $\mcc_{\Le}$ is the ``$\Le$ diagram'' cluster for the top positroid cell. 

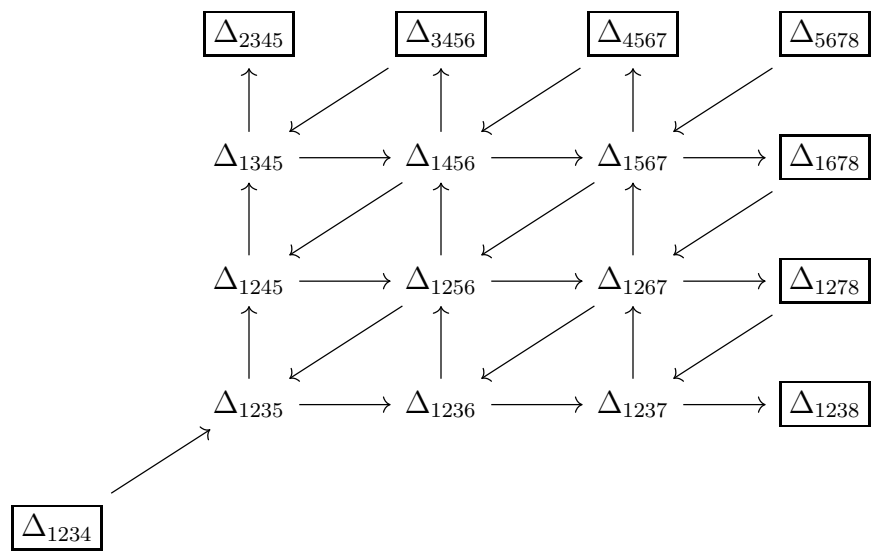
\begin{figure}
\begin{tikzcd}
 & \boxed{\Delta_{2345}} & \boxed{\Delta_{3456}} \arrow[ld] & \boxed{\Delta_{4567}} \arrow[ld] & \boxed{\Delta_{5678}} \arrow[ld] \\
 & \Delta_{1345} \arrow[r] \arrow[u] & \Delta_{1456} \arrow[r] \arrow[u] \arrow[ld] & \Delta_{1567} \arrow[r] \arrow[u] \arrow[ld] & \boxed{\Delta_{1678}} \arrow[ld] \\
 & \Delta_{1245} \arrow[r] \arrow[u] & \Delta_{1256} \arrow[r] \arrow[u] \arrow[ld] & \Delta_{1267} \arrow[r] \arrow[u] \arrow[ld] & \boxed{\Delta_{1278}} \arrow[ld] \\
 & \Delta_{1235} \arrow[r] \arrow[u] & \Delta_{1236} \arrow[r] \arrow[u] & \Delta_{1237} \arrow[r] \arrow[u] & \boxed{\Delta_{1238}} \\
\boxed{\Delta_{1234}} \arrow[ru] &  &  &  & 
\end{tikzcd}
\caption{The $\Le$ diagram seed for $\bbc[\tGr(4,8)]$. Boxed variables are frozen. \label{fig:stdseed48}
}
\end{figure}

We refer to the frozen variable $[1,k]$ as the ``extra frozen variable,'' occurring when $a=k$ in \eqref{eq:stdseed}. Besides the extra frozen variable, we associate the Pl\"ucker coordinate indexed by $a$ and $b$ in \eqref{eq:stdseed} with the lattice point $(b-a,k-a) \in \bbz^2$. This identifies $\mcc_{\Le} \setminus \{[1,k]\}$ with the entries in a $k \times n-k$ rectangular array, cf.~Figure~\ref{fig:stdseed48}. We place the extra frozen variable at the origin in $\bbz^2$. The extended quiver $\tQ(\mcc_{\Le})$ consists of all eastward, northward, and southwest arrows in the rectangular array, as well as a northeast arrow from the extra frozen variable, cf.~Figure~\ref{fig:stdseed48}.

\section{Braid groups}\label{secn:BraidGroup}
\begin{defn}\label{defn:braidgroup}
The \emph{braid group on $k$ strands} is the group $B_k$ with generators $\sigma_1,\dots,\sigma_{k-1}$, subject to the relations 
\begin{equation}\label{eq:braidrelns}
\sigma_i \sigma_{i+1} \sigma_i = \sigma_{i+1} \sigma_{i} \sigma_{i+1} \text{ for $i = 1,\dots,k-1$}, \quad \sigma_i \sigma_{j}  = \sigma_{j} \sigma_{i} \text{ for $|i-j| \geq 2$}.
\end{equation}
\end{defn}

Imagine connecting the points $(1,0),\dots,(k,0)$ in $\bbr^2$ to the points $(1,1),\dots,(k,1)$ by straight lines known as \emph{strands}. The Artin generator $\sigma_i$ acts on such a picture by crossing the $i$th strand over the $i+1$th strand (we indicate the under-crossing using dashed lines). A \emph{braid diagram} is a picture that can be obtained from the initial straight lines picture by a finite sequence of applications of the Artin generators. Such a diagram should be thought of as a two-dimensional projection of $k$ strands in $\bbr^3$. Two braid diagrams are identified if the corresponding braids are related by an ambient isotopy of $\bbr^3$.

We recall the following well known algebraic properties of braid groups \cite{Birman,KasselTuraev}.  
\begin{itemize}
\item The element $\Delta = \sigma_1(\sigma_2 \sigma_1)(\sigma_3 \sigma_2 \sigma_1) \cdots  (\sigma_{k-1} \cdots \sigma_1) \in B_k$ is known as the half-twist. There is also an element $\delta = \sigma_{k-1} \cdots \sigma_1$ which we call the $\frac{1}{k}$-twist. These elements are related by $\Delta^2 = \delta^k$.  
\item The center $Z(B_k)$ is infinite cyclic, and is generated by a full-twist $\Delta^2 = \delta^k$. 
\item Conjugation by $\delta$ has order $k$ in $\Aut(B_k)$. Furthermore, we have $\delta^{-1} \sigma_{i} \delta= \sigma_{i+1}$ for $1 \leq i \leq k-2$. Setting $\sigma_k = \delta^{-1} \sigma_{k-1} \delta \in B_k$, we obtain elements $\sigma_1,\dots,\sigma_k \in B_k$ that are cyclically permuted under conjugation by $\delta$.
\end{itemize}

%

The group $B_1$ is the trivial group, and $B_2 = \bbz$. The special case of $B_3$ will be especially important to us. It is classically known that $B_3$ modulo its center is a free product of cyclic groups:  
\begin{equation}\label{eq:PSL2Z}
B_3 / Z(B_3) \cong \bbz / 2 \bbz \ast \bbz / 3 \bbz \cong \PSL_2(\bbz). 
\end{equation} 
The elements $\Delta,\delta \in B_3 / Z(B_3)$ generate the two factors of this free product.

\subsection{Extended affine braid groups}
Our reference is \cite{ExtendedBraid}. 
\begin{defn}\label{defn:extdbraidgroup}
The \emph{extended affine braid group on $d$ strands} is the group $\hat{B}_{\hat{A}_{d-1}}$ with generators $\sigma_1,\dots,\sigma_{d-1},\rho$, presented by the relations 
\begin{align}\label{eq:extdbraidrelns}
\sigma_i \sigma_{i+1} \sigma_i = \sigma_{i+1} \sigma_{i} \sigma_{i+1} &\text{ for $i = 1,\dots,d-1$, and } \quad \sigma_i \sigma_{j}  = \sigma_{j} \sigma_{i} \text{ for $|i-j| \geq 2$} \\
\rho^{-1} \sigma_i \rho = \sigma_{i+1} &\text{ for $i = 1,\dots,d-2$, and} \quad \rho^{-1} \sigma_{d-1} \rho = \rho \sigma_{1} \rho^{-1}.
\end{align}
\end{defn}
It also goes by the name of \emph{circular braid group} or \emph{annular braid group}.
As in the ordinary braid group, one can define an element $\sigma_d = \rho^{-1} \sigma_{d-1} \rho = \rho \sigma_{1} \rho^{-1}$, to obtain elements $\sigma_1,\dots,\sigma_d \in \hat{B}_{\hat{A}_{d-1}}$ that are cyclically permuted under conjugation by~$\rho$.
The center of $\hat{B}_{\hat{A}_{d-1}}$ is infinite cyclic, generated by $\rho^d$. When $d=1$, we interpret $\hat{B}_{\hat{A}_{d-1}} = \la \rho \ra \cong \bbz $.

The following centralizer interpretation of the extended affine braid group is very natural from our point of view. Let $d$ be a divisor of a positive integer $k$. One can define a group homomorphism $\iota_{d,k} \colon \hat{B}_{\hat{A}_{d-1}} \to B_k$ on the generators by
\begin{align}\label{eq:iotadefn}
\sigma_i \in \hat{B}_{\hat{A}_{d-1}}  & \mapsto \prod_{j=0}^{\frac{k}{d}-1} \sigma_{i+jd} \in B_k \text{ for $i = 1,\dots,d-1$} \\
\rho \in \hat{B}_{\hat{A}_{d-1}} & \mapsto \delta \in B_k,
\end{align} 
where we have used the same symbols $\sigma_i$ for the generators in both $\hat{B}_{\hat{A}_{d-1}}$ and $B_k$. Using the braid relations in $B_k$, and the properties of $\delta$ mentioned above, one checks that this is a well-defined homomorphism. 

When $d<k$ this map is injective, and in fact, $\iota_{d,k}(\hat{B}_{\hat{A}_{d-1}})$ is the centralizer subgroup $Z(\delta^d) \subset B_k$ \cite[Theorem 0.2]{BessisBraid}. On the other hand, when $d=k$, the homomorphism $\iota_{k,k} \colon \hat{B}_{\hat{A}_{k-1}} \twoheadrightarrow B_k$ is the quotient map imposing the relation $\rho  = \sigma_{k-1} \cdots \sigma_1 \in B_k$, which does not hold in the extended affine braid group. This qualitatively different behavior is also reflected in the centralizer point of view, because $Z(\delta^k) = B_k$ is a braid group (not an extended affine braid group).

\section{The braid group action}\label{secn:BraidGroupAction}
We introduce the braid group action on $\tGr^\circ(k,n)$ and then state and prove our main theorem. 

Let $\rho \colon V^n \to V^n$ denote the \emph{twisted cyclic shift}
\begin{equation}\label{eq:cyclicshift}
(v_1,\dots,v_n) \,  \overset{\rho}{\mapsto}  \,  (v_2,\dots,v_n,(-1)^{k-1} v_1). 
\end{equation}
 
We denote also by $\rho$ the induced twisted cyclic shift maps on $(V^\circ)^n$, $\tGr(k,n)$, and $\tGr^\circ(k,n)$. Notice that $\rho^n$ is the identity map on configurations because each vector is multiplied by $(-\text{Id})^{k-1}$. 

Though less important to our story, the cyclic group action on configurations can be enriched to a dihedral group action via the \emph{twisted reflection} map $\theta: V^n \to V^n$, defined as $(v_1,\dots,v_n) \mapsto (-1)^{\binom k 2}(v_n,\dots,v_1)$.

The signs $(-1)^{k-1}$ and $(-1)^{\binom k 2}$ above are necessary to ensure that the pullbacks $\rho^*$ and $\theta^*$ send Pl\"ucker coordinates to Pl\"ucker coordinates. For example in $\bbc[\tGr(4,8)]$ we have $\rho^*(\Delta_{1238}) = -\Delta_{2341} = \Delta_{1234}$. The reader is welcome to ignore these signs henceforth. 

A consequence of the combinatorial recipes from Section~\ref{secn:Grassmannians} is that both the twisted cyclic shift and the twisted reflection are \emph{cluster automorphisms} of the cluster structure on $\bbc[\tGr(k,n)]$ \cite{ADS,ASS} -- i.e. both of these maps permute the frozen variables, the cluster variables, and the clusters, respectively, in $\bbc[\tGr(k,n)]$. Furthermore, $\rho^*$ is orientation-preserving whereas $\theta^*$ is orientation-reversing. 

Now let $d$ be a divisor of $n$. A function $f \colon V^n \to V^n$ is  
\emph{$d$-periodic} if it satisfies $f \circ \rho^d = \rho^d \circ f$. We can understand a $d$-periodic map $f$ by breaking it up into $\frac{n}{d}$ \emph{windows}, each of which is a list of $d$ of the coordinate functions of $f$ written in terms of the input vectors $(v_1,\dots,v_n)$. By $d$-periodicity, the map $f$ is determined in any such window. We therefore specify $d$-periodic maps by writing down the coordinate functions for the first window, referring to this as the \emph{window notation} for the map $f$. For example, the twisted cyclic shift map has $d$-periodic window notation $[v_2,\dots,v_{d+1}]$.

The following suggests why the concept of $d$-periodicity is a relevant one. It follows from results in \cite{GoncharovShenDT,MarshScott}, as we explain in Section~\ref{secn:groups}. 
\begin{prop}\label{prop:periodicity} Let $f \colon \bbc[\tGr(k,n)] \to \bbc[\tGr(k,n)]$ be an orientation-preserving quasi-automorphism and $d = \gcd(k,n)$. Then $f \circ \rho^d \propto \rho^d \circ f$. 
\end{prop}

The notion of $d$-periodicity is stronger than that of satisfying $f \circ \rho^d \propto \rho^d \circ f$. In principle there can be quasi-automorphisms of $\bbc[\tGr^\circ(k,rk)]$ that commute with $\rho^d$ up to proportionality, but that do not come from a $d$-periodic maps $V^n \to V^n$ on vectors. 

Recall the space $(V^n)^\circ$ of consecutively generic $n$-tuples in $V$. We now make our main definition. 
\begin{defn}[Artin generator]\label{defn:sigmai}
Let $d = \gcd(k,n) > 1$. For $i=1,\dots,d-1$, define a map $d$-periodic map $\sigma_i \colon  (V^n)^\circ \to  (V^n)^\circ$ as follows. The first window of $\sigma_i$ is $[v_1,\dots,v_{i-1},v_{i+1},w_1,v_{i+2},\dots,v_d]$, with the vector $w_1$ defined uniquely by the conditions 
\begin{equation}\label{eq:sigmai}
v_iv_{i+1} = v_{i+1}w_1 \in \bigwedge^2(V) \text{ and } w_1 \in \Span\{v_{i+2},\dots,v_{i+k}\}.
\end{equation}
The $\ell$th window is defined by the same recipe by $d$-periodically augmenting indices $v_i \mapsto v_{(\ell-1)d+i}$ and defining the new vector $w_\ell$ with respect to these augmented indices. 

Define also a $d$-periodic map $\sigma_i^{-1} \colon  (V^n)^\circ \to  (V^n)^\circ$, whose first window is \\
$[v_1,\dots,v_{i-1},u_1,v_i,v_{i+2}, \dots v_d]$,
with $u_1$ defined uniquely by the conditions 
\begin{equation}\label{eq:sigmaiinverse}
v_iv_{i+1} = u_1 v_{i} \in \bigwedge^2(V) \text{ and } u_1 \in  \Span\{v_{n-k+i+1},\dots,v_n,v_1,\dots,v_{i-1}\}\end{equation}
and with the other windows defined by augmenting indices $d$-periodically. 
\end{defn}

Let us clarify why \eqref{eq:sigmai} indeed defines $w_1$ uniquely. From the equality of subspaces $\ov{v_iv_{i+1}} = \ov{v_{i+1}w_1}$ (cf.~\eqref{eq:tensortosubspace}) it follows that $w_1 \in \Span\{v_{i},v_{i+1}\} \cap \Span\{v_{i+2},\dots,v_{i+k}\}$, which is a line by consecutive genericity. 
This defines $w_1$ up to a scalar, and the scalar is fixed by the normalization $v_iv_{i+1} = v_{i+1}w_1$.

The maps \eqref{eq:sigmai} satisfy $\rho^{-1} \circ \sigma_{i} \circ \rho = \sigma_{i+1}$, where $\rho$ is the twisted cyclic shift. This matches the corresponding relation in the extended affine braid group. The maps $\sigma_i$ and $\sigma_i^{-1}$ are related by $\theta \circ \sigma_i \circ \theta  = \sigma_{d-i}^{-1}$. 

Each map $\sigma_i$ commutes with the $(\bbc^*)^n$ action up to a permutation of the factors. For example, if $\vec{v} \in (V^n)^\circ$ and $\mathbf{t} = (t_1,\dots,t_n) \in (\bbc^*)^n$, then 
$\sigma_1(\vec{v} \cdot \mathbf{t})$ equals \\
$\sigma_1(\vec{v}) \cdot (t_2,t_1,t_3,\dots,t_d,t_{d+2},\dots,t_n)$. Thus, $\sigma_i$ descends to a map on $\tGr^\circ(k,n) / (\bbc^*)^n$.

The following is our main theorem. 
\begin{thm}\label{thm:braidgroupacts}
Let $d = \gcd(k,n)>1$ and $\sigma_1,\dots,\sigma_{d-1} \colon (V^n)^\circ \to  (V^n)^\circ$ be the $d$-periodic maps just defined. 
\begin{enumerate}
\item The maps $\sigma_i$ and $(\sigma_i^{-1})$ are inverse regular automorphisms of $(V^n)^\circ$. They commute with the $\SL(V)$ action, determining inverse regular automorphisms of $\tGr^\circ(k,n)$. 
\item Together with the twisted cyclic shift $\rho$, the $\sigma_i$ satisfy the extended affine braid relations \eqref{eq:extdbraidrelns}, determining group homomorphisms $\hat{B}_{\hat{A}_{d-1}} \to \Aut((V^n)^\circ)$ and $\hat{B}_{\hat{A}_{d-1}} \to \Aut(\tGr^\circ(k,n))$.
\item Each pullback $\sigma_i^*$ is an orientation-preserving quasi-automorphism of $\bbc[\tGr^\circ(k,n)]$, determining an (anti)homomorphism $\hat{B}_{\hat{A}_{d-1}} \to \CMG^+(\tGr^\circ(k,n))$ into the cluster modular group. 
\end{enumerate}
\end{thm}

The proof of Theorem~\ref{thm:braidgroupacts} occupies the end of this section. Before getting there, we state and prove an important lemma, followed by some remarks. Recall that two maps are \emph{proportional} if they agree on cluster variables up to Laurent monomials in frozens. 

\begin{lem}\label{lem:consecutivecomp} If $k$ divides $n$, then the pullback $(\sigma_{k-1} \cdots \sigma_1)^*$ is proportional to the cyclic shift $\rho^*\in \bbc[\tGr^\circ(k,rk)]$. 
\end{lem}

Thus, in the special case that $k$ divides $n$,  the extended affine braid group action on clusters reduces to an ordinary $B_k$ action. We do not expect that the cyclic shift is proportional to a composition of Artin generators when $d <k$. 


\begin{proof} To simplify notation, we write out the proof in the case $k=5$. When we apply $\sigma_1$ to an $n$-tuples of vectors $(v_1,\dots,v_n)$, we create a new vector in each window, let us call this vector $a_\ell$ (rather than $w_\ell$ as in \eqref{eq:sigmai}). We then apply $\sigma_2$ to $\sigma_1(v_1,\dots,v_n)$, creating a vector $b_\ell$ in the $\ell$th window. Then we create a vector $c_\ell$ after applying $\sigma_3$, and $d_\ell$ after applying $\sigma_4$. In window notation: 
\begin{equation}
[v_1,\dots,v_4,v_5] \overset{\sigma_1}{\mapsto} [v_2,a_1,v_3,v_4,v_5] \overset{\sigma_2}{\mapsto} [v_2,v_3,b_1,v_4,v_5] \overset{\sigma_3}{\mapsto}  [v_2,v_3,v_4,c_1,v_5] \overset{\sigma_4}{\mapsto}  [v_2,v_3,v_4,v_5,d_1]  \label{eq:321}.
\end{equation}
The defining conditions for these vectors are the span conditions
\begin{align*}
a_1 \in & \Span \{v_3,v_4,v_5,v_6\} \\ 
b_1 \in & \Span \{v_4,v_5,v_7,a_2,\} =\Span \{v_4,v_5,v_6,v_7\} \\
c_1 \in &\Span \{v_5,v_7,v_8,b_2\} = \Span \{v_5,v_7,a_2,v_8\} =\Span \{v_5,v_6,v_7,v_8\} \\
d_1 \in & \Span \{v_7,v_8,v_9,c_2\} = \cdots  = \Span \{v_6,v_7,v_8,v_9\},
\end{align*}
as well as the normalization conditions $v_1v_2 = v_2 a_1$, $a_1v_3 = v_3 b_1$, $b_1v_4 = v_4 c_1$, and $c_1v_5 = v_5 d_1$. In performing the simplifications above, we used the equality $\Span\{v_7,a_2\} = \Span\{v_6,v_7\}$, and other analogous equalities. Since $a_1 \in \Span \{v_3,v_4,v_5,v_6\}$ and $b_1 \in \Span\{a_1,v_3\} \cap \{v_4,v_5,v_6,v_7\} $ we conclude that $b_1 \in \Span\{v_4,v_5,v_6\}$. In a similar fashion, we conclude that $c_1 \in \Span\{v_5,v_6\}$ and then that $d_1 \in \Span\{v_6\}$. So $d_1$ and $v_6$ are related by a scalar multiple. To compute this scalar multiple, we use the normalization conditions 
\begin{equation}
v_2v_3v_4v_5d_1 = v_2 v_3v_4c_1v_5 =v_2v_3b_1v_4v_5 = v_2 a_1 v_3v_4v_5 = v_1v_2 v_3v_4v_5.
\end{equation}
Thus $d_1 = \frac{\omega^*(v_1\cdots v_5)}{\omega^*(v_2\cdots v_6)}v_6$. This scalar multiple is a Laurent monomial in the frozen variables for $\tGr(k,n)$. Thus the maps $\sigma_{k-1} \cdots \sigma_1$ and $\rho$ agree as maps on $(V^n)^\circ$, up to rescaling certain vectors by Laurent monomials in the frozens. By the homogeneity of cluster variables with respect to the $(\bbc^*)^n$ rescaling the vectors, it follows that the two maps are proportional. 
\end{proof}

\begin{rmk}\label{rmk:kdvidesn} In some sense, the extended braid group action when $d<k$ can be deduced from the ordinary braid group action when $d=k$. Recall the inclusion $\iota_{d,k} \colon \hat{B}_{\hat{A}_{d-1}} \to B_k$, identifying $\hat{B}_{\hat{A}_{d-1}}$ with the centralizer $Z(\delta^d)$. Set $r = \frac{n}{d}$ so that lcm$(k,n) = rk$. We think of $V^{rk}$ as the space of $k \times rk$ matrices, and ``diagonally'' embed $V^{n} \subset V^{rk}$ as the subspace of matrices whose columns are $n$-periodic. 
A composition of the Artin generators $\sigma_1,\dots,\sigma_{k-1} \colon V^{rk} \to V^{rk}$ preserves the subspace $V^n$ precisely when it is $n$-periodic. Since each Artin generator is $k$-periodic, a composition of them is $n$-periodic precisely when it is $d$-periodic. From Lemma~\ref{lem:consecutivecomp}, the maps $\delta^d$ and $\rho^d$ differ only by rescaling certain vectors by frozen variables, so commuting with $\rho^d$ is closely related to commuting with $\delta^d$. By construction the braid $\prod_{j=0}^{\frac{k}{d}-1} \sigma_{i+jd} \in B_k$ from \eqref{eq:iotadefn}) 
determines a $d$-periodic map on $V^{rk}$, so descends to $V^n$. This is exactly the definition of the $d$-periodic Artin generator for $(V^n)^\circ$. Conversely, the braids $\prod_{j=0}^{\frac{k}{d}-1} \sigma_{i+jd} \in B_k$ generate the centralizer $Z(\delta^d)$, so one expects that these are the \emph{only} compositions of Artin generators for $V^{rk}$ that might give rise to quasi-automorphisms of $\tGr(k,n)$. 
\end{rmk}

\begin{rmk}[Renormalized Artin generators]\label{rmk:renormd}
There is a more explicit description of the vector $w_1$ in Definition~\ref{defn:sigmai}, namely 
\begin{equation}\label{eq:sigmaidiagram}
w_1 = \frac{\omega^*(v_i,v_{i+2},\dots,v_{i+k})}{\omega^*(v_{i+1},v_{i+2},\dots,v_{i+k})}v_{i+1}-v_i
\end{equation}
Clearly this vector satisfies $v_iv_{i+1} =v_{i+1}w_1$, and it is in $\Span\{v_{i+2},\dots,v_{i+k}\}$ because its exterior product with $v_{i+2} \cdots v_{i+k}$ vanishes. It is occasionally convenient to clear denominators in \eqref{eq:sigmaidiagram}, replacing $w_1$ by $\tilde{w}_1 = \omega^*(v_i,v_{i+2},\dots,v_{i+k})v_{i+1}-\omega^*(v_{i+1},v_{i+2},\dots,v_{i+k})v_i$, and replacing $w_\ell$ by its corresponding $\tilde{w}_\ell$. This recipe is defined on $V^n$, not just on $(V^n)^\circ$. The resulting \emph{renormalized Artin generators} $\tilde{\sigma}_i \colon \tGr(k,n) \to \tGr(k,n)$ only satisfy the braid relations up to monomials in the frozen variables. The pullbacks $\tilde{\sigma_i}^* \propto \sigma_i^*$ are proportional. 
\end{rmk}

\begin{rmk}[Cluster variables and webs]
The renormalized Artin generator can be understood ``diagrammatically'' in terms of $\SL_k$ webs (cf.~\cite{SkewHowe} for the definition of $\SL_k$ webs). The Artin generator $\sigma_i$ can be thought of as a map $V^{\otimes n} \to V^{\otimes n}$ that can be encoded as a web, because the renormalized map $(v_i,\dots,v_{i+k}) \mapsto \tilde{w}_1$ from Remark~\ref{rmk:renormd} has a pictorial interpretation in terms of $\SL_k$ tensor diagrams. In this way, one can convert any braid diagram to an $\SL_k$ tensor diagram that computes the corresponding quasi-automorphism. This tensor diagram will not be planar, but it turns out that one it can be ``planarized'' using the $\SL_k$ crossing removal relation (also known as the \emph{braiding} skein relation cf~\cite[Corollary 6.2.3]{SkewHowe}). We illustrate these ideas in the case of $\SL_3$ webs and $\SL_4$ webs in Section~\ref{secn:proofs}. Since each Pl\"ucker coordinate is an $\SL_k$ web, Theorem~\ref{thm:braidgroupacts} can be used to construct (conjecturally infinitely many) new cluster variables in $\bbc[\tGr(k,n)]$ that are proportional to $\SL_k$ webs, provided $\gcd(k,n) > 1$.  
\end{rmk}


We close this section by giving a proof of Theorem~\ref{thm:braidgroupacts} in steps, starting with claim 1), then claim 2), then claim 3). The proof of claim 3) is the most involved. The details of these proofs are not used in subsequent sections. 
\begin{proof}[Proof Theorem~\ref{thm:braidgroupacts} claim 1)]
Since $\rho^{-1} \circ \sigma_i \circ \rho = \sigma_{i+1}$ and $\sigma_{i}^{-1} = \tau \circ \sigma_{d-i} \circ \tau$, it suffices to prove our claims for $\sigma_1$. 
First we argue that the pullback of any frozen variable along the map $\sigma_1 \colon (V^n)^\circ \to (V^n)^\circ$ is a Laurent monomial in frozen variables. Explicitly:
\begin{equation}\label{eq:frozenpullback}
\sigma_1^*(\Delta_{i,\dots,i+k-1}) = \begin{cases}
\Delta_{i,\dots,i+k-1} \text{ if $i \not \equiv 2 \mod d$} \\
\frac{\Delta_{i-1,\dots,i+k}\Delta_{i+1,\dots,i+k}}{\Delta_{i,\dots,i+k-1}} \text{ otherwise.}
\end{cases}
\end{equation}
To prove this formula, we need to calculate the exterior product of the vectors in locations $i,\dots,i+k-1$ of $\sigma_1(v_1,\dots,v_n)$. We can group the terms in this exterior product window by window. Since  $v_{(\ell-1)d+1}v_{(\ell-1)d+2} = v_{(\ell-1)d+2}w_\ell$ in the $\ell$th window, for any window in which the numbers $(\ell-1)d+1$ and $(\ell-1)d+2$ are both in our interval [i,i+k-1], the exterior product of the vectors in this window of $\sigma_1(v_1,\dots,v_n)$ agrees with the exterior product of the the vectors $(v_1,\dots,v_n)$. By this reasoning, one concludes that $\Delta_{i,\dots,i+k-1}(\sigma_1(v_1,\dots,v_n)) = \Delta_{i,\dots,i+k-1}(v_1,\dots,v_n)$ whenever $i \not \equiv 2 \mod d$. For the other cases, we can focus on $i=2$. Then 
\begin{align*}
\Delta_{2,\dots,k+1}(\sigma_1(v_1,\dots,v_n)) &= \omega^*(w_1 v_3 \cdots v_d v_{d+2} w_2 v_{d+3} \cdots v_{2d} \cdots v_{k} v_{k+2}) \\ 
&= \omega^*(w_1 v_3 \cdots v_{k} v_{k+2}).
\end{align*}
By definition, $w_1 \in \Span\{v_3,\dots,v_{k+1}\}$. So the tensors $w_1 v_3 \cdots v_{k} v_{k+2}$ and $v_3 \cdots v_{k+2}$ agree up to a scalar multiple (namely, the coefficient of $v_{k+1}$ when $w_1$ is expanded in terms of $\{v_3,\dots,v_{k+1}\}$). To compute this scalar multiple, we note that it is the same as the scalar multiple relating the tensors $v_2 w_1 v_3 \cdots v_k$ and $v_2 v_3 \cdots v_{k+1}$. From the normalization condition on $w_1$, one has $v_2 w_1 v_3 \cdots v_k = v_1 \cdots v_k$. This establishes $\Delta_{2,\dots,k+1}(\sigma_1(v_1,\dots,v_n)) = \frac{\Delta_{1,\dots,k}\Delta_{3,\dots,k+2}}{\Delta_{2,\dots,k+1}}|_{(v_1,\dots,v_n)}$
as claimed. 

Second, the regularity of $\sigma_1$ follows from the manifestly polynomial formula \eqref{eq:sigmaidiagram}. Third, the map on $(V^n)^\circ$ descends to a regular map on the open positroid stratum $\tGr^\circ(k,n)$ because the conditions in \eqref{eq:sigmai} are $\SL(V)$-equivariant. Fourth and finally, we compute the composite 
$[v_1,\dots,v_d] \overset{\sigma_1}{\mapsto} [v_2,w_1,v_3,\dots,v_d] \overset{\sigma_1^{-1}}{\mapsto}[u_1,v_2,\dots,v_d].$
The definitions of $w_1$ and $u_1$ imply that $u_1 \in \Span \{w_{\frac{n}{d}},v_{n-k+3},\dots,v_n\} \cap \Span \{v_2,w_1\} = \Span \{w_{\frac{n}{d}},v_{n-k+3},\dots,v_n\} \cap \Span \{v_1,v_2\},$
and also that $w_{\frac{n}{d}} \in \Span\{v_{n-k+3},\dots,v_n,v_1\}$. Hence $u_1 \in \Span\{v_{n-k+3},\dots,v_1\} \cap \Span\{v_1,v_2\} = \Span\{v_1\}$.  So $u_1$ and $v_1$ are related by a scalar multiple, and this scalar multiple equals one using $u_1 v_2 = v_2 w_1 = v_1v_2$. The other composite is similar. 
\end{proof}

\begin{proof}[Proof Theorem~\ref{thm:braidgroupacts} claim 2)] By the reasoning in Remark~\ref{rmk:kdvidesn}, the $i$th Artin generator when $d < k$ are 
is constructed as $\prod_{j=0}^{\frac{k}{d}-1} \sigma_{i+jd}$ where $\sigma_{i}$'s are Artin generators for $\tGr(k,\frac{n}{d}k)$, and we diagonally embed $\tGr(k,n) \subset \tGr(k,\frac{n}{d}k)$. These braids $\prod_{j=0}^{\frac{k}{d}-1} \sigma_{i+jd}$ clearly satisfy the braid relations provided the $\sigma_i$ do. 

Thus we henceforth assume that $k$ divides $n$, i.e. $d=k$. Let us check that $\sigma_i \circ \sigma_j = \sigma_j \circ \sigma_i$ when $|i-j| \geq 2$. As in the proof of Lemma~\ref{lem:consecutivecomp}, we can define a vector $w_\ell$ to be the new vector in the $\ell$th window when we apply $\sigma_i$ to $(v_1,\dots,v_n)$, and $u_\ell$ to be the new vector in the $\ell$th window when we apply $\sigma_j$ to $\sigma_i(v_1,\dots,v_n)$. The vectors $w_1$ and $u_1$ in the first window are defined by conditions 
\begin{align*}
w_1 \in & \Span\{v_{i+2},\dots,v_{i+k}\} \\ 
u_1 \in & \Span\{v_{j+2},\dots,v_k,v_{k+2},w_2,v_{k+j}\} =\Span\{v_{j+2},\dots,v_{k+j}\},
\end{align*}
together with the normalizations $v_iv_{i+1} = v_{i+1}w_1$ and $v_jv_{j+1} = v_{j+1}u_1$. After performing the above simplification, these conditions are manifestly symmetric in $i$ and $j$, proving the claim. 

Next, we check that $\sigma_1 \circ \sigma_2 \circ \sigma_1= \sigma_2 \circ \sigma_1 \circ \sigma_2$. The general result holds after conjugating by the twisted cyclic shift. As before, we define vectors $w_\ell,u_\ell$, and $z_\ell$ to be the new vectors created by applying $\sigma_1$, then $\sigma_2$, and then $\sigma_1$. In the first window: 
\begin{equation}
[v_1,\dots,v_k] \overset{\sigma_1}{\mapsto} [v_2,w_1,v_3,\dots,v_k] \overset{\sigma_2}{\mapsto} [v_2,v_3,u_1,v_4\dots,v_k] \overset{\sigma_1}{\mapsto} [v_3,z_1,u_1,v_4\dots,v_k]. \label{eq:121}
\end{equation}
The conditions defining the vectors $w_1,u_1,z_1$ are the following: 
\begin{align*}
w_1 \in &  \Span\{v_{3},\dots,v_{k+1}\}, \\
u_1 \in & \Span\{v_{4},\dots,v_{k},v_{k+2},w_2\} = \Span\{v_{4},\dots,v_{k+2}\} \\
z_1 \in & \Span\{u_1,v_{4},\dots,v_{k},v_{k+2}\} = \Span\{w_2,v_{4},\dots,v_{k},v_{k+2}\} =  \Span\{v_{4},\dots,v_{k+2}\},
\end{align*}
together with the normalizations $v_1v_2 = v_2w_1$, $w_1v_3 = v_3u_1$, and $v_2v_3 = v_3z_1$. 
 

For the other composition, we define vectors $z'_\ell$, $x'_\ell$, and $u'_\ell$ created when applying $\sigma_2$, then $\sigma_1$, then $\sigma_1$. In the first window: 
\begin{equation}
[v_1,\dots,v_k] \overset{\sigma_2}{\mapsto} [v_1,v_3,z'_1,\dots,v_k] \overset{\sigma_1}{\mapsto} [v_3,x'_1,z'_1,v_4\dots,v_k] \overset{\sigma_2}{\mapsto} [v_3,z'_1,u'_1,v_4\dots,v_k]. \label{eq:212}
\end{equation}
By comparing the conditions, we see immediately that $z'_\ell = z_\ell$. The vectors $x'_1$ and $u'_1$ are defined by 
\begin{align*}
x'_1 \in & \Span\{z_1,v_{4},\dots,v_{k+1}\}  =    \Span\{v_{4},\dots,v_{k+2}\} \\  
u'_1 \in & \Span\{v_{4},\dots,v_k,v_{k+3},x'_2\} = \Span\{v_{4},\dots,v_k,v_{k+1},v_{k+3}\},
\end{align*}
together with the normalizations $v_1v_3 = v_3 x'_1$ and $x'_1z_1 = z_1 u'_1$. 

We need to establish that $u_1 = u'_1$. First, we explain that both vectors are in $\Span\{v_1,v_2,v_3\} \cap \Span\{v_{4},\dots,v_{k+1}\}$, which is a line by consecutive genericity. Indeed, from the conditions defining $x'_1$ and $z_1$ we conclude that $u'_1 \in \Span\{v_1,v_2,v_3\}$, and also that $u'_1 \in \Span\{v_4,\dots,v_{k+2}\}$. But then  $u'_1 \in \Span\{v_4,\dots,v_{k+2}\} \cap \Span \{v_{4},\dots,v_k,v_{k+1},v_{k+3}\} =\Span\{v_4,\dots,v_{k+1}\}$, as claimed. On the other hand,  $u_1 \in \Span\{v_1,v_2,v_3\}$ since this is true of $w_1$, and furthermore $u_1 \in \Span\{v_3,\dots,v_{k+1}\} \cap \Span\{v_4,\dots,v_{k+2}\} = \Span\{v_4,\dots,v_{k+1}\}$. Thus, $u_1$ and $u'_1$ differ by a scalar multiple. This scalar multiple equals one by the normalization conditions: 
$$v_2v_3u_1 = v_2w_1v_3 = v_1v_2v_3 = v_1 v_3 z_1 = v_3 x'_1 z_1 = v_3 z_1 u'_1 = v_2v_3u'_1.$$
\end{proof}

The proof of claim 3) in Theorem~\ref{thm:braidgroupacts} is more involved, because one has to show explicitly that the Artin generators can be implemented by mutations. Let us begin proving this. 

\subsection{The Artin generators are quasi-automorphisms}
Since $\rho^{-1} \sigma_i \rho = \sigma_{i+1}$ and $\sigma_i^{-1} = \tau \sigma_{d-i} \tau$, we can deduce that each Artin generator is a quasi-automorphism once we prove this is true of $\sigma_1^*$. Our first step is to describe a cluster on which the action of $\sigma_1^*$ is especially convenient. As usual, we let $d = \gcd(k,n)$ and assume $d \geq 2$. 

For a $k$-subset $I = [1,a]\} \cup \{[b+1,b+k-a] \in \mcc_{\Le}$, we define its modified version $\Imod$ as follows: 
\begin{equation}
\label{eq:Imod}
\Imod = 
\begin{cases}
I, \text{ if $\Delta_I$ is frozen} \\
I, \text{ if $b+1 \not \equiv 2 \mod d$} \\
[1,a-1] \cup \{b\} \cup [b+2,b+k-a+1] \text{ if $b+1 \equiv 2 \mod d$.}
\end{cases}
\end{equation}

We let $\mccmod = \{\Imod \colon I \in \mcc_{\Le}\}$ be the modified version of the $\Le$ cluster. 
In the $\Gr(4,8)$ example (cf.~Figure~\ref{fig:stdseed48}), the modified version of $\Delta_{1236}$ is $\Delta_{1257}$, the modified version of $\Delta_{1267}$ is $\Delta_{1578}$, and no other Pl\"ucker coordinates are modified.  

By construction, the cluster $\mccmod$ has the following key property: if $I \in \mccmod$ is a non-frozen Pl\"ucker coordinate, and if $i \in I$ satisfies $i \equiv 2 \mod d$, then $i-1 \in I$. We make frequent use of this property in our subsequent lemmas.  

Let $\mfs_n$ denote the symmetric group on $n$ symbols. Define $\pi \in \mfs_n$ as the product of commuting transpositions
$\pi = \prod_{j=0}^{\frac{n}{d}-1}(jd+1,jd+2).$ This permutation acts by switching adjacent numbers that are equivalent to $1$ and $2$ $\mod d$.

\begin{lem}\label{lem:sigmaispi}
If $\Imod \in \mccmod$ is a non-frozen $k$-subset, then $\sigma_1^*(\Delta_{\Imod}) = \Delta_{\pi(\Imod)}$.  
\end{lem}

\begin{proof}
This follows from the key property of the collection $\mccmod$ alluded to above, and the normalization condition $v_{(\ell-1)d+1}v_{(\ell-1)d+2} = v_{(\ell-1)d+2}w_\ell$ in the $\ell$th window. 
\end{proof}

That is, on our convenient choice of weakly separated collection $\mccmod$, the Artin generator $\sigma_1$ acts on the non-frozen Pl\"ucker coordinates by permuting indices according to~$\pi$. 

We let $\pi(\mccmod)$ denote the collection of Pl\"ucker coordinates obtained from $\mccmod$ by applying $\pi$ to each of the non-frozen Pl\"ucker coordinates, and doing nothing to the frozen coordinates.

To show that $\sigma_1^*$ is a quasi-automorphism, we must show that both collections $\mccmod$ and $\pi(\mccmod)$ are weakly separated (hence are clusters in $\tGr(k,n)$), and furthermore check that $\sigma_1^*$ maps exchange ratios in the seed $\Sigma(\mccmod)$ to those in $\Sigma(\pi(\mccmod))$. This is the content of the next several lemmas, each of which is an exercise in weak separation combinatorics, and may be skipped without affecting the later sections. 

\begin{lem}\label{lem:mccmodisws} The collection $\mccmod$ is a maximal weakly separated collection. 
\end{lem}

\begin{proof}
We give an explicit sequence of square moves $\mcc_{\Le} \to \mccmod$. Recall that elements of $\mcc_{\Le} \setminus \{1,\dots,k\}$ are identified with points in a rectangular array. We decompose this rectangular array into $n-1$ diagonals, each of which starts in either the first column or last row of the array, and moves northeast. The first diagonal is the entry in the first row and column, and the last diagonal is the entry in the last row and column. 

Then the Pl\"ucker coordinates $\Delta_I \in \mcc_{\Le}$ that are affected by performing the modification $I \mapsto \Imod$ \eqref{eq:Imod} are exactly those sitting on the diagonals $d,2d,3d,\dots,n-d$. Since $d \geq 2$, if two Pl\"ucker coordinates $\Delta_I$ and $\Delta_J$ are in different affected diagonals, then they are not adjacent to each other in the quiver $\tQ(\mcc_{\Le})$. 

Then the promised sequence of square moves $\mcc \to \mccmod$ is the one in which we mutate once along each of these affected diagonals, starting with southwest entry in this diagonal, continuing to the penultimate entry (because the final entry is frozen). 
It is straightforward to see that in a given diagonal, when we mutate a given Pl\"ucker coordinate  $I$, it is replaced by its modified version $\Imod$ (the corresponding exchange relation is a three-term Pl\"ucker relation). Since affected diagonals are not connected by arrows in $Q(\mcc_{\Le})$, the mutations in the various diagonals commute.  
\end{proof}

\begin{lem}\label{lem:pimodisws}  The collection $\pi(\mccmod)$ is a maximal weakly separated collection. 
\end{lem}

\begin{proof} Let $I,J \in \mcc_{\Le}$ be a pair of non-frozen $k$-subsets; we need to show that $\pi(\Imod)$ and $\pi(\Jmod)$ are weakly separated. Let us suppose the contrary, so there are cyclically ordered elements $a < b < c < d \in [n]$, with $a,c \in \pi(\Imod) \setminus\pi(\Jmod)$ and $b,d \in \pi(\Jmod) \setminus\pi(\Imod)$. Then we claim that the elements $\pi(a),\dots,\pi(d)$ remain cyclically sorted: indeed, the only concern would be that a pair of entries are swapped by $\pi$, e.g. $b = a+1$ and $b \equiv 2 \mod d$. But by the key property of $\mccmod$, the assumption $a \in \pi(\Imod)$ would imply that $a+1 \in \pi(\Imod)$ also, a contradiction.  Furthermore, it is clear that one still has $\pi(a),\pi(c) \in \Imod \setminus \Jmod$ and $\pi(b),\pi(d) \in \Jmod \setminus \Imod$ (because $\pi$ is an involution). Thus, the elements $\pi(a),\dots,\pi(d)$ are a witness to the fact that the sets $\Imod$ and $\Jmod$ are not weakly separated, contradicting Lemma~\ref{lem:mccmodisws}. So our initial assumption on the existence of $a,b,c,d$ was incorrect.    
\end{proof}

\begin{lem}\label{lem:yhatsidentified} The map $\sigma_1^*$ sends exchange ratios in $\Sigma(\mccmod)$ to the corresponding ones in $\Sigma(\pi(\mccmod))$, i.e. 
$\sigma_1^*(\hat{y}_{\Sigma(\mccmod)}(\Imod)) = \hat{y}_{\Sigma(\pi(\mccmod))}(\pi(\Imod))$ for each non-frozen $\Imod \in \mccmod$. 
\end{lem}

\begin{proof} 
We use the description of the extended quiver $\tQ(\mccmod)$ in terms of nontrivial cliques. Let $\mcw(S)$ be a nontrivial clique in $\mccmod$. Let us enumerate the non-frozen subsets in $\mcw(S)$ as $I_1,\dots,I_t$ so that the singletons $S \setminus I_1,\dots,S \setminus I_t$ are sorted in increasing order. Since any pair of elements in a  clique intersect in $k-1$ elements, $\mcw(S)$ has either zero, one, or two frozen variables. Likewise, we have that $t \geq 1$. 

First, we observe that the singletons $\pi(S) \setminus \pi(I_1),\dots,\pi(S) \setminus \pi(I_t)$ remain sorted. Indeed, since $\pi$ is a product of transpositions, it could only swap the relative positions of adjacent singletons, which were congruent to 1 and 2 $\mod d$ respectively. But the key property of the cluster $\mccmod$ disallows this. 

Second, we observe that a clique $\mcw(S)$ contains zero frozen variables if and only if $\mcw(\pi(S))$ contains zero frozen variables. Indeed, suppose $[b+1,\dots,b+k] \subset \pi(S)$. Then $[b+1,\dots,b+k] \subset S$ follows unless $b \equiv 1 \mod d$ and $b+1 \notin S$. But if $b \equiv 1 \mod d$, then $b+k \equiv 1 \mod d$, so $b+k \notin \pi(S)$ unless $b+k+1 \in \pi(S)$. But then $[b+2,\dots,b+k+1] \subset S$ contradicting the zero frozen variable assumption. The reverse argument is similar. 

For any nontrivial clique $\mcw(S) \subset \mccmod$ with $t\geq 2$, the second observation implies that $\mcw(\pi(S))$ is also a nontrivial clique with $t \geq 2$, and vice versa. Since mutable arrows in $\tQ(\mccmod)$ only come from cliques in which $t \geq 2$, the first and second observations together imply that mutable arrows in $\tQ(\mccmod)$ exactly correspond to mutable arrows in $\tQ(\pi(\mccmod))$. That is, we have so far verified that $\pi$ induces an isomorphism $Q(\mccmod) \cong Q(\pi(\mccmod))$. 

It remains to address the frozen variables. We let $f_i = [i,\dots,i+k-1]$ denote the $i$th frozen subset. We recall that $\sigma_1^*$ fixes the frozen variable $\Delta_{f_i}$, unless $i \equiv 2 \mod d$, cf.~\eqref{eq:frozenpullback}.  

By the key property of $\mccmod$, if $i \equiv 1 \mod d$, then a nontrivial clique for $\mccmod$ cannot contain the subset $f_{i+1}$ unless it also contains the frozen subset $f_{i}$. Furthermore, if a nontrivial clique does not contain a frozen variable $f_{i+1}$ where $i \equiv 1 \mod d$, then one verifies that any frozen variables that are in $\mcw(S)$ are also in $\mcw(\pi(S))$. Dually, if a nontrivial clique in $\pi(\mccmod)$ contains the subset $f_{i+1}$ where $i \equiv 1 \mod d$, then it must contain the frozen subset $f_{i+2}$. And if it does not contain such an $f_{i+1}$, then any frozen variables in $\mcw(\pi(S))$ are also in $\mcw(S)$. 

By the explicit description \eqref{eq:Imod}, we see that if $i \equiv 1 \mod d$ and $S = [i,\dots,i+k]$, then the clique 
$\mcw(S) \subset \mccmod$ is nontrivial, and consists of three elements $\mcw(S) = \{f_i,f_{i+1},\Imod\}$ where $\Imod = \{i,i+2,\dots,i+k\}$ is the modified version of $I =\{1,i+1,\dots,i+k-1\} \in \mcc_{\Le}$. We have $\pi(\Imod) = [i+1,\dots,i+k-1,i+k+1]$, so the clique $\{f_{i+1},f_{i+2},\pi(\Imod)\} = \mcw([i+1,\dots,i+k+1])$ is a nontrivial clique in $\pi(\mccmod)$. This is compatible with \eqref{eq:frozenpullback}: the frozen part of 
$\hat{y}(\Imod)$ is $\frac{f_i}{f_{i+1}}$, and we have $\sigma_1^*(\frac{f_i}{f_{i+1}}) = \frac{f_i f_{i+1}}{f_{i}f_{i+2}}$, which agrees with the frozen part of $\hat{y}(\pi(\Imod))$. 

On the other hand, for nontrivial cliques that do not contain $f_{i+1}$ where $i \equiv 1 \mod d$, the frozen part of each $\hat{y}(\Imod)$ becomes the the frozen part of $\hat{y}(\pi(\Imod))$, which matches what happens when we apply $\sigma_1^*$. 
\end{proof}

\begin{proof}[Proof of Theorem~\ref{thm:braidgroupacts} claim 3)]
Since the various $\sigma_i^*$ are related to each other by cluster automorphisms, and we already know they satisfy the braid relations, we only need to explain why 
$\sigma_1^*$ is a quasi-automorphism. First, one has that $\sigma_1^*(\bbp) \subset \bbp$ from the argument at the start of the proof of claim 1) of Theorem~\ref{thm:braidgroupacts}. Second, by Lemma~\ref{lem:mccmodisws} and Lemma~\ref{lem:pimodisws}, the collections $\mccmod$ and $\pi(\mccmod)$ are maximal weakly separated collections, thus they determine seeds in the cluster structure on $\tGr(k,n)$. And by Lemma~\ref{lem:sigmaispi} and \ref{lem:yhatsidentified}, $\sigma_1^*$ takes the cluster $\mccmod$ to the cluster $\pi(\mccmod)$ while preserving the exchange ratios. Then by definition, $\sigma_1^*$ is an orientation-preserving quasi-automorphism. 
\end{proof}



\section{Fock-Goncharov cluster algebras}\label{secn:FockGoncharovSpaces}
We retain all the conventions and notation concerning the exterior algebra $\bigwedge(V)$ from the previous sections. If $v \in \bigwedge^a(V)$ and $w \in \bigwedge^{a+1}(V)$ are simple tensors, we say that \emph{$v$ divides $w$} if $w = vv'$ for some vector    $v' \in V$. 
 
\begin{defn}\label{defn:affineflag} An \emph{affine flag in $V$} is a sequence 
\begin{equation}\label{eq:affineflag}
F_\bullet = F_{(1)}, F_{(2)}, \dots, F_{(k)} 
\end{equation} of simple anti-symmetric tensors $F_{(a)} \in \bigwedge^a(V)$, each one dividing the next, whose top flag $F_{(k)}$ is the volume form~$\omega$. We let $\textnormal{Aff}(V)$ denote the space of affine flags in $V$. A pair of affine flags $F_{\bullet},G_{\bullet}$ meet \emph{generically} if each their tensors of complementary size meet generically, i.e. $\omega^*(F_{(a)}G_{(k-a)}) \neq 0$ for $a=1,\dots,k-1$. We denote by $(\textnormal{Aff}(V)^r)^\circ \subset \textnormal{Aff}(V)^r$ the space of consecutively generic $r$-tuples, i.e. those in which adjacent affine flags meet generically. We let $\Conf(k,r) = \SL(V) \backslash \textnormal{Aff}(V)^r$ denote the space of configurations of $r$ affine flags, and $\Conf^\circ(k,r)$ the subset of consecutively generic configurations. 
\end{defn}

Geometrically, an affine flag is an ordinary flag together with a choice of volume form at each step of the flag. The action of $\SL(V)$ on affine flags is transitive, and the stabilizer of any particular affine flag is a maximal unipotent subgroup~$U$. Thus we can identify $\textnormal{Aff}(V)$ with $\SL(V) / U$.  

Fock and Goncharov described a cluster algebra structure in the field of rational functions $\bbc(\Conf(k,r))$ .~\cite[Section 9]{ModuliSpaces}. We will summarize this construction. One begins with the space $\Conf(k,3)$, i.e. configurations of three affine flags in $V$. Consider a configuration, represented by affine flags $(F_{1,\bullet},F_{2,\bullet},F_{3,\bullet})$. Let $(a,b,c)$ be a triple of nonnegative integers satisfying $a+b+c = k$, at least  two of which are  positive. We define the \emph{Fock-Goncharov coordinate} by
\begin{equation}\label{eq:FGCoord}
\Delta_{a,b,c}(p) = \omega^*(F_{1,(a)}  F_{2,(b)} F_{3,(c)}), 
\end{equation} 
noting that the right-hand side of \eqref{eq:FGCoord} does not depend on the choice affine flags~$F_{i,\bullet}$ representing the configuration. 

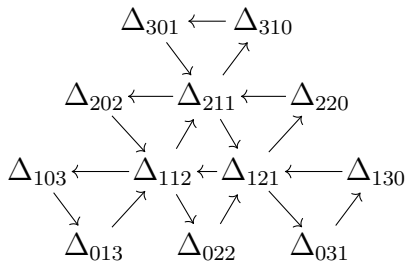
\begin{figure}[ht]
\begin{center}
\begin{tikzpicture}[scale = 1]
\coordinate (AAAC) at (-.75,3);
\coordinate (AAAB) at (.75,3);
\coordinate (AACC) at (-1.5,2);
\coordinate (AABC) at (0,2);
\coordinate (AABB) at (1.5,2);
\coordinate (ACCC) at (-2.25,1);
\coordinate (ABCC) at (-.585,1);
\coordinate (ABBC) at (.585,1);
\coordinate (ABBB) at (2.25,1);
\coordinate (BCCC) at (-1.5,0);
\coordinate (BBCC) at (0,0);
\coordinate (BBBC) at (1.5,0);

\node at (AAAC) {$\Delta_{301}$};
\node at (AAAB) {$\Delta_{310}$};
\node at (AACC) {$\Delta_{202}$};
\node at (AABC) {$\Delta_{211}$};
\node at (AABB) {$\Delta_{220}$};
\node at (ACCC) {$\Delta_{103}$};
\node at (ABCC) {$\Delta_{112}$};
\node at (ABBC) {$\Delta_{121}$};
\node at (ABBB) {$\Delta_{130}$};
\node at (BCCC) {$\Delta_{013}$};
\node at (BBCC) {$\Delta_{022}$};
\node at (BBBC) {$\Delta_{031}$};

\draw [shorten >=0.5cm,shorten <=0.5cm,->] (AAAB)--(AAAC);
\draw [shorten >=0.45cm,shorten <=0.45cm,->] (AABB)--(AABC);
\draw [shorten >=0.45cm,shorten <=0.45cm,->] (AABC)--(AACC);
\draw [shorten >=0.45cm,shorten <=0.45cm,->] (ABBB)--(ABBC);
\draw [shorten >=0.45cm,shorten <=0.45cm,->] (ABBC)--(ABCC);
\draw [shorten >=0.45cm,shorten <=0.45cm,->] (ABCC)--(ACCC);

\draw [shorten >=0.35cm,shorten <=0.35cm,->] (AAAC)--(AABC);
\draw [shorten >=0.35cm,shorten <=0.35cm,->] (AABC)--(AAAB);

\draw [shorten >=0.35cm,shorten <=0.35cm,->] (AACC)--(ABCC);
\draw [shorten >=0.35cm,shorten <=0.35cm,->] (ABCC)--(AABC);
\draw [shorten >=0.35cm,shorten <=0.35cm,->] (AABC)--(ABBC);
\draw [shorten >=0.35cm,shorten <=0.35cm,->] (ABBC)--(AABB);

\draw [shorten >=0.35cm,shorten <=0.35cm,->] (ACCC)--(BCCC);
\draw [shorten >=0.35cm,shorten <=0.35cm,->] (BCCC)--(ABCC);
\draw [shorten >=0.35cm,shorten <=0.35cm,->] (ABCC)--(BBCC);
\draw [shorten >=0.35cm,shorten <=0.35cm,->] (BBCC)--(ABBC);
\draw [shorten >=0.35cm,shorten <=0.35cm,->] (ABBC)--(BBBC);
\draw [shorten >=0.35cm,shorten <=0.35cm,->] (BBBC)--(ABBB);
\end{tikzpicture}
\caption{A triangular array of Fock-Goncharov coordinates for $\Conf(4,3)$. The three ``corners'' of the triangle are not considered part of the array. \label{fig:FGQuiver}}
\end{center}
\end{figure}

The Fock-Goncharov cluster structure on $\Conf(k,3)$ is obtained from an initial seed whose initial extended cluster consists of all such Fock-Goncharov coordinates \eqref{eq:FGCoord}. The frozen variables in this extended cluster are the $\Delta_{a,b,c}$ in which one of $a$,$b$, or $c$ is $0$. To describe the extended quiver in this initial seed, we arrange the Fock-Goncharov coordinates in a triangular array, drawing directed arrows between adjacent entries so that every small triangle in the diagram is oriented counterclockwise, cf.~Figure~\ref{fig:FGQuiver}. 

Now we return to the general case of $r$ affine flags.  

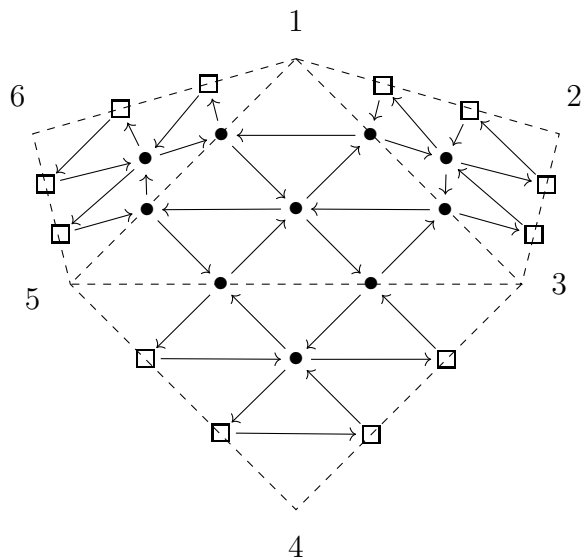
\begin{figure}
\begin{center}
\begin{tikzpicture}
\coordinate (A) at (0,3);
\node at (0,3.5) {$1$};
\coordinate (D) at (0,-3);
\node at (0,-3.5) {$4$};
\coordinate (B) at (3.5,2);
\node at (3.7,2.5) {$2$};
\coordinate (F) at (-3.5,2);
\node at (-3.7,2.5) {$6$};
\coordinate (C) at (3,0);
\node at (3.5,0) {$3$};
\coordinate (E) at (-3,0);
\node at (-3.5,-.2) {$5$};

\draw [dashed] (A)--(B)--(C)--(D)--(E)--(F)--(A);
\draw [dashed] (A)--(C)--(E)--(A);

\coordinate (AAB) at (.66*0+.33*3.5,.66*3+.33*2);
\node at (AAB) {$\boxed{}$};
\coordinate (ABB) at (.33*0+.66*3.5,.33*3+.66*2);
\node at (ABB) {$\boxed{}$};

\coordinate (BBC) at (.666*3.5+.333*3,.66*2+.33*0);
\node at (BBC) {$\boxed{}$};
\coordinate (BCC) at (.333*3.5+.666*3,.33*2+.66*0);
\node at (BCC) {$\boxed{}$};

\coordinate (CCD) at (.666*3+.333*0,.66*0-.333*3);
\node at (CCD) {$\boxed{}$};
\coordinate (CDD) at (.333*3+.666*0,.33*0-.666*3);
\node at (CDD) {$\boxed{}$};

\coordinate (DDE) at (.666*0-.333*3,-.66*3+.333*0);
\node at (DDE) {$\boxed{}$};
\coordinate (DEE) at (.333*0-.666*3,-.33*3+.666*0);
\node at (DEE) {$\boxed{}$};

\coordinate (EEF) at (-.66*3-.33*3.5,.666*0+.333*2);
\node at (EEF) {$\boxed{}$};
\coordinate (EFF) at (-.333*3-.666*3.5,.333*0+.666*2);
\node at (EFF) {$\boxed{}$};

\coordinate (AFF) at (-.666*3.5+.333*0,.666*2+.333*3);
\node at (AFF) {$\boxed{}$};
\coordinate (AAF) at (-.333*3.5,.333*2+.666*3);
\node at (AAF) {$\boxed{}$};

\coordinate (AAC) at (.66*0+.33*3,.66*3+.33*0);
\node at (AAC) {$\bullet$};
\coordinate (ACC) at (.33*0+.66*3,.33*3+.66*0);
\node at (ACC) {$\bullet$};

\coordinate (AAE) at (.66*0-.33*3,.66*3+.33*0);
\node at (AAE) {$\bullet$};
\coordinate (AEE) at (.33*0-.66*3,.33*3+.66*0);
\node at (AEE) {$\bullet$};

\coordinate (CCE) at (1,0);
\node at (CCE) {$\bullet$};
\coordinate (CEE) at (-1,0);
\node at (CEE) {$\bullet$};

\coordinate (ACE) at (0,1);
\node at (ACE) {$\bullet$};

\coordinate (CDE) at (0,-1);
\node at (CDE) {$\bullet$};

\coordinate (ABC) at (2,1.66);
\node at (ABC) {$\bullet$};

\coordinate (AEF) at (-2,1.66);
\node at (AEF) {$\bullet$};

\draw [shorten >=0.2cm,shorten <=.2cm,->] (AAC) -- (AAE);
\draw [shorten >=0.2cm,shorten <=.2cm,->] (AAE) -- (ACE);
\draw [shorten >=0.2cm,shorten <=.2cm,->] (ACE) -- (AAC);
\draw [shorten >=0.2cm,shorten <=.2cm,->] (ACC) -- (ACE);
\draw [shorten >=0.2cm,shorten <=.2cm,->] (ACE) -- (CCE);
\draw [shorten >=0.2cm,shorten <=.2cm,->] (CCE) -- (ACC);
\draw [shorten >=0.2cm,shorten <=.2cm,->] (ACE) -- (AEE);
\draw [shorten >=0.2cm,shorten <=.2cm,->] (AEE) -- (CEE);
\draw [shorten >=0.2cm,shorten <=.2cm,->] (CEE) -- (ACE);

\draw [shorten >=0.2cm,shorten <=.2cm,->] (AAB) -- (AAC);
\draw [shorten >=0.2cm,shorten <=.2cm,->] (AAC) -- (ABC);
\draw [shorten >=0.2cm,shorten <=.2cm,->] (ABC) -- (AAB);
\draw [shorten >=0.2cm,shorten <=.2cm,->] (ABB) -- (ABC);
\draw [shorten >=0.2cm,shorten <=.2cm,->] (ABC) -- (BBC);
\draw [shorten >=0.2cm,shorten <=.2cm,->] (BBC) -- (ABB);
\draw [shorten >=0.2cm,shorten <=.2cm,->] (ABC) -- (ACC);
\draw [shorten >=0.2cm,shorten <=.2cm,->] (ACC) -- (BCC);
\draw [shorten >=0.2cm,shorten <=.2cm,->] (BCC) -- (ABC);

\draw [shorten >=0.2cm,shorten <=.2cm,->] (AAE) -- (AAF);
\draw [shorten >=0.2cm,shorten <=.2cm,->] (AAF) -- (AEF);
\draw [shorten >=0.2cm,shorten <=.2cm,->] (AEF) -- (AAE);
\draw [shorten >=0.2cm,shorten <=.2cm,->] (AEE) -- (AEF);
\draw [shorten >=0.2cm,shorten <=.2cm,->] (AEF) -- (EEF);
\draw [shorten >=0.2cm,shorten <=.2cm,->] (EEF) -- (AEE);
\draw [shorten >=0.2cm,shorten <=.2cm,->] (AEF) -- (AFF);
\draw [shorten >=0.2cm,shorten <=.2cm,->] (AFF) -- (EFF);
\draw [shorten >=0.2cm,shorten <=.2cm,->] (EFF) -- (AEF);

\draw [shorten >=0.2cm,shorten <=.2cm,->] (CDE) -- (DDE);
\draw [shorten >=0.2cm,shorten <=.2cm,->] (DDE) -- (CDD);
\draw [shorten >=0.2cm,shorten <=.2cm,->] (CDD) -- (CDE);
\draw [shorten >=0.2cm,shorten <=.2cm,->] (CDE) -- (CEE);
\draw [shorten >=0.2cm,shorten <=.2cm,->] (CEE) -- (DEE);
\draw [shorten >=0.2cm,shorten <=.2cm,->] (DEE) -- (CDE);
\draw [shorten >=0.2cm,shorten <=.2cm,->] (CDE) -- (CCD);
\draw [shorten >=0.2cm,shorten <=.2cm,->] (CCD) -- (CCE);
\draw [shorten >=0.2cm,shorten <=.2cm,->] (CCE) -- (CDE);

\end{tikzpicture}
\caption{A Fock-Goncharov seed $\Sigma(T)$ for $\bbc(\Conf(3,6)$, i.e. for configurations of $6$ affine flags in $3$-space. $T$ is the triangulation of the hexagon indicated in dashed lines. The extended quiver $\tQ_3(T)$ is drawn in solid lines. There are 10 cluster variables, and 12 frozen variables sitting on the boundary of the hexagon. 
\label{fig:FGSeed}
}
\end{center}
\end{figure}

\begin{defn}\label{defn:FGSeed} Let $D$ be a disk with marked points $1,\dots,r$ in clockwise order on the boundary. For each triangulation $T$ of $D$, we define a seed $\Sigma(T) = (\tilde{\mathbf{x}}(T), \tilde{Q}_k(t))$ in $\bbc(\Conf(k,r))$ as follows. The extended cluster $\tilde{\mathbf{x}}(T)$ is the union of the Fock-Goncharov coordinates coming from the various triangles in $T$. Notice that if $e$ is an internal edge of $T$ then it lies on two triangles, but the Fock-Goncharov coordinates associated to the edge $e$ in either triangle agree as functions on $\FG(k,r)$. The Fock-Goncharov coordinates sitting on the boundary edges of $D$ serve as frozen variables. The extended quiver $\tQ_k(T)$ for this seed is obtained by gluing together the quiver fragments from each triangle in $T$, using the directed edges indicated Figure~\ref{fig:FGQuiver}. See Figure~\ref{fig:FGSeed} for an example. 
\end{defn}

Fock and Goncharov proved that each extended cluster $\tilde{\mathbf{x}(T)}$ coming from a triangulation provides a rational coordinate system on $\Conf(k,r)$. Furthermore, the seeds $\Sigma(T)$, as $T$ varies over all triangulations of the $r$-gon, are related to each other by sequences of mutations. Consequently, the seeds in Definition~\ref{defn:FGSeed} give rise to a well-define cluster structure in $\bbc(\Conf(k,r))$. We let $P \colon \Conf(k,r) \to \Conf(k,r)$ denote the 
twisted cyclic shift $(F_{1,\bullet},\dots,F_{r,\bullet}) \mapsto    (F_{2,\bullet},\dots,F_{r,\bullet},(-\text{Id})^{k-1} (F_{1,\bullet}))$ of affine flags. 

\begin{rmk}\label{rmk:AGS}
The spaces we have called $\Conf(k,r)$ are part of a more general family of spaces $\mca_{G,S}$ defined in \cite{ModuliSpaces} in the context of Higher Teichm\"uller theory. The space $\mca_{G,S}$ is the moduli space of \emph{decorated twisted $G$-local systems on $S$}, for a semisimple Lie group $G$ and a \emph{bordered marked surface} $S$. By the latter we mean mean an orientable Riemann surface $S$, possibly with boundary and punctures, and with marked points on each boundary component. The space $\Conf(k,r)$ is the space $\mca_{\SL_k,S}$ when $S$ is a disk with $r$ marked points on its boundary. 
\end{rmk}


\section{Quasi-isomorphism of configuration spaces}\label{secn:themaps}
We define a quasi-isomorphism between the cluster structures on $\tGr(k,rk)$ and $\Conf(k,2r)$. 

Recall the twisted cyclic shifts $\rho$ on the Grassmannian and $P$ on the configuration space of affine flags. We say a map $\Psi \colon  V^{rk} \to (\textnormal{Aff}(V))^{2r}$ is \emph{$k$-periodic} if it satisfies $\Psi \circ  \rho^k = P^2 \circ \Psi$. This descends to a notion of $k$-periodicity for $\SL(V)$-equivariant maps $\tGr(k,rk) \to \Conf(k,2r)$, and for maps $\tGr^\circ(k,rk) \to \Conf^\circ(k,2r)$ preserving the consecutive generic loci. Similarly a map $\Phi \colon   (\textnormal{Aff}(V))^{2r} \to V^{rk}$ is $k$-periodic if it satisfies $\Phi \circ  P^2  = \rho^k \circ \Phi$. 

Now let $(v_1,\dots,v_n) \in (V^n)^\circ$ be a consecutively generic $n$-tuple in $V$. Define a $k$-periodic map as follows. From the first $k$ vectors, produce a pair of ``opposite'' affine flags  
\begin{align}\label{eq:Psidefnstart}
F_{1,\bullet} = &v_1, v_1v_2, \dots , v_1v_2\cdots v_{k-1},\omega \\
F_{2,\bullet} = &v_k, v_{k-1}v_k, \dots , v_2v_3\cdots v_{k},\omega \label{eq:Psidefnstart2},
\end{align}
and extend by $k$-periodicity. This procedure is $\SL(V)$-equivariant, and determines a rational map
\begin{equation}\label{eq:Psidefn}
\Psi \colon \tGr^\circ(k,rk) \dashrightarrow \Conf^\circ(k,2r). 
\end{equation}
The implicit claim that this map lands in $\Conf^\circ$ is part of Theorem~\ref{thm:PsiPhyXY}. 

\begin{defn}\label{defn:weakmeet}
If $v_1,\dots,v_{k+1}$ are vectors such that both $\{v_1,\dots,v_k\}$ and $\{v_2,\dots,v_{k+1}\}$ are bases for $V$, we define the vector $v_1 \cdots v_i \cap v_{i+1} \cdots v_{k+1} \in V$ to be the unique vector satisfying 
\begin{align}\label{eq:meetofexcess}
v_1 \cdots v_{i-1}(v_1 \cdots v_i \cap v_{i+1} \cdots v_{k+1}) =& \omega^*(v_1 \cdots v_{i-1} v_{i+1} \cdots v_{k+1}) v_1 \cdots v_i \in \bigwedge^i(V) \text{ and }\\
(v_1 \cdots v_i \cap v_{i+1} \cdots v_{k+1})v_{i+1} \cdots v_{k} =& \omega^*(v_1 \cdots v_k) v_{i+1} \cdots v_{k+1} \in \bigwedge^{k+1-i}(V).
\end{align}
This vector is unique, because it lies in the line $\Span\{v_1,\dots,v_i\} \cap \Span\{v_{i+1},\dots,v_{k+1}\}$ and is determined by (either of) the normalizations above. If $F_{\bullet}$ and $G_{\bullet}$ are affine flags that meet generically, we define the vector $F_{(i)} \cap G_{(k-i+1)}  \in V$ as $v_1 \cdots v_i \cap v_{i+1} \cdots v_{k+1}$ where $F_{(i)} = v_ 1 \cdots v_i$ and $G_{(k-i+1)} = v_{i+1} \cdots v_{k+1}$. 
\end{defn}

One can give an explicit formula for the vector $v_1 \cdots v_i \cap v_{i+1} \cdots v_{k+1}$ as a linear combination of $\{v_1,\dots,v_i\}$, and also as a linear combination of $\{v_{i+1},\dots,v_{k+1}\}$ (cf.~\cite[Equation 3.3.6 and Theorem 3.3.2a]{Sturmfels}). The renormalized vector $\tilde{w}_1$ in Remark~\ref{rmk:renormd} is a special instance of this formula. 

We now use the operation $\cap$ to define a quasi-inverse to $\Psi$. For $(F_{1,\bullet}, \dots,F_{2r,\bullet}) \in (\textnormal{Aff}(V)^{2r})^\circ$, from the first two flags, produce the $k$-tuple of vectors 
\begin{equation}\label{eq:kvectorsfromtwoflags}
F_{1,(1)}, F_{1,(2)} \cap F_{2,(k-1)}, F_{1,(3)} \cap F_{2,(k-2)}, \cdots, F_{1,(k-1)} \cap F_{2,(2)}, F_{2,(1)}. 
\end{equation}
Extending $k$-periodically and quotienting by $\SL(V)$, we obtain a map  
\begin{equation}\label{eq:Phidefn}
\Phi \colon \Conf^\circ(k,2r) \to \tGr^\circ(k,rk),   
\end{equation}
with the same implicit claim about consecutive genericity. 

\begin{thm}\label{thm:PsiPhyXY}
The pullbacks $\Psi^* \colon \bbc(\Conf^\circ(k,2r)) \to \bbc[\tGr^\circ(k,rk)]$ and $\Phi^* \colon  \bbc[\tGr(k,rk)] \to \bbc(\Conf(k,2r))$ define a quasi-isomorphism of $\bbc(\Conf^\circ(k,2r))$ and $ \bbc[\tGr^\circ(k,rk)]$. 
\end{thm}

It follows that there is a sequence of mutations $Q_k(T) \to Q(\mcc_{\Le})$ where $Q_k(T)$ is (the mutable subquiver of) a Fock-Goncharov quiver for a triangulation $T$, and $Q(\mcc_{\Le})$ is the mutable subquiver of the rectangular grid in Figure~\ref{fig:stdseed48}. Our proof of Theorem~\ref{thm:PsiPhyXY} is algebraic, i.e. we do not exhibit any such sequence of mutations.

Before proving Theorem~\ref{thm:PsiPhyXY}, we will identify a convenient nerve in the Fock-Goncharov cluster structure. First, we extend the definition of Fock-Goncharov coordinate \eqref{eq:FGCoord} in the obvious way to quadruples of flags, giving rise to functions $\Delta_{a,b,c,d}(w,x,y,z) \in \bbc(\Conf(k,r))$ for any choices of $a+b+c+d = k$ and $\{w < x < y < z \} \subset [1,r]$.  We use the same name (\emph{Fock-Goncharov coordinate}) for these more general functions on $\bbc(\Conf(k,r))$. 

\begin{lem}\label{lem:ConfNerve} Suppose $r \geq 4$. There is a nerve $\mcn$ for $\bbc(\Conf(k,r))$ on which every cluster variable on $\mcn$ is a Fock-Goncharov coordinate (involving four or fewer flags). 
\end{lem}

\begin{proof}
Let $T$ be a triangulation of the $r$-gon. Fock and Goncharov described some Pl\"ucker-like relations amongst quadruple invariants \cite[Equation 10.3]{ModuliSpaces}. Letting $\ov{a} = (a,b,c,d)$ denote a nonnegative integer solution to $a+b+c+d = k-2$, these relations are  
\begin{equation}\label{eq:flagrelation}
\Delta_{\ov{a}+(1,0,1,0)}\Delta_{\ov{a}+(0,1,0,1)} = \Delta_{\ov{a}+(1,1,0,0)}\Delta_{\ov{a}+(0,0,1,1)}+\Delta_{\ov{a}+(1,0,0,1)}\Delta_{\ov{a}+(0,1,1,0)},
\end{equation}

If a Fock-Goncharov coordinate $X \in \mathbf{x}(T)$ lies on a shared edge between two triangles in $T$, then its corresponding vertex in $Q_k(T)$ has valence four. 
The exchange relation mutating $X$ out of $\mathbf{x}(T)$ is of the form \eqref{eq:flagrelation}.

If a Fock-Goncharov coordinate $X \in \mathbf{x}(T)$ lies in a triangle of $T$, pick a quadrilateral $Q$ containing this coordinate, so that $X$ is in the interior of the quadrilateral. This is possible because $r \geq 4$. The quadrilateral $Q$ has two triangulations, one which is used in $T$. Let $T'$ be the triangulation obtained from $T$ by flipping the diagonal in the quadrilateral $Q$. Fock and Goncharov described an explicit sequence of mutations between the seeds $\Sigma(T)$ and $\Sigma(T')$ \cite[Section 10]{ModuliSpaces}. Each exchange relation in this sequence is of the form \eqref{eq:flagrelation} -- thus every cluster variable that arises during this mutation sequence is a quadruple invariant -- and every Fock-Goncharov coordinate that is in the interior of the quadrilateral is exchanged at least once. In particular, $X$ is mutated during this quadrilateral flip. Starting with an initial triangulation $T_0$, one can perform a finite sequence of quadrilateral flips so that every triangulation of the $r$-gon is visited at least once. We take as our nerve $\mcn$ all intermediate clusters that arise while performing the mutations necessary to carry out this sequence of quadrilateral flips. 
\end{proof}

\begin{proof}[Proof of Theorem~\ref{thm:PsiPhyXY}] We follow the blueprint in Lemma \ref{lem:QHnerves}. First we show that $\Psi^*$ and $\Phi^*$ pull back frozen variables to monomials in frozen variables. This also establishes the implicit claim that these maps preserve the consecutively generic loci. 

The map $\Psi$ pulls back Fock-Goncharov coordinates to Pl\"ucker coordinates. By inspection, the $\Psi^*$ pullback of a frozen Fock-Goncharov coordinate a frozen Pl\"ucker coordinate. For example, $\Psi^*(\Delta_{3,k-3}(F_{1,\bullet},F_{2,\bullet}))$ is equal to $\Delta_{1,\dots,k}$. 
On the other hand, by the defining property \eqref{eq:meetofexcess} of the vectors $F_{i,(a)}\cap F_{i+1,(k-a+1)}$, one sees that $\Phi^*$ pulls back each Grassmannian frozen variable to a product of $k-1$ Fock-Goncharov frozen variables. For example, when $r = 2, k = 4$ we have that 
\begin{align}
\Phi^*(\Delta_{1234}) &= \Delta_{1,3}(1,2)\Delta_{2,2}(1,2)\Delta_{3,1}(1,2) \\
\Phi^*(\Delta_{2345}) &= \Delta_{2,2}(1,2)\Delta_{3,1}(1,2)\Delta_{3,1}(2,3) \\
\Phi^*(\Delta_{3456}) &= \Delta_{3,1}(1,2)\Delta_{2,2}(2,3)\Delta_{1,3}(3,4) \\
\Phi^*(\Delta_{4567}) &= \Delta_{1,3}(2,3)\Delta_{1,3}(3,4)\Delta_{2,2}(3,4),
\end{align} 
and so on. 

Next, we show that $\Psi^*$ sends cluster variables to cluster variables on the nerve $\mcn$ from Lemma~\ref{lem:ConfNerve}. To begin, let $T_0$ be an initial triangulation of the $2r$-gon and $\mathbf{x}(T)$ the cluster in $\bbc(\Conf(k,2r))$. Then $\Psi^*(\mathbf{x}(T))$ is a collection of Pl\"ucker coordinates in $\bbc[\tGr^\circ(k,rk)]$. We claim this collection is weakly separated, hence forms a cluster. Indeed, for any Fock-Goncharov coordinate $X$ in a given triangle, the subset labeling the Pl\"ucker coordinate $\Psi^*(X)$ consists of at most three disjoint cyclic intervals.  If $Y$ is another coordinate in this triangle, then the cyclic intervals in $\Psi^*(X)$ and $\Psi^*(Y)$ are nested (meaning each cyclic interval in $\Psi^*(X)$ either \emph{contains} or \emph{is contained in} a cyclic interval for $\Psi^*(Y)$). This nestedness prevents the existence of $a<b<c<d$ as in the definition of weak separation. Notice that this would be false if we were considering 4 disjoint cyclic intervals. For Fock-Goncharov coordinates $X$ and $Y$ lying in different triangles, weak separation is even more clear: the three disjoint cyclic intervals will be ``far away'' from each other and can't lead to a witness.

The nerve $\mcn$ is constructed by performing a sequence of quadrilateral flips from $T_0$. During each step in this sequence, every cluster variable is a Fock-Goncharov coordinate which is mapped to a Pl\"ucker coordinate by $\Psi^*$. Applying $\Psi^*$ to the Fock-Goncharov exchange relation \eqref{eq:flagrelation} produces the corresponding $3$-term Pl\"ucker relation; thus clusters in $\mcn$ are sent to clusters in $\bbc[\tGr^\circ(k,rk)]$ in a way that is compatible with mutation. The conditions of Lemma~\ref{lem:QHnerves} are satisfied (in fact: $\Psi^*(x)$ is \emph{equal} to a cluster variable on the nerve, not merely proportional to a cluster variable), so $\Psi^*$ is a quasi-homomorphism. 

To complete the proof we need to check that $\Phi^* \circ \Psi^* = (\Psi \circ \Phi)^*$ is proportional to the identity map on $\bbc(\Conf(k,2r))$. Indeed, let $F'_{1,\bullet},F'_{2,\bullet}$ be the first two affine flags created when evaluating $\Psi \circ \Phi$ on $(F_{1,\bullet},\dots,F_{2r,\bullet})$. Using the defining property \eqref{eq:meetofexcess} repeatedly, one sees that the tensor $F'_{1,(a)}$ is a scalar multiple of the tensor $F_{1,(a)}$, and furthermore this scalar multiple is a product of frozen variables. The same relationship holds for $F'_{2,(a)}$ and $F_{2,(a)}$. 

Continuing with the example $r=2,k=4$, the affine flags $F'_1$ and $F'_2$ are given by the following sequence of simple tensors
\begin{align}
F'_1 &=  F_{1,(1)}, \omega^*(F_{1,(1)} F_{2,(3)}) F_{1,(2)}, \omega^*(F_{1,(1)} F_{2,(3)}) \omega^*(F_{1,(2)} F_{2,(2)}) F_{1,(3)}\\
F'_2 &=  F_{2,(1)}, \omega^*(F_{1,(3)} F_{2,(1)})  F_{2,(2)}, \omega^*(F_{1,(3)} F_{2,(2)})\omega^*(F_{1,(2)} F_{2,(2)})\cdot F_{2,(3)}.
\end{align}

Since the corresponding tensors only differ by scalar multiples in the frozen variables, the corresponding quadruple invariants only differ by a monomial in the frozen variables, i.e. $\Phi^* \circ \Psi^*(x) \propto x$ on the nerve.
\end{proof}


\section{Group theory results and conjectures}\label{secn:groups}
We summarize what we know about the cluster modular groups for $\Gr(k,n)$ and $\FG(k,r)$ and state some conjectures. 

After identifying the cluster structures on $\tGr(k,rk)$ and $\Conf(k,2r)$ by a quasi-isomorphism, we get an identification of their cluster modular groups $\CMG(\Gr(k,rk) \cong \CMG(\Conf(k,2r))$. It is fruitful to study the group $\CMG$ from either side of this identification. 

The well known symmetries of the Grassmannian cluster structure are the twisted cyclic shift $\rho$, the reflection symmetry $\theta$, and finally the twist map $\tau$ defined by Marsh and Scott \cite{MarshScott}. Both $\rho$ and $\theta$ are cluster automorphisms, while $\tau$ is a quasi-automorphism but not a cluster automorphism \cite[Proposition 8.10]{MarshScott}. As will be very relevant to us, the square of the twist is proportional to a power of the cyclic shift $\tau^2 \propto \rho^{-k}$ \cite[Corollary 4.2]{MarshScott}.

The well-known symmetries of the cluster structure on $\Conf(k,r)$ are the twisted cyclic shift $P$ of affine flags, the reflection symmetry $\Theta$ of affine flags, and the \emph{Hodge star map} (or \emph{duality map}) $\ast$ (cf.~\cite{GoncharovShenDT,HenriquesCactus,LeClassicalGroups} for different treatments of $\ast$). All three of these are cluster automorphisms of $\Conf(k,r)$.

\subsection{Cluster DT transformations and the center}
We recall a general method for constructing central elements of $\CMG^+(k,n)$, and thereby motivate the concept of $d$-periodicity.  

For any cluster algebra, Goncharov and Shen introduced the concept of a \emph{Donaldson-Thomas transformation} of the cluster structure \cite[Definition 3.5]{GoncharovShenDT}. When the Donaldson-Thomas transformation can be realized by mutations, it determines an element of the cluster modular group. Furthermore, this DT element is always in the center $Z(\CMG^+)$ of the group of orientation-preserving elements \cite[Corollary 3.7]{GoncharovShenDT}. In the case of the cluster structure on the space $\mca_{G,S}$ of decorated $G$-local systems (cf~Remark~\ref{rmk:AGS}), Goncharov-Shen gave an explicit construction of a DT-transformation and explained that it can be realized by mutations. For $\Conf(k,r)$, their DT-transformation is the composition $\ast \circ P^{-1} = P^{-1} \circ \ast$ \cite[Theorem 1.3]{GoncharovShenDT}. 

Marsh and Scott \cite[Theorem 11.17]{MarshScott} showed that the twist map could be implemented by particular type of mutation sequence known as a \emph{maximal green sequence}, which implies that the Marsh-Scott twist map is a DT-transformation of $\Gr(k,n)$. 
In particular, the Marsh-Scott twist map is always central in $\CMG^+(\Gr(k,n))$. Since $\tau^2 = \rho^{-k}$, it follows that every orientation-preserving quasi-automorphism commutes with $\rho^k$, and hence, with $\rho^d$. This motivates the notion of $d$-periodicity in Section~\ref{secn:BraidGroupAction}. Notice also that when $\gcd(k,n) = 1$, it follows that any quasi-automorphism of $\tGr(k,n)$ commutes with $\rho$ up to proportionality. This is a very restrictive requirement, and we view it as a suggestion that these particular Grassmannians should have a small cluster modular group.

Intertwining with the maps \eqref{eq:Psidefn} and \eqref{eq:Phidefn}, symmetries of Fock-Goncharov spaces become interesting symmetries of Grassmannians. 

\begin{prop}\label{prop:compnisP} Let $f_\Delta$ be the quasi-automorphism of $\tGr(k,rk)$ corresponding to the half-twist braid $\Delta$. Then $f_\Delta  \propto (\Phi \circ P \circ \Psi)^*$. That is, the twisted cyclic shift of affine flags in $\FG(k,2r)$ intertwines to the half-twist braid on the Grassmannian. 

The Donaldson-Thomas transformation on $\FG(k,2r)$ intertwines to the  Marsh-Scott twist map $\tau$ on $\tGr(k,rk)$. In fact, $\tau= \Phi \circ (P^{-1} \circ \ast) \circ \Psi$, not merely $\tau \propto \Phi \circ (P^{-1} \circ \ast) \circ \Psi$. 
\end{prop}

Note that the relation $\Delta^2 = \rho^k \in B_k$ reflects the $k$-periodicity of the maps \eqref{eq:Psidefn} and \eqref{eq:Phidefn}. 

Our proof of the formula $\tau= \Phi \circ (P^{-1} \circ \ast) \circ \Psi$ uses a well-known relationship between the Hodge star map $\ast \colon \bigwedge^a(V) \to \bigwedge^{k-a}(V)$ and the $\cap$ operation from Definition~\ref{defn:weakmeet}. First, we extend $\cap$ bilinearly to a map $\bigwedge^{a}(V) \otimes \bigwedge^{k-a+1}(V) \to V$. Then $\ast$ interchanges the exterior product for the $\cap$ operation, i.e. 
$\ast(vw) = \ast(v) \cap \ast(w)$ for tensors $v \in \bigwedge^a(V), w \in \bigwedge^{k-a-1}(V)$ (cf., for example, \cite[Theorem 6.3]{RotaHodge}).

\begin{proof}
The window notation for the map $\Phi \circ P \circ \Psi$ is
\begin{equation}\label{eq:phippsi}
[v_k, v_{k-1}v_k \cap v_{k+1} \cdots v_{2k-1},v_{k-2}v_{k-1}v_k \cap v_{k+1} \cdots v_{2k-2}, \cdots, v_{2}\cdots v_{k} \cap v_{k+1} v_{k+2}, v_{k+1}].
\end{equation}

To evaluate the half-twist braid $f_\Delta$, we break the computation into steps, by first evaluating $\sigma_{k-1} \cdots \sigma_1$, then evaluating $\sigma_{k-2} \cdots \sigma_1$, and so on, ending with $\sigma_1$. In the proof of Lemma~\ref{lem:consecutivecomp}, we we wrote down \eqref{eq:321} the window notation for each of these compositions $\sigma_i \cdots \sigma_1$ (in the case $k=5$). Continuing with the case $k=5$ for concreteness, the arguments in the proof of 
Lemma~\ref{lem:consecutivecomp} established that $a_1 \in \Span\{v_1,v_2\} \cap \Span\{v_3,v_4,v_5,v_6\}$, $b_1 \in \Span\{v_1,v_2,v_3\} \cap \Span\{v_4,v_5,v_6\}$, and $c_1 \in \Span\{v_1,v_2,v_3,v_4\} \cap \Span\{v_5,v_6\}$, where $a_1,b_1,c_1$ are the vectors from \eqref{eq:321}. 
Thus $a_1$ differs from the vector $v_1v_2 \cap v_3v_4v_5v_6$ by a scalar multiple, and we claim that this scalar multiple is a ratio of two frozen variables. Indeed, from the definition of $\cap$ one has $v_2(v_1v_2 \cap v_3v_4  v_5v_6)v_3v_4v_5 = \omega^*(v_1\cdots v_5)v_2 \cdots v_6$. On the other hand, 
$v_2a_1v_3v_4v_5 = v_1 \cdots v_5$.  So $a_1 = \frac{1}{\omega^*(v_2 \cdots v_6)} v_1v_2 \cap v_3 \cdots v_6$ as claimed. By a similar style of argument, one see that $b_1$ agrees with $v_1v_2v_3 \cap v_4v_5v_6$ up to a frozen Laurent monomial, and likewise for $c_1$ and $v_1v_2v_3v_4 \cap v_5v_6$. Now we plug each of these simplifications into the corresponding window notation from \eqref{eq:321}, while dropping frozen monomial factors, to see that $f_\Delta$ is proportional to a map with window notation
\begin{align*}
[v_1,\dots,v_5] & \overset{\sigma_4 \sigma_3 \sigma_2\sigma_1}{\mapsto} [v_2,v_3,v_4,v_5,v_6] \overset{\sigma_3 \sigma_2\sigma_1}{\mapsto}[v_3,v_4,v_5,v_2v_3v_4v_5 \cap v_6v_7,v_6] \\
& \overset{\sigma_2\sigma_1}{\mapsto}  [v_4,v_5, v_3v_4v_5 \cap v_6v_7v_8 ,v_2v_3v_4v_5 \cap v_6v_7,v_6]  \\
&\overset{\sigma_1}{\mapsto}  [v_5,v_4v_5 \cap v_6 v_7v_7v_9, v_3v_4v_5 \cap v_6v_7v_8 ,v_2v_3v_4v_5 \cap v_6v_7,v_6].
\end{align*}
Comparing with \eqref{eq:phippsi} establishes the claim when $k=5$. The general argument is similar.

As for the statement about the twist: the composition $\Phi \circ P^{-1} \circ \ast \circ \Psi$ has window notation 
\begin{equation}\label{listofstuff}
[\ast (v_{n-k+2} \cdots v_n ), \ast \left(v_{n-k+3} \cdots v_n \right) \cap \ast (v_{1}),\dots,  \ast (v_n)  \cap \ast (v_{1}\cdots v_{k-2} ), \ast(v_1 \cdots v_{k-1})] 
\end{equation}
Recalling that $\ast$ interchanges the exterior product and $\cap$, one simplifies $\ast(v_{n-k+3} \cdots v_n) \cap \ast(v_{1}) = \ast(v_{n-k+3} \cdots v_n v_{1})$, and so on. The $i$th term in the window notation becomes $\ast(v_{i-k+1} \cdots v_{i -1})$, which agrees with the definition \cite[Definition 2.1]{MarshScott} of the (left) Marsh-Scott twist. 
\end{proof}

\subsection{Conjectures}
To treat the cases $d=k$ and $d<k$ uniformly, we set $G' = \hat{B}_{\hat{A}_{d-1}}$ when $d<k$ and $G' = B_k$ when $d=k$. In either case, we have a group homomorphism $G' \to \CMG^+(\Gr(k,n))$. In the latter case $d=k$, we use the symbol $\rho$ for the braid $\sigma_{k-1} \cdots \sigma_1 \in B_k$, remembering that this braid acts by the cyclic shift in $\CMG^+$.  

\begin{conj}\label{conj:groups} The orientation-preserving cluster modular group $\CMG^+(\tGr(k,n))$ is generated by the twist map, the cyclic shift, and the Artin generators $\sigma_1,\dots,\sigma_{d-1}$. 
For fixed $k$ and $n >>k$,  the kernel of the map $G' \to \CMG^+(\Gr(k,n))$ is generated by $\rho^n$. 
\end{conj}

The hypothesis $n>>k$ is necessary, because we are aware of extra elements in the kernel of the map $G' \to \CMG^+(\Gr(k,2k))$ (see the next section). The conjecture 
is fairly optimistic, because our only evidence is that it is compatible with what is known for the finite type Grassmannians (cf.~Remark~\ref{rmk:finitetype}), and with what we prove for the other finite mutation type Grassmannians $\Gr(3,9)$ and $\Gr(4,8)$ in Section~\ref{secn:Webs}. In the case that $d = \gcd(k,n) = 1$, then no Artin generators are defined, and the conjecture says that $\CMG^+(\Gr(k,n))$ is a finite cyclic group generated by the twist. The conjecture could be disproved in these cases by finding an element of the cluster modular group that is not a power of the twist. 

When $k=2$, one has $\tau = \rho^{-1}$. When $k>2$, we expect that the twist is not proportional to any element of the extended affine braid group. 

\begin{rmk}[Groups for finite type Grassmannians]\label{rmk:finitetype}
In \cite{ASS}, the authors computed the cluster modular groups for acyclic cluster algebras, and in particular for all finite or affine Dynkin types. The Grassmannians $\Gr(2,n)$, which have cluster type $A_{n-3}$, satisfy $\CMG^+(\Gr(2,n)) \cong \bbz / n \bbz $. 
The cluster structures on $\Gr(3,7),\FG(3,5)$,  and $\Gr(3,8)$ have finite Dynkin types $E_6,E_7,E_8$. One has $\CMG^+(\Gr(3,7)) \cong \bbz / 14 \bbz$, $\CMG^+(\FG(3,5)) \cong \bbz / 10 \bbz$, and $\CMG(\Gr(3,8)) \cong \bbz / 16 \bbz$, with generator the Donaldson-Thomas transformation in each case.  
The only other finite type Grassmannian or Fock-Goncharov space is $\bbc[\Gr(3,6)] \sim \bbc[\FG(3,4)]$ which has Dynkin type $D_4$. In this case $\CMG^+(3,6) \cong \mfs_3 \times \bbz / 4 \bbz$ is larger than the group $\la \rho,\tau\ra \cong \bbz_3 \times \bbz_4$, with the extra elements coming from braids.   
\end{rmk}

\subsection{The case $n=2k$} We explain why we the hypothesis $n>>k$ is needed in Conjecture~\ref{conj:groups}. 
\begin{prop}\label{prop:upandown} The pullback of the composition $\sigma_1 \cdots \sigma_{k-1}^2 \cdots \sigma_1$ is proportional to the identity map 
on $\bbc[\tGr^\circ(k,2k)]$. %
\end{prop}

\begin{rmk} The quotient $B_k / \la \sigma_1 \cdots \sigma_{k-1}^2 \cdots \sigma_1 \ra$ is the \emph{spherical braid group}, defined as the fundamental group $\pi_1(\text{Conf}(S^2,k))$ of the configuration space of $k$ unlabeled points on the 2-sphere. Its center $Z(\pi_1(\text{Conf}(S^2,k)))$ is generated by the image of the full-twist $\Delta^2 \in B_k$. This element $\Delta$ has order $4$ in $\pi_1(\text{Conf}(S^2,k))$, so that $Z(\pi_1(\text{Conf}(S^2,k))) \cong \bbz / 2 \bbz$. The quotient of the spherical braid group by this center is the mapping class group $\MCG(S^2,k)$ of a sphere with $k$ punctures. We add to Conjecture~\ref{conj:groups} by conjecturing that $\CMG^+(\Gr(k,2k)) / Z(\CMG^+(\Gr(k,2k))) \cong \MCG(S^2,k)$, i.e. that the relation $\sigma_1 \cdots \sigma_{k-1}^2 \cdots \sigma_1$ is the only surprising relation in the cluster modular group. We speculate that the assumption $n>>k$ in Conjecture~\ref{conj:groups} can be replaced by $n\geq 3k$.  \end{rmk}

\begin{rmk}\label{rmk:exotic} The relation $\sigma_1 \cdots \sigma_{k-1}^2  \cdots \sigma_1 \propto \text{id} \in \bbc[\tGr(k,2k)]$ is exotic in a sense we explain now. 
Recall the (right) torus action $(V^n)^\circ \curvearrowleft (\bbc^*)^n$ by rescaling the $v_i$'s. Suppose that $f$ and $g$ are a pair of regular maps on $(V^n)^\circ$.  
Then the typical explanation for the statement that $f^* \propto g^*$ as maps on $\tGr^\circ(k,n)$ is to observe that each of the coordinate functions of $f$ and $g$ agree up to multiplication by a scalar function, and this scalar function is a Laurent monomial in the frozens. This was how we observed the proportionality statements in Lemma~\ref{lem:consecutivecomp} and Proposition~\ref{prop:compnisP}. This is simply not true for the present relation: $\sigma_1 \cdots \sigma_{k-1}^2  \cdots \sigma_1$ is not the identity map on $(V^{2k})^\circ / (\bbc^*)^{2k}$, but it becomes the identity map after quotienting by the left $\SL(V)$ action. For an explicit calculation, take the $4 \times 8$ matrix $M = \begin{pmatrix}
     1   &   0  &    0   &   0  &   -1   &  -3  &-\frac{19}{2} &  -\frac{151}{8} \\
     0   &   1  &    0   &   0  &   1    &  2   &   5          & \frac{19}{2} \\  
     0   &   0  &    1   &   0  &   -2   &  -2  &   -2         &-3 \\
     0   &   0  &    0   &   1  &    3   &   2  &    1         &1,
\end{pmatrix}$ for whom 
$\sigma_1\sigma_2\sigma_3^2\sigma_2\sigma_1(M) = 
\begin{pmatrix}
   -1 & -\frac{11}{4} &    -3  &-\frac{7}{4}    &-1   &  0  &   0 &    0 \\
     1 &  \frac{7}{4}   &  3   & \frac{7}{4}     &0    & -2  &  -5 & -\frac{19}{2} \\
   -2 & -\frac{11}{2}  &  -7  &  -\frac{7}{2}   &  0  &   2 &    2 &    3\\
    3 & \frac{33}{4}   &  9  &  \frac{17}{4}    & 0   & -2  &  -1  &  -1.
\end{pmatrix}$ (we omit the calculation). All frozen Pl\"ucker coordinates of $M$, and hence of $\sigma_1\sigma_2\sigma_3^2\sigma_2\sigma_1(M)$, equal one. All Pl\"ucker coordinates of these two matrices agree (in agreement with Proposition~\ref{prop:upandown}). 
On the other hand, these matrices clearly are different elements of $(V^\circ)^8 / (\bbc^*)^8$. We expect that relations between Artin generators of this sort can only happen for small $n$. 
\end{rmk}

\begin{proof}[Proof of Proposition \ref{prop:upandown}]
To show that $\sigma_1 \cdots \sigma_{k-1}^2 \cdots \sigma_1$ acts trivially up to frozens, by Lemma~\ref{lem:consecutivecomp} we can equivalently show that 
$\sigma_1 \circ \cdots \circ \sigma_{k-1}$ is proportional to $\rho^{-1}$. 

By similar calculations as in Lemma~\ref{lem:consecutivecomp} and Proposition~\ref{prop:compnisP}, one sees that
$\sigma_1 \circ \cdots \circ \sigma_{k-1}$ has window notation
\begin{equation}\label{eq:upandown}
\propto [v_k,v_1v_k \cap v_{k+1} \cdots v_{2k-1},v_2v_k \cap v_{k+1} \cdots v_{2k-1},\dots,v_{k-1}v_k \cap v_{k+1} \cdots v_{2k-1}].
\end{equation}
These two windows do not resemble the windows for $\rho^{-1}$. To show these two maps are proportional when acting on $\tGr(k,2k)$, 
we evaluate both on the $\Le$ diagram cluster $\mcc_{\Le}$ \eqref{eq:stdseed}. For each non-frozen Pl\"ucker coordinate $\Delta_S \in \mcc_{\Le}$, we explicitly verify $(\sigma_1 \circ \cdots \circ \sigma_{k-1})^*(\Delta_S) \propto \rho^{-1}(\Delta_S)$ using the window notation \eqref{eq:upandown}. This establishes the claim because a quasi-automorphism is determined by where it sends a cluster. The reader may wish to check our calculations in $\Gr(4,8)$ (cf.~Figure~\ref{fig:stdseed48}).
 
Let $S = [1,a] \cup [b+1,b+k-a)]$ be given. We want to prove that if we take the wedge of the first $a$ vectors in \eqref{eq:upandown}, and then take the wedge of the vectors in locations $b+1,\dots,b+k-a$, we obtain an element of $\bbp$ times the tensor $v_1 \cdots v_{a-1}v_{b}\cdots v_{b+k-a-1}v_{2k}$.
Since $1 \in S$, the vector $v_k$ is certainly in our outcome. In the presence of this vector, the defining property of $\cap$ \eqref{eq:meetofexcess} allows us to simplify $v_k (v_i v_k\cap v_{k+1} \cdots v_{2k-1})$ as $v_iv_k$ times a frozen variable, and we ignore this frozen variable. Thus, we can summarize the exterior product of the vectors in columns $S \cap [1,k]$ of \eqref{eq:upandown} as follows: $1\in S$ becomes $k \in (\sigma_1 \circ \cdots \circ \sigma_{k-1})^*(\Delta_S)$, and $i \in S \cap [2,k]$ becomes $i-1 \in (\sigma_1 \circ \cdots \circ \sigma_{k-1})^*(\Delta_S)$. Since $\Delta_S$ is non-frozen, $S \setminus [1,k]$ is nonempty.
There are two cases, either $b \leq k$ (thus $k+1 \in S$) or $b >k$ (thus $k+1 \notin S$). The first case is easier: applying the defining property of $\cap$ again on this second window, $k+1 \in S$ becomes $2k \in (\sigma_1 \circ \cdots \circ \sigma_{k-1})^*(\Delta_S)$ and $i \in S \cap [k+2,2k]$ becomes $i-1 \in (\sigma_1 \circ \cdots \circ \sigma_{k-1})^*(\Delta_S)$. Altogether, the union of the calculations in the two windows shows that $(\sigma_1 \circ \cdots \circ \sigma_{k-1})^*(\Delta_S) \propto \rho^{-1}(\Delta_{S})$ in this first case. 

The second case is quite subtle. The wedge of the $a$ vectors in the first window produces the tensor $v = v_kv_1 \cdots v_{a-1}$. We denote by $w \in \bigwedge^{k-a}(V)$ the wedge of the vectors in locations $b+1,\dots,b+k-a$ of (the second window of) \eqref{eq:upandown}. We denote by $u_i = v_iv_{2k} \cap v_1 \cdots v_{k-1}$ the vector in the $i+1$th column. 

Each such $u_i$ can be expressed as a linear combination of the vectors $v_1,\dots,v_{k-1}$, and consequently we can simplify 
\begin{align*}
vw &= (\text{ coefficient of $v_a\cdots v_{k-1}$ in $w$}) \cdot v_1 \cdots v_k \in \bigwedge^k(V) \\ 
   &= (\text{ coefficient of $v_{2k}v_1 \cdots v_{k-1}$ in $v_{2k}v_1 \cdots v_{a-1}w$}) \cdot v_1 \cdots v_k
\end{align*}
Both of the tensors $v_{2k}v_1 \cdots v_{k-1}$ and $v_{2k}v_1 \cdots v_{a-1}w$ are in the top exterior power $\bigwedge^k(V)$, so the coefficient in the second line above is just the scalar multiple relating them. Using the defining property of $\cap$, we can simplify $v_{2k}u_i = \omega^*(v_{2k} v_1 \cdots v_{k-1}) v_i$ which is a frozen variable times $v_i$. Ignoring this frozen variable factor and repeating this simplification, we can replace each $u_i$ in the tensor by its corresponding $v_i$, which has the effect of performing $\rho^{-1}$ on the corresponding column. When the dust clears, the tensor $vw$ above factors as a product of frozen variables times $\omega^*(\rho^{-1}(v_1,\dots,v_{2k}))$ times the tensor $v_1 \cdots v_k$. Applying $\omega^*$ to this equality, the latter tensor becomes a frozen variable, and the result follows. 
\end{proof}

\begin{rmk}
The Grassmannians $\Gr(k,2k)$ bear a cluster automorphism induced by the complementation map on Pl\"ucker coordinates.  By Conjecture~\ref{conj:groups}, this element should be expressible as a composition of our conjectural generators. We calculated that the complementation map is proportional to the composition $ \theta \circ \tau \circ \sigma_1 \circ \rho^{-2}$. We omit the details. 
\end{rmk}

\section{Webs and the Fomin-Pylyavskyy conjectures}\label{secn:Webs}
Fomin and Pylyavskyy \cite{tensors,tensorsII} proposed a description of the cluster combinatorics for $\bbc[\tGr(3,n)]$ (and in fact, for a wider class of cluster algebras associated with algebras of $\SL_3$ invariants) in terms of Kuperberg's basis of \emph{non-elliptic} webs. We review this combinatorics, and related versions for $\FG(3,r)$ and $\Gr(4,n)$, in the following sections. 

\subsection{Fomin-Pylyavskyy conjectures}
\begin{defn}\label{defn:tensordiagram}
Let $D$ be a disk with $n$ marked points $1,\dots,n$ labeled clockwise on its boundary. A \emph{tensor diagram} is a finite bipartite graph drawn in $D$, with a fixed bipartition of its vertex set into black and white color sets, subject to the following additional requirements:
\begin{itemize}
\item the \emph{boundary vertices} are in the vertex set and are colored black, 
\item the remaining \emph{interior vertices} are in the interior of the disk and are trivalent. 
\end{itemize}
\end{defn}

Such a tensor diagram is considered up to isotopy of the disk fixing the boundary. The edges in a tensor diagram are allowed to intersect transversely. We emphasize that the boundary vertices can have arbitrary valence, including valence zero. 

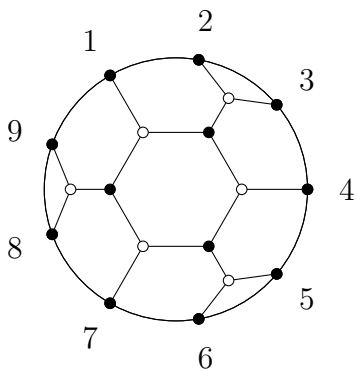
\begin{figure}[ht]
\begin{center}
\begin{tikzpicture}[scale = 1.75]
\def \nnn{9};
\def  \tic{360/9};
\draw [black] (1,0) arc [radius = 1, start angle = 0, end angle = 360];
\foreach \s in {1,...,\nnn}
{
  \node at ({\tic* (-\s+4)}:1.3 cm) {$\s$};
  \draw [fill= black] ({\tic * (\s - 1)}:1 cm) circle [radius = .04];
}

\coordinate (A) at (1.5*\tic:.5);
\coordinate (B) at (3*\tic:.5);
\coordinate (C) at (4.5*\tic:.5);
\coordinate (D) at (6*\tic:.5);
\coordinate (E) at (7.5*\tic:.5);
\coordinate (F) at (9*\tic:.5);
\coordinate (G) at (1.5*\tic:.8);
\coordinate (H) at (4.5*\tic:.8);
\coordinate (I) at (7.5*\tic:.8);

\draw (A)--(B)--(C)--(D)--(E)--(F)--(A);
\draw (A)--(G);
\draw (C)--(H);
\draw (E)--(I);
\draw (B)--(3*\tic:1);
\draw (D)--(6*\tic:1);
\draw (F)--(9*\tic:1);
\draw (1*\tic:1)--(G)--(2*\tic:1);
\draw (4*\tic:1)--(H)--(5*\tic:1);
\draw (7*\tic:1)--(I)--(8*\tic:1);

\draw [fill= white] (G) circle [radius = .04];
\draw [fill= white] (H) circle [radius = .04];
\draw [fill= white] (I) circle [radius = .04];
\draw [fill= white] (B) circle [radius = .04];
\draw [fill= white] (D) circle [radius = .04];
\draw [fill= white] (F) circle [radius = .04];
\draw [fill= black] (A) circle [radius = .04];
\draw [fill= black] (C) circle [radius = .04];
\draw [fill= black] (E) circle [radius = .04];

\begin{scope}
\draw [black] (1,0) arc [radius = 1, start angle = 0, end angle = 360];
\foreach \s in {1,...,\nnn}
{
  \draw [fill= black] ({\tic * (\s - 1)}:1 cm) circle [radius = .04];
}

\end{scope}

\end{tikzpicture}
\caption{A tensor diagram for $\Gr(3,9)$. \label{fig:SCW}}
\end{center}
\end{figure}

Every tensor diagram~$T$ with $n$ boundary vertices defines an \emph{invariant} $[T] \in \bbc[\tGr(3,n)]$. We will define this invariant in a somewhat indirect way via the $\Gr(3,n)$ skein algebra, and then explain how the skein algebra can be identified with $\bbc[\tGr(3,n)]$. The validity of our description is based on \cite{tensors, Kuperberg}. 

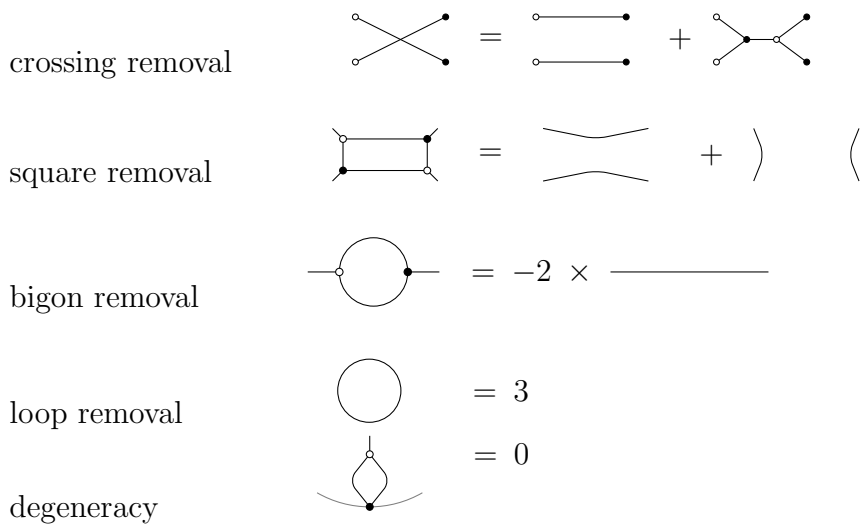
\begin{figure}[ht]
\begin{center}
\begin{tabular}{lc}
crossing removal &\begin{tikzpicture}[scale = .6]
\coordinate (A) at (0,1);
\coordinate (B) at (2,1);
\coordinate (C) at (2,0);
\coordinate (D) at (0,0);
\draw (A)--(C);
\draw (B)--(D);

\draw [fill= white] (A) circle [radius = .065];
\draw [fill= white] (D) circle [radius = .065];
\draw [fill= black] (B) circle [radius = .065];
\draw [fill= black] (C) circle [radius = .065];
\node at (3,.5) {$=$};
\node at (7.2,.5) {$+$};

\begin{scope}[xshift = 4cm]
\coordinate (AA) at (0,1);
\coordinate (BB) at (2,1);
\coordinate (CC) at (2,0);
\coordinate (DD) at (0,0);
\draw (AA)--(BB);
\draw (CC)--(DD);

\draw [fill= white] (AA) circle [radius = .065];
\draw [fill= white] (DD) circle [radius = .065];
\draw [fill= black] (BB) circle [radius = .065];
\draw [fill= black] (CC) circle [radius = .065];
\end{scope}

\begin{scope}[xshift = 8cm]
\coordinate (AAA) at (0,1);
\coordinate (BBB) at (2,1);
\coordinate (CCC) at (2,0);
\coordinate (DDD) at (0,0);
\coordinate (E) at (.67,.5);
\coordinate (F) at (1.33,.5);
\draw (AAA)--(E)--(F)--(BBB);
\draw (CCC)--(F);
\draw (DDD)--(E);
\draw [fill= white] (AAA) circle [radius = .065];
\draw [fill= white] (DDD) circle [radius = .065];
\draw [fill= black] (BBB) circle [radius = .065];
\draw [fill= black] (CCC) circle [radius = .065];
\draw [fill= black] (E) circle [radius = .065];
\draw [fill= white] (F) circle [radius = .065];
\node at (2.5,0) {\hfill};

\end{scope}
\end{tikzpicture} \\
\vspace{.05cm} & \vspace{.05cm} \\ 
square removal& 
\begin{tikzpicture}[scale = .7]
\coordinate (A) at (.2,.8);
\coordinate (B) at (1.8,.8);
\coordinate (C) at (1.8,.2);
\coordinate (D) at (.2,.2);
\draw (A)--(B)--(C)--(D)--(A);
\draw (A)--(0,1);
\draw (B)--(2,1);
\draw (C)--(2,0);
\draw (D)--(0,0);

\draw [fill= white] (A) circle [radius = .065];
\draw [fill= white] (C) circle [radius = .065];
\draw [fill= black] (B) circle [radius = .065];
\draw [fill= black] (D) circle [radius = .065];
\node at (3,.5) {$=$};
\node at (7.2,.5) {$+$};

\begin{scope}[xshift = 4cm]
\draw [rounded corners] (0,1)--(1,.8)--(2,1);
\draw [rounded corners] (0,0)--(1,.2)--(2,0);
\end{scope}

\begin{scope}[xshift = 8cm]
\draw [rounded corners] (0,1)--(.2,.5)--(0,0);
\draw [rounded corners] (2,1)--(1.8,.5)--(2,0);
\end{scope}
\end{tikzpicture} \\
\vspace{.05cm} & \vspace{.05cm} \\ 
bigon removal & 
\begin{tikzpicture}[scale = .7]
\draw  (0,-2) circle [radius = .65];
\draw (-1.25,-2)--(-.65,-2);
\draw (1.25,-2)--(.65,-2);
\draw [fill= white] (-.65,-2) circle [radius = .07];
\draw [fill= black] (.65,-2) circle [radius = .07];
\node at (3,-2) {$= \hspace{.1cm} -2 \hspace{.2cm} \times$}; 
\draw (4.5,-2)--(7.5,-2);
\node at (9.5,-2) {\hfill};
\end{tikzpicture} \\
\vspace{.025cm} & \vspace{.025cm} \\ 
loop removal & 
\begin{tikzpicture}[scale = .7]
\draw (.5,0) circle [radius = .6];
\node at (-1.4,0) {\hfill};
\node at (3,0) {$= \hspace{.1cm} 3$}; 
\node at (11,0) {\hfill};
\end{tikzpicture} \\
degeneracy& 
\begin{tikzpicture}[scale = .7]
\draw [gray] (0,0) arc [radius=2, start angle=-90, end angle= -60];
\draw [gray] (0,0) arc [radius=2, start angle=-90, end angle= -120];
\draw [rounded corners] (0,0)--(-.4,.5)--(0,1);
\draw [rounded corners] (0,0)--(.4,.5)--(0,1);
\draw (0,1)--(0,1.35);
\draw [fill= black] (0,0) circle [radius = .065];
\draw [fill= white] (0,1) circle [radius = .065];
\node at (2.5,1) {$= \hspace{.1cm} 0$};
\node at (9.4,0) {\hfill};
\end{tikzpicture}
\end{tabular}
\end{center}
\caption{Skein relations amongst $\SL_3$ web invariants. The degeneracy relation only applies at boundary vertices.  \label{fig:SL3skeinrelations}}
\end{figure}

\begin{defn}[Skein algebra]\label{defn:skeinalgebra}
Consider the space of $\bbc$-linear combinations of tensor diagrams for $\tGr(3,n)$. The $\Gr(3,n)$ \emph{skein algebra} is the quotient of this algebra by the \emph{skein relations} in Figure~\ref{fig:SL3skeinrelations}. If $T$ is a tensor diagram with $n$ boundary vertices, the \emph{invariant}~$[T]$ associated to $T$ is its image in the $\Gr(3,n)$ skein algebra. If $T$ and $T'$ are tensor diagrams, then their product $[T][T']$ in the skein algebra is the tensor diagram obtained by superimposing the two tensor diagrams (one on top of the other). 
\end{defn}

Each skein relation is a local diagrammatic rule: if one locates a local fragment in a tensor diagram that looks like the left hand side of a skein relation, one can replace it with the corresponding fragment on the right hand side. If the right hand side is a sum of two fragments, then applying the skein relation produces two tensor diagrams, one in which the first fragment is used and another in which the second fragment is used. We illustrate the use of skein relations in Example~\ref{eg:TwoTripods}.

There will often be several (inisotopic) ways of superimposing two tensor diagrams as in Definition~\ref{defn:skeinalgebra}, but it is a fact that the resulting element of the skein algebra is well-defined (this is a consequence of Kuperberg's theorem below, but there are more intrinsic explanations).   

We define a map from the skein algebra to $\bbc[\Gr(3,n)]$ by mapping the tripod joining boundary vertices $i,j,k$ to the Pl\"ucker coordinate $\Delta_{ijk} \in \bbc[\Gr(3,n)]$. One can deduce the skein relations from the Pl\"ucker relations in $\bbc[\Gr(3,n)]$, and vice versa, so this map identifies the skein algebra with the coordinate ring (note the map is surjective since Pl\"ucker coordinates generate the coordinate ring). 

\begin{defn}[Web invariants]\label{defn:Web} A \emph{web} is a planar tensor diagram. A web is \emph{non-elliptic} if it contains no $2$-cycles based at a boundary vertex,  and if all of its faces formed by interior vertices are bounded by at least six sides. A \emph{web invariant} is an element $[W] \in \bbc[\tGr(3,n)]$ for a non-elliptic web $W$. Two web invariants are \emph{compatible} if their product is again a web invariant. A web invariant is \emph{indecomposable} if it does not factor as a product of web invariants. 
\end{defn}

By repeatedly applying the skein relations, it is easy to see that web invariants span the skein algebra. Kuperberg \cite{Kuperberg} established that in fact, the web invariants are a \emph{basis} for $\bbc[\tGr(3,n)]$. In particular, any web invariant is equal to $[W]$ for a unique non-elliptic web, so we often identify non-elliptic webs with their invariants.

\begin{example}\label{eg:TwoTripods} Consider the Pl\"ucker coordinates $\Delta_{124}, \Delta_{135} \in \bbc[\tGr(3,5)]$. We represent their product~$\Delta_{124}\Delta_{135}$ as a tensor diagram by superimposing two tripods. Figure~\ref{fig:expressproduct} expands this product as a sum of non-elliptic webs via the crossing removal relation. 
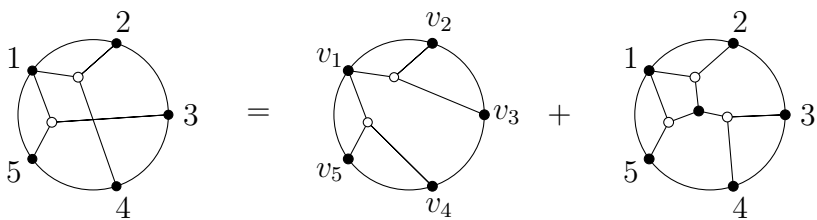
\begin{figure}
\begin{tikzpicture}[scale = 1]
\def \nnn{5};
\def  \tic{360/5};
\node at (2.2,0) {$=$};
\node at (6.2,0) {$+$};
\draw [black] (1,0) arc [radius = 1, start angle = 0, end angle = 360];
\foreach \s in {1,...,\nnn}
{
  \node at ({\tic* (-\s+3)}:1.3 cm) {$\s$};
  \draw [fill= black] ({\tic * (\s - 1)}:1 cm) circle [radius = .06];
}

\coordinate (A) at (2*\tic:1);
\coordinate (B) at (1*\tic:1);
\coordinate (C) at (0:1);
\coordinate (D) at (-1*\tic:1);
\coordinate (E) at (-2*\tic:1);
\coordinate (F) at (-.2,.5);
\coordinate (G) at (-.55,-.1);
\draw (A)--(F)--(B)--(F)--(D);
\draw (A)--(G)--(C)--(G)--(E);
\draw [fill= white] (F) circle [radius = .065];
\draw [fill= white] (G) circle [radius = .065];

\begin{scope}[xshift = 4.2cm]
\draw [black] (1,0) arc [radius = 1, start angle = 0, end angle = 360];
\foreach \s in {1,...,\nnn}
{
  \node at ({\tic* (-\s+3)}:1.3 cm) {$v_{\s}$};
  \draw [fill= black] ({\tic * (\s - 1)}:1 cm) circle [radius = .065];
}

\coordinate (AA) at (2*\tic:1);
\coordinate (BB) at (1*\tic:1);
\coordinate (CC) at (0:1);
\coordinate (DD) at (-1*\tic:1);
\coordinate (EE) at (-2*\tic:1);
\coordinate (FF) at (-.2,.5);
\coordinate (GG) at (-.55,-.1);
\draw (AA)--(FF)--(BB)--(FF)--(CC);
\draw (AA)--(GG)--(DD)--(GG)--(EE);
\draw [fill= white] (FF) circle [radius = .06];
\draw [fill= white] (GG) circle [radius = .06];
\end{scope}

\begin{scope}[xshift = 8.2cm]
\draw [black] (1,0) arc [radius = 1, start angle = 0, end angle = 360];
\foreach \s in {1,...,\nnn}
{
  \node at ({\tic* (-\s+3)}:1.3 cm) {$\s$};
  \draw [fill= black] ({\tic * (\s - 1)}:1 cm) circle [radius = .065];
}

\coordinate (AAA) at (2*\tic:1);
\coordinate (BBB) at (1*\tic:1);
\coordinate (CCC) at (0:1);
\coordinate (DDD) at (-1*\tic:1);
\coordinate (EEE) at (-2*\tic:1);
\coordinate (FFF) at (-.2,.5);
\coordinate (GGG) at (-.55,-.1);
\coordinate (H) at (-.15,.05);
\coordinate (I) at (.23,-.03);
\draw (AAA)--(FFF)--(H)--(GGG)--(AAA);
\draw (FFF)--(BBB);
\draw (GGG)--(EEE);
\draw (H)--(I)--(CCC)--(I)--(DDD);
\draw [fill= white] (FFF) circle [radius = .065];
\draw [fill= white] (GGG) circle [radius = .065];
\draw [fill= black] (H) circle [radius = .065];
\draw [fill= white] (I) circle [radius = .065];
\end{scope}

\end{tikzpicture}
\caption{An application of the crossing removal skein relation. \label{fig:expressproduct}}
\end{figure}
\end{example}

\begin{rmk}\label{rmk:invariantdefn} 
The invariant $[T]$ associated to a tensor diagram~$T$ has an intrinsic definition as a repeated contraction of certain basic $\SL_3$-invariant tensors, namely the volume form $V^{\otimes 3 } \to \bbc$, the dual form $(V^*)^{\otimes 3 } \to \bbc$,  and the pairing $V \otimes V^* \to \bbc$ (cf.~\cite[Section 4]{tensors}). These correspond to an internal white vertex, an internal black vertex, and an edge connecting a black vertex to a white vertex, respectively.  One can also give an explicit expression \cite[Equation (4.1)]{tensors} for the invariant~$[T]$ has  as an $\SL(V)$-invariant polynomial in the coordinates on $V^n$. 
\end{rmk}

Fomin and Pylyavskyy conjectured that every cluster variable in $\bbc[\tGr(3,n)]$ is a web invariant. They also proposed the following combinatorial procedure for determining which web invariants are cluster variables.  

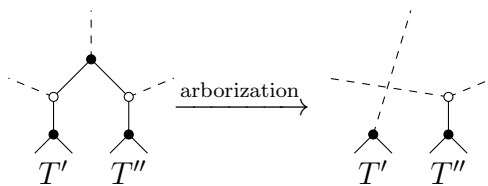
\begin{figure}[ht]
\begin{center}
\begin{tikzpicture}[scale = .5]
\coordinate (A) at (0,0);
\coordinate (B) at (1,-1);
\coordinate (C) at (-1,-1);
\coordinate (D) at (1,-2);
\coordinate (E) at (-1,-2);

\draw (E)--(C)--(A)--(B)--(D);
\draw (-1.5,-2.5)--(E)--(-.5,-2.5);
\draw (1.5,-2.5)--(D)--(.5,-2.5);
\node at (-1,-3) {$T'$};
\node at (1,-3) {$T''$};
\draw [dashed] (A)--(0,1.3);
\draw [dashed] (C)--(-2.2,-.5);
\draw [dashed] (B)--(2.2,-.5);

\draw [fill= black] (A) circle [radius = .12];
\draw [fill= black] (D) circle [radius = .12];
\draw [fill= black] (E) circle [radius = .12];
\draw [fill= white] (B) circle [radius = .12];
\draw [fill= white] (C) circle [radius = .12];

\node at (4,-1) {$\xrightarrow{\text{arborization}}$};

\begin{scope}[xshift = 8.5cm]
\coordinate (DD) at (1,-2);
\coordinate (EE) at (-1,-2);
\coordinate (BB) at (1,-1);

\draw (-1.5,-2.5)--(EE)--(-.5,-2.5);
\draw (1.5,-2.5)--(DD)--(.5,-2.5);
\draw [dashed] (EE)--(0,1.3);
\draw (DD)--(BB);
\draw [dashed] (BB)--(-2.2,-.5);
\draw [dashed] (BB)--(2.2,-.5);
\node at (-1,-3) {$T'$};
\node at (1,-3) {$T''$};

\draw [fill= black] (DD) circle [radius = .12];
\draw [fill= black] (EE) circle [radius = .12];
\draw [fill= white] (BB) circle [radius = .12];
\end{scope}
\end{tikzpicture}

\caption{An arborization step indicated schematically. The tensor diagram $T$ has two copies, $T'$ and $T''$, of the same binary tree, connecting to the boundary in the same way. These trees are joined by a path of length~$4$. The dashed lines indicate how $T'$ and $T''$ are connected to the rest of $T$. The arborization step removes the path of length $4$ and connects $T'$ and $T''$ to the rest of the diagram as indicated on the right-hand side. A similar arborization step holds with all of the colors reversed. 
\label{fig:arborization}}
\end{center}
\end{figure}

\begin{defn}[Arborization] The \emph{arborization algorithm} $T \mapsto \Arb(T)$ transforms a tensor diagram~$T$ by repeating applying \emph{arborization steps}, each of which is a a diagrammatic move given in Figure~\ref{fig:arborization}, until such steps are no longer possible. To perform such a step, locate two rooted binary trees $T'$ and $T''$ in $T$, which are isomorphic rooted binary trees connecting to the same set of boundary vertices, and whose roots are joined by a path of length $4$ (cf.~Figure \ref{fig:arborization}). The arborization step removes this path as indicated by Figure~\ref{fig:arborization}. As a degenerate case of this, we allow the case that $T'$ and $T''$ are the same boundary vertex, in which case the path of length $4$ is a $4$-cycle based at a boundary vertex. A web invariant~$[W]$ is \emph{arborizable} if $\Arb(W)$ has no interior cycles.   
\end{defn}


The process of applying several arborization steps is \emph{confluent} \cite[Theorem 10.5]{tensors}, meaning the resulting tensor diagram $\Arb(T)$ does not depend on the choice of arborization steps that are used. Each arborization step does not change the value of the invariant defined by~$T$, and consequently~$[\Arb(T)] = [T]$. This follows by applying the crossing removal skein relation to the right hand side of Figure~\ref{fig:arborization} -- one of the two resulting terms vanishes by the degeneracy skein relation.

\begin{example}\label{eg:arborizablewebs} For a \emph{planar} tensor diagram, the arborization algorithm can only start if there is a four-cycle based at a boundary vertex. Since the web~$W$ in Figure \ref{fig:SCW} has no such 4-cycles, it is equal to its own arborization (the same is true for any Pl\"ucker coordinate). This web is indecomposable, but it is \emph{not} arborizable since it has an interior cycle. On the other hand, every Pl\"ucker coordinate is an example of an indecomposable arborizable web.  

The third web in Figure \ref{fig:expressproduct} is not indecomposable. By applying an arborization step to the $4$-cycle at the boundary, one sees that its arborized form is a union of two tripods. 
Thus, Figure \ref{fig:expressproduct}~expresses the three-term Pl\"ucker relation $\Delta_{124}\Delta_{135} = \Delta_{123}\Delta_{145}+\Delta_{134}\Delta_{125}$
in the language of tensor diagrams. 
\end{example}

The following conjecture summarizes the (predicted) cluster combinatorics of $\bbc[\tGr(3,n)]$ (cf.~\cite[Sections 9 and 10]{tensors}): 
\begin{conj}[Fomin-Pylyavskyy]\label{conj:GrThreeNine} In the cluster algebra $\bbc[\tGr(3,n)]$:
\begin{enumerate}
\item The set of cluster (and frozen) variables coincides with the set of indecomposable arborizable web invariants. 
\item Two cluster variables lie in the same cluster if and only if they are compatible web invariants. 
\item If $n \geq 9$, there are infinitely many indecomposable non-arborizable web invariants. 
\end{enumerate}
\end{conj}

Conjecture \ref{conj:GrThreeNine} has been verified in the finite type examples, i.e. for $n < 9$.

\begin{rmk}
The third part of Conjecture~\ref{conj:GrThreeNine} is relevant in light of the expected link between cluster algebras and canonical bases. It has long been expected that there is a naturally defined linear basis for any cluster algebra that contains the cluster monomials, cf.~recent breakthroughs in \cite{IKLP, GHKK, KKKO}. The non-arborizable webs in Conjecture \ref{conj:GrThreeNine} are not expected to be cluster monomials. These webs should play a distinguished role in comparing the different versions of canonical bases for $\bbc[\tGr(3,n)]$. 
\end{rmk}

The following theorem is one of our main applications of the braid group action. 

\begin{thm}\label{thm:GrThreeNine} In the cluster algebra $\bbc[\tGr(3,9)]$:
\begin{enumerate}
\item Every cluster variable is an indecomposable arborizable web invariant.
\item Every cluster monomial is a web invariant (thus, cluster variables in any cluster are pairwise compatible). 
\item There are infinitely many indecomposable non-arborizable web invariants. 
\end{enumerate}
\end{thm}

We emphasize that we do not prove the reverse implications of Fomin-Pylyavskyy's conjectures, namely that every indecomposable arborizable non-elliptic web invariant \emph{is} a cluster variable, and also, that two cluster variables are in a cluster whenever they are compatible web invariants. 

Our proof of part (3) in Theorem~\ref{thm:GrThreeNine} establishes that the ``single cycle web'' in Figure~\ref{fig:SCW} has an infinite orbit with respect to the dot action of $B_3$, and every web in this orbit is indecomposable and non-arborizable. We conjecture that these are the only indecomposable non-arborizable webs.   

Our second main theorem is a presentation for the cluster modular group. The four generators are the cyclic shift $\rho$, the cyclic shift $P$ on $\FG(3,6)$ brought over to $\Gr(3,9)$, the twist map $\tau$, and the reflection $\theta$. 
\begin{thm}\label{thm:CMG3presentation} The cluster modular group $\CMG(\Gr(3,9))$ has the presentation
\begin{equation}\label{eq:CMG3presentation}
\begin{split}
\CMG = \la \rho,P,\tau,\theta \colon \, & \rho^3 = P^2 = \iota^{-2}, \rho^9 = 1 ,  \:  \tau \rho = \rho \tau , \: \tau P = P \tau  \\
 \:  \: \: \:  & \theta^2=1, \theta\rho \theta = \rho^{-1}, \: \theta P \theta= P^{-1}, \: \theta \tau \theta = \tau^{-1} \ra.
\end{split}
\end{equation}
We obtain an isomorphism $\CMG^+ / Z(\CMG^+) \cong \PSL_2(\bbz).$ 
\end{thm}

That is, besides the order of the cyclic shift, the statement that $\tau$ is central, and the description of how conjugation by $\theta$ permutes the generators, the only interesting relation is that $\rho^3 = P^2 = \tau^{-2}$. 
The proofs of Theorems~\ref{thm:GrThreeNine} and~\ref{thm:CMG3presentation} are in Section~\ref{secn:proofs}.

\subsection{Webs for $\FG(3,r)$}
To prove Theorems~\ref{thm:GrThreeNine} and \ref{thm:CMG3presentation}, we will move back and forth between the spaces $\Gr(3,9)$ and $\Conf(3,6)$. Following \cite[Section 12]{tensorsII}, it is convenient to have a notion of tensor diagrams, webs, and arborization, for the Fock-Goncharov cluster algebras. Our proof of Theorem~\ref{thm:GrThreeNine} establishes an analogous statement for $\Conf(3,6)$. 

\begin{figure}[ht]
\begin{center}
\begin{tabular}{ll}
boundary 3-cycles &  
\begin{tikzpicture}[scale = .7]
\draw [gray] (0,0) arc [radius=2, start angle=-90, end angle= -60];
\draw [gray] (0,0) arc [radius=2, start angle=-90, end angle= -120];
\draw (.15,.02)--(.5,1)--(-.5,1)--(0,0);
\draw [dashed] (.5,1)--(1,1.5);
\draw [dashed] (-.5,1)--(-1,1.5);
\draw [fill= black] (0,0) circle [radius = .065];
\draw [fill= white] (.15,.02) circle [radius = .065];
\draw [fill= black] (.5,1) circle [radius = .065];
\draw [fill= white] (-.5,1) circle [radius = .065];
\begin{scope}[xshift = 5.25cm]
\node at (-2,0) {$=$};
\draw [gray] (0,0) arc [radius=2, start angle=-90, end angle= -60];
\draw [gray] (0,0) arc [radius=2, start angle=-90, end angle= -120];
\draw [dashed] (.15,.02)--(1,1.5);
\draw [dashed] (0,0)--(-1,1.5);
\draw [fill= black] (0,0) circle [radius = .065];
\draw [fill= white] (.15,.02) circle [radius = .065];
\end{scope}[xshift = 3cm]
\end{tikzpicture}
\\
degeneracy& 
\begin{tikzpicture}[scale = .7]
\draw [gray] (0,0) arc [radius=2, start angle=-90, end angle= -60];
\draw [gray] (0,0) arc [radius=2, start angle=-90, end angle= -120];
\draw (.15,.02) arc [radius = .3, start angle =-80,end angle = 245];
\draw [fill= black] (0,0) circle [radius = .065];
\draw [fill= white] (.15,.02) circle [radius = .065];
\node at (4.2,.5) {$= \hspace{1cm} 0$};
\node at (3.5,1.1) {\hfill };
\end{tikzpicture}
\end{tabular}
\end{center}
\caption{Additional skein relations for bivariant tensor diagrams. Both take place at a black and white vertex pair. \label{fig:SL3skeinrelationsII}}
\end{figure}
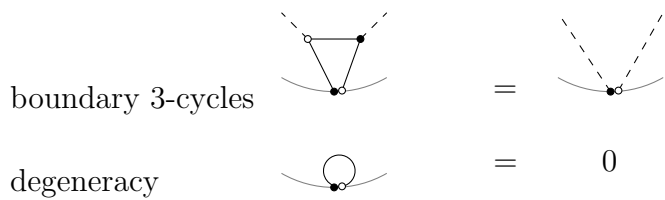

\medskip
A \emph{bivariant tensor diagram} is a tensor diagram drawn in the disk with $r$ ``colorless'' boundary vertices, subject to the same constraints as before (interior trivalent and bipartite). Each half-edge connecting to a boundary vertex is decorated either black or white, and the colorless boundary vertex is considered black or white accordingly in the definition of bipartiteness when it is being used in such a half-edge. In pictures, we draw the colorless boundary vertices as a black and white vertex pair glued to each other, but the cyclic orientation of the black and white part of a colorless boundary vertex is not relevant.  The $\Conf(3,r)$ \emph{skein algebra} is the space spanned by bivariant tensor diagrams, subject to the skein relations in Figure~\ref{fig:SL3skeinrelations} and the additional skein relations in Figure~\ref{fig:SL3skeinrelationsII}. The \emph{bi-invariant} $[T]$ associated to $T$ is the image of $T$ in the skein algebra. 

A \emph{non-elliptic biweb} is defined as before, with the additional requirement that there are no loops or $3$-cycles based at boundary vertices. A \emph{biweb invariant} is the image of a non-elliptic biweb in the skein algebra. Kuperberg's theorem still holds: biweb invariants are a basis for the skein algebra. The arborization algorithm for bivariant tensor diagrams allows for the steps in Figure~\ref{fig:arborization}, as well as the boundary 3-cycle removal step. 
 
Through a more intrinsic definition, one can see that any bivariant tensor diagram $T$ defines a function $[T]$ on a configuration of affine flags $(F_{1,\bullet},\dots,F_{r,\bullet})$ (cf.~Example~\ref{eg:FGwebs}). The $i$th affine flag $F_{i,\bullet} = (F_{i,(1)},F_{i,(2)})$ corresponds to the $i$th boundary vertex. The black part of the bi-colored vertex corresponds to the vector $F_{i,(1)}$, and the white part of the vertex corresponds to $F_{i,(2)} \in \bigwedge^2(V)$. We have this more intrinsic definition in mind, but it will not play a very important role in our proofs.

\begin{example}\label{eg:FGwebs}
The Fock-Goncharov coordinates \eqref{eq:FGCoord} for $\bbc(\Conf(3,r))$ correspond to either 1) bipods, i.e. an arc connecting a black part of one boundary vertex to the white part of another boundary vertex, or 2) tripods connecting the black parts of three boundary vertices. The duality map $\ast$ acts on bivariant tensor diagrams by globally swapping the colors black and white. The frozen variables are the $2r$ bipods joining consecutive boundary vertices. 
\end{example}

\subsection{Webs for $\SL_4$}
By a similar approach to the one we carry out in detail for $\Gr(3,9)$, the braid group action allows us to give a description of all cluster variables for $\Gr(4,8)$. We give a very quick introduction to tensor diagrams in this case. 

Tensor diagrams for $\Gr(4,n)$ are drawn in a disk with $n$ black boundary vertices. Unlike tensor diagrams for $\Gr(3,n)$, each edge is labeled with multiplicity (either $1$ or $2$), so that the sum of multiplicities around every interior vertex is $4$. All boundary edges have multiplicity 1. The tensor diagram should be bipartite, with every interior vertex either bivalent or trivalent. We draw multiplicity 2 edges as ``double bonds,'' and multiplicity 1 edges as ordinary. 

As before, there is a set of skein relations amongst such tensor diagrams for $\Gr(4,n)$ (the crossing removal \cite[Corollary 6.2.3]{SkewHowe} together with the relations amongst planar diagrams \cite[Chapter 4]{KimThesis}), and there is a way of interpreting each such tensor diagram $T$ as an invariant $[T] \in \bbc[\tGr(4,8)]$. We note that, at least in the most naive way, the invariant $[T] \in \bbc[\tGr(4,n)]$ represented by a tensor diagram $T$ is only well-defined up to a sign. For brevity's sake, we will forgo a careful discussion of these signs, referring the reader to \cite{SkewHowe} or \cite{FraserLamLe} for some possible choices of conventions. We illustrate the smallest non-Pl\"ucker $\SL_4$ webs (both of which are cluster variables), called ``octapods'': 
\begin{equation}\label{eq:octapod}
\begin{tikzpicture}
\draw (45:1cm)--(90:.5cm)--(135:1cm);
\draw (180:1cm)--(225:.5cm)--(270:1cm);
\draw (90:1cm)--(90:.5cm);
\draw (225:1cm)--(225:.5cm);
\draw (315:1cm)--(337.5:.5cm)--(0:1cm);
\draw (0,0)--(90:.5cm);
\draw (0,0)--(225:.5cm);
\draw (.05cm,.05cm)--(.5119cm,-.1413cm);
\draw (-.05cm,-.05cm)--(.4119cm,-.2413cm);
\filldraw[white] (90:.5cm) circle (0.1cm);
\draw (90:.5cm) circle (0.1cm);
\filldraw[white] (225:.5cm) circle (0.1cm);
\draw (225:.5cm) circle (0.1cm);
\filldraw[black] (0,0) circle (0.1cm);
\draw (0,0) circle (0.1cm);
\filldraw[white] (337.5:.5cm) circle (0.1cm);
\draw (337.5:.5cm) circle (0.1cm);
\begin{scope}[xshift = 3cm]
\draw (45:1cm)--(90:.5cm)--(150:1cm);
\draw (150:1cm)--(225:.5cm)--(270:1cm);
\draw (90:1cm)--(90:.5cm);
\draw (225:1cm)--(225:.5cm);
\draw (315:1cm)--(337.5:.5cm)--(0:1cm);
\draw (0,0)--(90:.5cm);
\draw (0,0)--(225:.5cm);
\draw (.05cm,.05cm)--(.5119cm,-.1413cm);
\draw (-.05cm,-.05cm)--(.4119cm,-.2413cm);
\filldraw[white] (90:.5cm) circle (0.1cm);
\draw (90:.5cm) circle (0.1cm);
\filldraw[white] (225:.5cm) circle (0.1cm);
\draw (225:.5cm) circle (0.1cm);
\filldraw[black] (0,0) circle (0.1cm);
\draw (0,0) circle (0.1cm);
\filldraw[white] (337.5:.5cm) circle (0.1cm);
\draw (337.5:.5cm) circle (0.1cm);
\end{scope}
\end{tikzpicture}
\end{equation}

Using the $\SL_4$ skein relations, one can write any element of the skein algebra as a linear combination of planar diagrams without $0$-cycles or $2$-cycles (once we have removed $2$-cycles, our double bond drawings are unambiguous). We call these latter planar diagrams \emph{web invariants}. The key difference 
between $\Gr(3,n)$ and $\Gr(4,n)$ is that web invariants are merely a spanning set. There is no convenient analogue of the notion of non-elliptic, thus no known natural way of identifying a smaller subset of basis web invariants. We have the following relationship between cluster variables and $\SL_4$ webs:

\begin{thm}\label{thm:GrFourEight} In the cluster algebra $\bbc[\tGr(4,8)]$:
\begin{itemize} 
\item Every cluster variable $x$ is proportional to a web invariant $[W]$ that can also be expressed as $[W] = [T]$, where $T$ is a tensor diagram $T$ with no interior cycles. 
\item Up to the action of the cluster modular group, every cluster variable is either a Pl\"ucker coordinate or an octapod (cf.~\eqref{eq:octapod}). 
\end{itemize}
\end{thm}

We call $T$ as above the \emph{tree form} of the cluster variable. Note that we believe that in fact, every cluster variable \emph{is} (not merely is proportional to) a web invariant with a tree form, but we do not carefully prove this. The property of having a tree form is an analogue of the condition of arborizability for $\SL_3$ webs. The following diagrammatic move is the analogue of Figure~\ref{fig:arborization} for $\Gr(4,8)$ webs:
\begin{equation}\label{eq:SL4arb}
\begin{tikzpicture}
\node at (2,.5) {\large $=$};
\draw [gray] (0,0) arc [radius=2, start angle=-90, end angle= -60];
\draw [gray] (0,0) arc [radius=2, start angle=-90, end angle= -120];
\draw (-.5,.5)--(0,0)--(.5,.5)--(0,1);
\draw (-.5,.55)--(0,1.05);
\draw (-.5,.45)--(0,.95);
\draw [dashed] (-.5,.5)--(-1.1,.9);
\draw [dashed] (0,1)--(-.3,1.5);
\draw [dashed] (.5,.55)--(1.1,.95);
\draw [dashed] (.5,.45)--(1.1,.85);
\draw [fill= black] (0,0) circle [radius = .06];
\draw [fill= black] (0,1) circle [radius = .06];
\draw [fill= white] (-.5,.5) circle [radius = .06];
\draw [fill= white] (.5,.5) circle [radius = .06];
\begin{scope}[xshift = 4cm]
\draw [gray] (0,0) arc [radius=2, start angle=-90, end angle= -60];
\draw [gray] (0,0) arc [radius=2, start angle=-90, end angle= -120];
\draw (0,0)--(.1,.5);
\draw [dashed] (0,0)--(-.3,1.5);
\draw [dashed] (.1,.5)--(-1.1,.9);
\draw [dashed] (.1,.55)--(1.1,.95);
\draw [dashed] (.1,.45)--(1.1,.85);
\draw [fill= white] (.1,.5) circle [radius = .06];
\draw [fill= black] (0,0) circle [radius = .06];
\end{scope}
\end{tikzpicture}.
\end{equation}
More precisely, the move~\eqref{eq:SL4arb} is the analogue of removing a 4-cycle at the boundary of an $\SL_3$ web, and the general arborization move is one in which the black boundary vertex in \eqref{eq:SL4arb} is replaced by two isomorphic copies of the same tree connecting to the boundary. It is a consequence of our proof of Theorem~\ref{thm:GrFourEight} that every cluster variable $x \in \bbc[\tGr(4,8)]$ has a planar form $x = [W]$ that can be converted to tree form by repeated applications of this arborization move. However, we are not bold enough to conjecture the analogue of (1) from Conjecture~\ref{conj:GrThreeNine} holds, i.e. a web invariant in $\bbc[\tGr(4,n)]$ is a cluster variable if and only if it can be transformed into tree form by repeating~\eqref{eq:SL4arb}. 

We have the following presentation for the cluster modular group. 
\begin{thm}\label{thm:CMG4presentation} The cluster modular group $\CMG(\Gr(4,8))$ has the presentation
\begin{equation}\label{eq:CMG4presentation}
\begin{split}
\CMG = \la \sigma_1,\sigma_2,\sigma_3,\tau,\theta: & \: \sigma_1 \sigma_2 \sigma_1 = \sigma_2 \sigma_1 \sigma_2, \: \sigma_2 \sigma_3 \sigma_2 = \sigma_3 \sigma_2 \sigma_3, \: \sigma_1 \sigma_3 = \sigma_3 \sigma_1, \\
& \: \sigma_1 \sigma_2 \sigma_3^2 \sigma_2 \sigma_1 = (\sigma_3 \sigma_2 \sigma_1)^8 =1,  \: \tau \sigma_i = \sigma_i \tau \text{ for $i = 1,\dots,3$ },  \\
&  \: \theta^2 = 1, \: \theta \tau \theta = \tau^{-1} , \:  \theta \sigma_i \theta = \sigma_{4-i}^{-1} \text{ for $i = 1,\dots,3$ }. \\
\end{split}
\end{equation}
Likewise, we obtain an isomorphism $\CMG^+/ Z(\CMG^+) \cong \MCG(S^2,4)$.
\end{thm}

Besides the order of the cyclic shift, the centrality of $\tau$, and the description of conjugation by $\theta$, the only interesting relations are the braid relations and the relation $\sigma_1 \sigma_2 \sigma_3^2 \sigma_2 \sigma_1=1$. 

\begin{rmk}[Groups in finite mutation type]\label{rmk:fmtype} The (skew-symmetric) quivers of finite mutation type and of rank $\geq 3$ were classified in \cite{FiniteType}. Besides the quiver mutation classes of surface type \cite{CATSI}, for whom the cluster modular group is closely related to the mapping class group of the surface \cite[Proposition 8.5]{BridgelandSmith}, there are 11 more quiver mutation classes of finite mutation type. Six of these extra 11 types are finite or affine type, with the groups computed in \cite{ASS}. Theorems~\ref{thm:CMG3presentation} and \ref{thm:CMG4presentation} address two of the remaining five cases -- namely the extended affine types $E_8^{(1,1)}$ and $E_7^{(1,1)}$ respectively. This leaves three remaining finite mutation types -- $E_6^{(1,1)}$ and $X_6,X_7$ \cite{DerksenOwen}. It remains an open problem to compute the cluster modular group, and to give a combinatorial description of the cluster variables and clusters, in these three cases. Significant progress has been made for $E_6^{(1,1)}$ \cite{TubularII} and for $X_6,X_7$ \cite{Ishibashi}. It would also be interesting to answer these questions for the skew-symmetrizable examples of finite mutation type (cf.~\cite{FiniteTypeII}).    
\end{rmk}

\section{Proofs for finite mutation type}\label{secn:proofs}
Our goal in this section is to describe all cluster variables in $\Gr(3,9)$ and $\Gr(4,8)$, and also to give a presentation for both cluster modular groups. We begin with some generalities and a summary of the computer verifications underlying our proofs.  

Let $f \colon \mca \to \mca $ be a quasi-automorphism of a cluster algebra. If $x \in \mca$ is a cluster variable, then $f(x)$ is proportional to a cluster variable $\ov{x} \in \mca$. We henceforth use the notation 
\begin{equation}\label{eq:dotaction}
f \cdot x = \ov{x}
\end{equation}
to summarize the situation $f(x) \propto \ov{x}$ and $\ov{x}$ is a cluster variable, and refer to this as the \emph{dot action} of $f$ on cluster variables. For $\SL_3$ webs or biwebs, we use the same dot action notation $f \cdot [W] = [W']$ to summarize the situation that $W,W'$ are indecomposable web invariants and $M$ is a monomial in the frozens.

The following lemma appears in \cite[Appendix 1]{FGXVarieties}. We rewrite a proof in our own language because the construction underlies the proof of Theorem~\ref{thm:CMG3presentation}.    
\begin{lem}\label{lem:finitegeneration}
If $Q_0$ is a quiver (without frozen vertices) of finite mutation type, then the cluster modular group for the cluster algebra $\mca(Q_0)$ is finitely generated.
\end{lem}

\begin{proof}
Let $\Gamma$ be the graph whose vertices are the (isomorphism classes of) quivers mutation equivalent to $Q_0$: two vertices $Q_1,Q_2 \in \Gamma$ are connected by an edge if 
\begin{equation}\label{eq:somek}
\mu_k(Q_1) \cong Q_2 \text{ for \emph{some} direction $k$.}
\end{equation}
Note that whether or not \eqref{eq:somek} holds for some $k$ is symmetric in $Q_1,Q_2$, but the number of such $k$'s is not always symmetric. 

Now consider the fundamental group of $\Gamma$ based at the quiver $Q_0$. Since $\Gamma$ is a finite graph, the fundamental group is finitely generated by generators $g_1,\dots,g_\ell$. For simplicity, we take a list of generators that is closed under taking inverses. For each such generating cycle $g_i$, lift it in all possible ways to a sequence of mutations in $\mca(Q_0)$. That is, for each edge $Q_1 \xrightarrow{e} Q_2$ in $g_i$, replace $e$ by every mutation $\mu_k$ as in \eqref{eq:somek}. We obtain a finite list of mutation sequences, $h_1,\dots,h_L$, each of which returns the quiver $Q_0$ to itself. Some of these $h_i$ could represent trivial elements of the cluster modular group (for an example: $\pi_1(\Gamma)$ will contain $4$-cycles arising from commuting mutations $\mu_i \circ \mu_j = \mu_j \circ \mu_i$). Furthermore, different $h_i$'s could represent the same element of the cluster modular group. 

Now pick a choice of initial seed $\Sigma_0$ for $\mca(Q_0)$, whose underlying quiver is $Q_0$. Let $g$ be an element of the cluster modular group for $\mca(Q_0)$. Then $g$ sends $\Sigma_0$ to some new seed $\Sigma_g$ in $\mca(Q_0)$. There exists a sequence $\mu$ of mutations from $\Sigma_0$ to $\Sigma_g$, and this sequence descends to a cycle in $\Gamma$. This cycle in $\Gamma$ can be expressed as a product of the generators $g_1,\dots,g_\ell$, hence $\mu$ can be written as a composition of the $h_1,\dots,h_\ell$. Thus, the elements of the cluster modular group coming from $h_1,\dots,h_\ell$ generate $\CMG$.
\end{proof}

The following lemma isolates the computer verifications underlying our results. The relevant software is included as an ancillary file with the arXiv version of this paper \cite{ArxivFraser}. For a cluster algebra of finite mutation type, a \emph{fundamental domain} is a specific choice of seeds exhausting the finitely many mutable subquivers. 
\begin{lem}\label{lem:computercheck} The following facts were checked on a computer: 
\begin{itemize}
\item In both cases $\Gr(3,9)$ and $\Gr(4,8)$, the cluster modular group is generated by the Artin generators $\sigma_i$ together with the reflection symmetry $\theta$ and the twist $\tau$. 
\item On a fundamental domain for $\Gr(3,9)$, every cluster monomial is a web invariant, and every cluster variable is indecomposable and arborizable.  
\item On a fundamental domain for $\Gr(4,8)$, every cluster variable is an $\SL_4$ web that can also be written in tree form.  
\end{itemize}
\end{lem}

\begin{proof}
To check the statement about generators for the cluster modular group, we implement the argument in Lemma~\ref{lem:finitegeneration}. For $\Gr(3,9)$ there are 5739 isomorphism classes of quivers, and the quiver exchange graph $\Gamma$ from Lemma~\ref{lem:finitegeneration} has 22007 independent cycles. For $\Gr(4,8)$ there are 506 isomorphism classes of quivers, and the quiver exchange graph $\Gamma$ from Lemma~\ref{lem:finitegeneration} has 1506 independent cycles. In either case, for each cycle $C$ in $\Gamma$, we perform the sequence of mutations described by $C$ starting from an initial seed $\Sigma_1$, and check that this sequence of mutations can be realized as a composition of the Artin generators $\sigma_i$ and the twist $\tau$. Then we do the same for each two-cycle. 

Let us move on to the statements about cluster variables and clusters on a fundamental domain. In $\Gr(3,9)$, a breadth-first search produces a choice of fundamental domain of seeds $\Sigma_1,\dots,\Sigma_{5739}$, one for each quiver isomorphism class. We computed the cluster variables in these seeds as webs using the program \cite{ArxivFraser}. To do this, we encoded webs of small degree (e.g. Pl\"uckers, hexapods, and so on) as explicit polynomials in the coordinates on $V^9$, so that each exchange relation can be verified directly. We wrote a compatibility tester for webs of small degree. To perform a mutation, we eliminate all webs of the correct degree that are not compatible with the current cluster, and check that the corresponding exchange relation holds. 
 
All of the $\SL_3$ webs that show up in $\Sigma_1,\dots,\Sigma_{5739}$ are arborizable. By construction, any two webs in a given cluster $\Sigma_i$ are pairwise compatible. For arborizable webs of small degree, it is not hard to convince oneself that if all webs in $\Sigma_i$ are pairwise compatible, then any  \emph{monomial} in $\Sigma_i$ is again a non-elliptic web.  

The computation for $\Gr(4,8)$ is similar, but weaker because we did not write a compatibility tester for $\SL_4$ webs (it is tedious to make sure one has identified all web invariants in a given degree because they are a spanning set rather than a basis; as a consequence, it is tedious to show that a given element in the skein algebra is not a web invariant). We checked that every $\SL_4$ web in the fundamental domain has a tree form, but we did not check that webs in a shared cluster are always pairwise compatible. We note that our chosen fundamental domain contains a handful of cluster variables that are neither Pl\"ucker coordinates nor octapods, but we checked by hand that each of these could be moved to a Pl\"ucker or octapod via the braid group. 
\end{proof}

\subsection{Proofs for $\Gr(3,9)$}
We now set out to prove Theorem~\ref{thm:GrThreeNine} and Theorem~\ref{thm:CMG3presentation}. By a \emph{fork} in a (bivariant) tensor diagram $T$ we will mean an interior vertex $v$ of $T$ that is adjacent to consecutive boundary vertices $i$ and $i+1$ (considered modulo $n$), both of which necessarily have the same color. We identify this fork with the subgraph of $T$ connecting the vertex $v$ to the two boundary vertices. If $i,i+1$ are black boundary vertices in a tensor diagram $T$ for $\Gr(3,9)$, then a fork between these vertices encodes the exterior product $v_iv_{i+1} \in \bigwedge^2(V)$ in the invariant $[T]$. Likewise, a fork between the black parts of boundary vertices $i,i+1$ in a tensor diagram for $\FG(3,6)$ encodes the exterior product $F_{i,(1)}F_{i+1,(1)} \in \bigwedge^2(V)$, whereas a fork between the white parts encodes the vector $F_{i,(2)}\cap F_{i+1,(2)} \in V$. These forks play an important role in our statement of Lemma~\ref{lem:HowToDrawII}.

By Lemma~\ref{lem:computercheck} the group $\CMG^+(\FG(3,6))$ is generated by the duality map $\ast^*$, the cyclic shift of affine flags $P^*$, and the map $\Phi^* \circ \rho^* \circ \Psi^*$, which is the cyclic shift on $\Gr(3,9)$ translated over to $\FG(3,6)$. Recall that $\ast$ acts by swapping the colors black and white, and $P$ acts by rotation. Let us explain how to evaluate the third map $\Phi^* \circ \rho^* \circ \Psi^*$ on a tensor diagram. 

\begin{figure}
\begin{center}
\begin{tikzpicture}[scale = 1.5]
\begin{scope}[xshift = -5cm]
\draw [black] (1,0) arc [radius = 1, start angle = 0, end angle = 360];
\foreach \s in {1,...,6}
{
  \node at ({360/6 * (-\s+4)-60}:1.3 cm) {$\s$};
  \draw [fill= black] ({360/6 * (\s - 1)}:1 cm) circle [radius = .04];
}

\draw (30:.8cm)--(30:.5cm);
\draw (90:.8cm)--(90:.5cm);
\draw (210:.8cm)--(210:.5cm);
\draw (270:.8cm)--(270:.5cm);
\draw (120:1cm)--(120:.5cm);
\draw (0:1cm)--(0:.5cm);
\draw (180:1cm)--(180:.5cm);
\draw (300:1cm)--(300:.5cm);

\node at (135:.7cm) {\small $F_{11}$};
\node at (-10:.7cm) {\small $F_{21}$};
\node at (170:.7cm) {\small $F_{41}$};
\node at (315:.7cm) {\small $F_{31}$};

\node at (45:.65cm) {\small $F_{22}$};
\node at (105:.7cm) {\small $F_{12}$};
\node at (225:.7cm) {\small $F_{42}$};
\node at (285:.7cm) {\small $F_{32}$};

\draw (120:1cm)--(90:.8cm);
\draw (60:1cm)--(90:.8cm);
\draw [fill= white] (90:.8cm) circle [radius = .04];
\draw (60:1cm)--(30:.8cm);
\draw (0:1cm)--(30:.8cm);
\draw [fill= white] (30:.8cm) circle [radius = .04];
\draw (240:1cm)--(210:.8cm);
\draw (180:1cm)--(210:.8cm);
\draw [fill= white] (210:.8cm) circle [radius = .04];
\draw (240:1cm)--(270:.8cm);
\draw (300:1cm)--(270:.8cm);
\draw [fill= white] (270:.8cm) circle [radius = .04];
\node at (0.1,-1.4) {$\Psi^* \colon \bbc[\FG(3,4)] \to  \bbc[\tGr(3,6)]$};
\end{scope}

\begin{scope}[xshift = 0cm]
\draw [black] (1,0) arc [radius = 1, start angle = 0, end angle = 360];

\foreach \s in {1,...,4}
{
  \node at ({360/4 * (-\s+4)-135}:1.3 cm) {$\s$};
  \draw [fill= black] ({360/4 * (\s - 1)-45}:1 cm) circle [radius = .04];
}

\draw (130:1cm)--(90:.6cm);
\draw (50:1cm)--(90:.6cm);
\draw (230:1cm)--(270:.6cm);
\draw (310:1cm)--(270:.6cm);

\draw (45:1cm)--(35:.6cm);
\draw (135:1cm)--(145:.6cm);
\draw (225:1cm)--(215:.6cm);
\draw (315:1cm)--(325:.6cm);
\draw (90:.6cm)--(90:.4cm) ;
\draw (270:.6cm)--(270:.4cm) ;

\node at (155:.7cm) {\small $v_1$};
\node at (25:.7cm) {\small $v_3$};
\node at (110:.3cm) {\small $v_2$};
\node at (205:.7cm) {\small $v_6$};
\node at (335:.7cm) {\small $v_4$};
\node at (290:.4cm) {\small $v_5$};

\draw [fill= white] (130:1 cm) circle [radius = .04];
\draw [fill= white] (50:1 cm) circle [radius = .04];
\draw [fill= white] (230:1 cm) circle [radius = .04];
\draw [fill= white] (310:1 cm) circle [radius = .04];
\draw [fill= black] (90:.6cm) circle [radius = .04];
\draw [fill= black] (270:.6cm) circle [radius = .04];
\node at (0.1,-1.4) {$\Phi^* \colon \bbc[\tGr(3,6)] \to  \bbc[\FG(3,4)]$};
\end{scope}

\end{tikzpicture}
\caption{To evaluate $\Psi^*$ on a tensor diagram $[T] \in \bbc[\FG(3,4)]$, reattach edges incident to $F_{i,(a)}$ in $T$ by the fragment labeled by $F_{ia}$ in the first figure, for $i=1,\dots,4$ and $a = 1,2$. The result is a tensor diagram in $\bbc[\Gr(3,6)]$. To evaluate $\Phi^*$ on a tensor diagram $[T] \in \bbc[\Gr(3,6)]$, reattach edges incident to boundary vertex $v_i$ in $T$ to according to the strand labeled by $v_i$ in the second figure, for $i=1,\dots,6$. 
\label{HowToDrawFig}}
\end{center}
\end{figure}
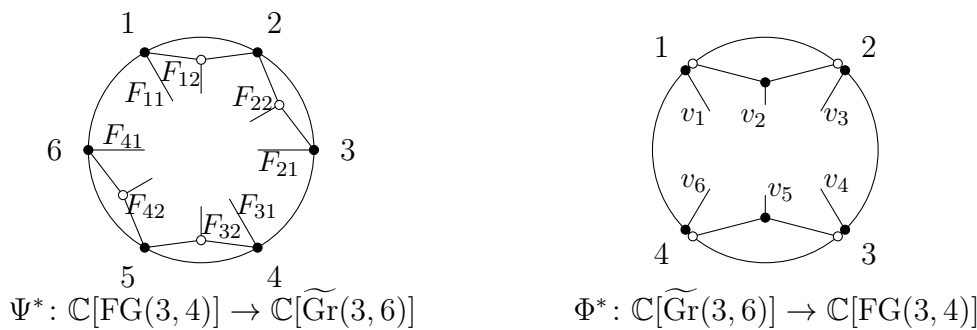

Since $\Psi$ is induced by a map from vectors to affine flags \eqref{eq:Psidefn}, evaluating $\Psi^*$ corresponds to an operation on the boundary edges of the tensor diagram: the interior edges of the tensor diagram are not affected, but boundary edges are reattached via fragments which encode the substitution \eqref{eq:Psidefn}. The first half of Figure~\ref{HowToDrawFig} illustrates this substitution for $\Psi^* \colon \bbc[\FG(3,4)] \to \bbc[\tGr(3,6)]$. For example, since $F_{1,\bullet}=v_1,v_1v_2$ in \eqref{eq:Psidefn}, every edge $e$ in $T$ that connects to the black boundary vertex $F_{1,(1)}$ gives rise to an edge connecting to the first boundary vertex in $\Psi^*([T])$, while an edge connecting to the white boundary vertex $F_{1,(2)}$ gives rise to a fork between the first two boundary vertices in $\Psi^*([T])$. The diagram encoding $\Phi^*$ has a similar flavor, pictured in the second diagram of Figure~\ref{HowToDrawFig}. For an example of how to carry out these reattaching operations, see Example~\ref{eg:FGForks} below.

Our first main assertion is that the maps $\Psi^*$ and $\Phi^*$ preserve the web combinatorics 
from Section~\ref{secn:Webs}, so we may translate freely between $\bbc[\tGr(3,3r)]$ and $\bbc[\FG(3,2r)]$. 

\begin{lem}\label{lem:HowToDraw}
Let $[W]$ be an indecomposable biweb invariant for $\Conf(3,2r)$, not equal to a frozen variable. Then $\Psi^*([W])$ factors as  $M [W']$ where $M$ is a product of frozens in $\bbc[\tGr(3,3r)]$, and $W'$ is an indecomposable non-elliptic web invariant for $\Gr(3,3r)$ not equal to a frozen variable. Furthermore, $W'$ is arborizable if and only if $W$ is. The analogous statements hold with the roles of $\Conf(3,2r)$ and $\tGr(3,3r)$ reversed. 
\end{lem}

\begin{proof}
Given a web invariant $[W] \in \bbc[\tGr(3,3r)]$, we evaluate $\Phi^*([W])$ by plugging in to the second diagram in Figure~\ref{HowToDrawFig}. The resulting bivariant tensor diagram will not necessarily be planar. However, the only crossings that are created come from plugging in to the fork labeled by $v_2$ several times (or one of the rotations $v_{2+3i}$ of this by 3 units). For each such crossing, applying the crossing removal skein relation produces two terms, one of which vanishes by the degeneracy relation. This planarization process adds a boundary 4-cycle to the diagram, but does not affect the interior faces (thus preserves the condition that all such interior faces have at least six sides). Resolving all of these crossings, we obtain a connected biweb whose invariant equals $\Phi^*([W])$. This diagram is not necessarily non-elliptic due to the potential presence of boundary 3-cycles. However, by inspection of Figure~\ref{HowToDrawFig}, the only possible boundary $3$-cycles created arise from plugging in forks $v_1v_2$, or plugging in forks $v_2v_3$, or rotations of $v_{i} \mapsto v_{i+3}$ of these by 3 units.  In both instances, after applying the boundary 3-cycle relation, the diagram factors into two pieces, one of which is a frozen variable and the other of which remains a connected planar diagram whose interior faces have at least six sides. After removing all such 3-cycles, the leftover diagram (with the frozen variables ignored) is a connected non-elliptic web $W'$. The removal of boundary 3-cycles did not affect the arborizability of the web (it did not affect the interior faces). To summarize, if $W$ is non-elliptic, then after applying crossing and boundary 3-cycle removals, $\Phi^*([W]) = M[W']$ where $M$ is a monomial in the frozens and 
$W'$ is a connected non-elliptic biweb invariant. 

Next we consider an indecomposable biweb $B \in \bbc(\FG(3,2r))$. We evaluate $\Psi^*([B])$ by plugging in to the first diagram in Figure~\ref{HowToDrawFig}. As above, the resulting tensor diagram might not be planar (due to crossings between the strands labeled $F_{11}$ and $F_{12}$, or by using the strand labeled by $F_{12}$ several times, or rotations of these crossings). However, it can be planarized by applying the crossing removal relation and then the degeneracy relation. Again, this does not affect the interior faces. The resulting diagram is already non-elliptic (however, it might have some boundary 4-cycles, so it might possibly factor). To summarize, if $[B]$ is an indecomposable biweb invariant, then after resolving crossings, $\Psi^*([B]) = [B']$  for a web invariant $B'$. If $B$ is not a frozen variable, then $\Psi^*([B'])$ is not a product of frozen variables (as follows by applying $\Phi^*$).  

Now we can argue that if the web invariant $W$ above is indecomposable, then so is $W'$. Since $\Psi^* \circ \Phi^* (W) = \Psi^*(M) \Psi^*(W')$, 
if $W'$ factored into a product of one or more non-frozen web invariants, then by the arguments in the second paragraph, it would follow that $\Psi^* \circ \Phi^* (W)$ does too. But from the proof of Theorem~\ref{thm:PsiPhyXY}, $\Psi^*(M) \Psi^*(W') = M' W$ for some product of frozens $M'$. We conclude that $W$ is not indecomposable (so our assumption that $W'$ factored was false). On the other hand, if $W$ factors a product of non-frozen web invariants, then $W'$ does clearly. This establishes Lemma~\ref{lem:HowToDraw} for $[W] \in \bbc[\tGr(3,3r)]$; the other direction follows by applying $\Phi^*$ and using $\Psi^* \circ \Phi^* \propto \text{Id}$.
\end{proof}

The first two parts of Theorem~\ref{thm:GrThreeNine} follow immediately:
\begin{proof}[Proof of Theorem~\ref{thm:GrThreeNine} (1) and (2)] 
By Lemma~\ref{lem:computercheck}, there are only finitely many clusters in $\bbc[\tGr(3,9)]$ up to the dot action of $B_3$. The $B_3$ action is generated by the cyclic shift $\rho$ and the composition $(\Phi \circ P \circ \Psi)^*$, and by Lemma~\ref{lem:HowToDraw} both of these generators preserves the adjectives indecomposable, arborizable, and non-elliptic, when acting by the dot action. Since every cluster variable on the fundamental domain is an indecomposable arborizable web invariant, it follows that all cluster variables are. Since every cluster monomial on the fundamental domain is a web invariant, it follows that all cluster monomials are. 
\end{proof}

\begin{figure}
\begin{center}
\begin{tikzpicture}[scale = 2.5]
\draw [black] (1,0) arc [radius = 1, start angle = 0, end angle = 360];

\foreach \s in {1,...,6}
{
  \node at ({360/6 * (-\s+6)-180}:1.3 cm) {$\s$};
  \draw [fill= black] ({360/6 * (\s - 1)}:1 cm) circle [radius = .04];
}

\draw (0:1cm)--(30:.6cm);
\draw (60:1cm)--(30:.6cm);
\draw (120:1cm)--(150:.6cm);
\draw (180:1cm)--(150:.6cm);
\draw (240:1cm)--(270:.6cm);
\draw (300:1cm)--(270:.6cm);

\draw (65:1cm)--(90:.6cm);
\draw (115:1cm)--(90:.6cm);
\draw (185:1cm)--(210:.6cm);
\draw (235:1cm)--(210:.6cm);
\draw (305:1cm)--(330:.6cm);
\draw (355:1cm)--(330:.6cm);

\draw (30:.6cm)--(30:.4cm);
\node at (45:.5cm) {\tiny $F_{22}$};
\draw (150:.6cm)--(150:.4cm);
\node at (165:.5cm) {\tiny $F_{62}$};
\draw (270:.6cm)--(270:.4cm);
\node at (285:.5cm) {\tiny $F_{42}$};
\draw (90:.6cm)--(90:.4cm);
\node at (105:.5cm) {\tiny $F_{11}$};
\draw (210:.6cm)--(210:.4cm);
\node at (225:.5cm) {\tiny $F_{51}$};
\draw (330:.6cm)--(330:.4cm);
\node at (345:.5cm) {\tiny $F_{31}$};

\draw (120:1cm)--(140:.5cm);
\node at (140:.45cm) {\tiny $F_{61}$};
\draw (240:1cm)--(260:.5cm);
\node at (260:.45cm) {\tiny $F_{41}$};
\draw (0:1cm)--(20:.5cm);
\node at (20:.45cm) {\tiny $F_{21}$};

\node at (75:.45cm) {\tiny $F_{12}$};
\draw (65:1cm)--(75:.55cm);
\node at (-45:.45cm) {\tiny $F_{32}$};
\draw (-55:1cm)--(-45:.5cm);
\node at (195:.45cm) {\tiny $F_{52}$};
\draw (185:1cm)--(195:.5cm);

\draw [fill= white] (115:1 cm) circle [radius = .04];
\draw [fill= white] (65:1 cm) circle [radius = .04];
\draw [fill= white] (355:1 cm) circle [radius = .04];
\draw [fill= white] (-55:1 cm) circle [radius = .04];
\draw [fill= white] (185:1cm) circle [radius = .04];
\draw [fill= white] (235:1cm) circle [radius = .04];
\draw [fill= black] (90:.6 cm) circle [radius = .04];
\draw [fill= black] (210:.6 cm) circle [radius = .04];
\draw [fill= black] (-30:.6 cm) circle [radius = .04];
\draw [fill= white] (30:.6 cm) circle [radius = .04];
\draw [fill= white] (150:.6 cm) circle [radius = .04];
\draw [fill= white] (-90:.6 cm) circle [radius = .04];

\begin{scope}[xshift = 3.5cm]
\draw [black] (1,0) arc [radius = 1, start angle = 0, end angle = 360];

\foreach \s in {1,...,6}
{
  \node at ({360/6 * (-\s+6)-180}:1.3 cm) {$\s$};
  \draw [fill= black] ({360/6 * (\s - 1)}:1 cm) circle [radius = .04];
}

\draw (0:1cm)--(30:.6cm);
\draw (60:1cm)--(30:.6cm);
\draw (120:1cm)--(150:.6cm);
\draw (180:1cm)--(150:.6cm);
\draw (240:1cm)--(270:.6cm);
\draw (300:1cm)--(270:.6cm);

\draw (65:1cm)--(90:.6cm);
\draw (115:1cm)--(90:.6cm);
\draw (185:1cm)--(210:.6cm);
\draw (235:1cm)--(210:.6cm);
\draw (305:1cm)--(330:.6cm);
\draw (355:1cm)--(330:.6cm);

\draw (30:.6cm)--(30:.4cm);
\node at (17:.43cm) {\tiny $F_{12}$};
\draw (150:.6cm)--(150:.4cm);
\node at (137:.5cm) {\tiny $F_{52}$};
\draw (270:.6cm)--(270:.4cm);
\node at (257:.5cm) {\tiny $F_{32}$};
\draw (90:.6cm)--(90:.4cm);
\node at (80:.35cm) {\tiny $F_{61}$};
\draw (210:.6cm)--(210:.4cm);
\node at (200:.35cm) {\tiny $F_{41}$};
\draw (330:.6cm)--(330:.4cm);
\node at (320:.45cm) {\tiny $F_{21}$};

\draw (60:1cm)--(45:.5cm);
\node at (45:.45cm) {\tiny $F_{11}$};
\draw (180:1cm)--(165:.5cm);
\node at (165:.5cm) {\tiny $F_{51}$};
\draw (300:1cm)--(285:.5cm);
\node at (285:.45cm) {\tiny $F_{31}$};

\node at (110:.43cm) {\tiny $F_{62}$};
\draw (115:1cm)--(105:.55cm);
\node at (350:.43cm) {\tiny $F_{22}$};
\draw (-5:1cm)--(345:.53cm);
\node at (230:.45cm) {\tiny $F_{42}$};
\draw (235:1cm)--(220:.5cm);

\draw [fill= white] (115:1 cm) circle [radius = .04];
\draw [fill= white] (65:1 cm) circle [radius = .04];
\draw [fill= white] (355:1 cm) circle [radius = .04];
\draw [fill= white] (-55:1 cm) circle [radius = .04];
\draw [fill= white] (185:1cm) circle [radius = .04];
\draw [fill= white] (235:1cm) circle [radius = .04];
\draw [fill= black] (90:.6 cm) circle [radius = .04];
\draw [fill= black] (210:.6 cm) circle [radius = .04];
\draw [fill= black] (-30:.6 cm) circle [radius = .04];
\draw [fill= white] (30:.6 cm) circle [radius = .04];
\draw [fill= white] (150:.6 cm) circle [radius = .04];
\draw [fill= white] (-90:.6 cm) circle [radius = .04];
\end{scope}
\end{tikzpicture}
\caption{The first diagram shows how to compute $\Phi^* \rho^* \Psi^*([B])$ by reattaching a bivariant tensor diagram $B$ for $\Conf(3,6)$ to the boundary via the indicated connections. The second diagram does so for $\Phi^* (\rho^2)^* \Psi^*([B])$.  
\label{HowToDrawFig2}}
\end{center}
\end{figure}
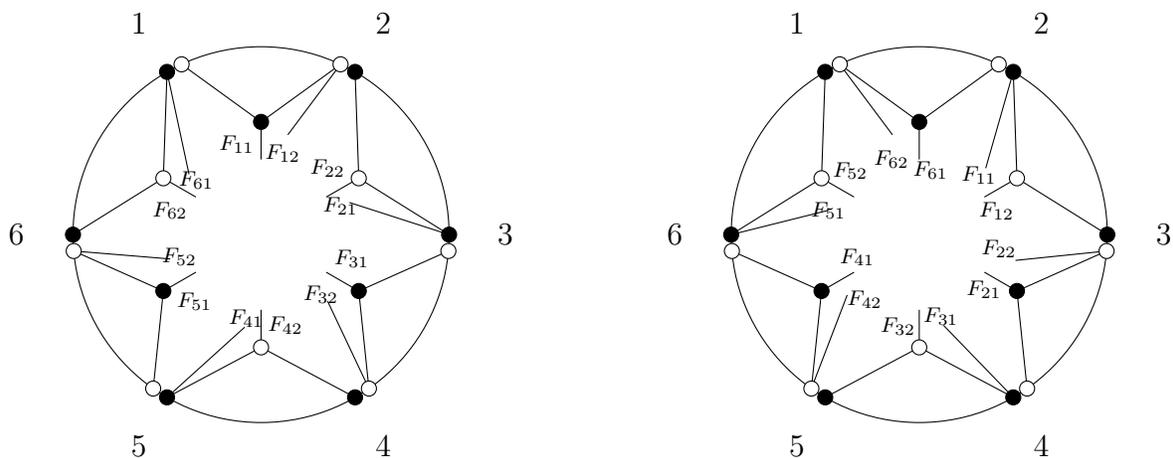

Establishing the presentation Theorem~\ref{thm:CMG3presentation} is more subtle, and for this part of the proof we work in $\FG(3,6)$ rather than $\tGr(3,9)$. 
On this side, $P^*$ acts by rotation, and the interesting functions are  $(\Psi \circ \rho \circ \Phi)^*$ and $(\Psi \circ \rho^2 \circ \Phi)^*$. Figure~\ref{HowToDrawFig2} and Example~\ref{eg:FGForks} shows how to evaluate these on tensor diagrams.

\begin{example}\label{eg:FGForks} Consider the Fock-Goncharov coordinate $\Delta_{2,3,5} \in \bbc[\Conf(3,6)]$ \eqref{eq:FGCoord}. As a biweb, $\Delta_{2,3,5}$ is a tripod joining the black parts of boundary vertices $2,3$ and $5$. We calculate $(\Phi^* \circ \rho^* \circ \Psi^*)(\Delta_{2,3,5})$ by creating a tripod on the strands labeled $F_{2,(1)},F_{3,(1)}$, and $F_{5,(1)}$ in the first diagram in Figure~\ref{HowToDrawFig2}. The result is a ``hexapod'' corresponding to the function  
\begin{equation}
(F_{1,\bullet},\dots,F_{6,\bullet}) \mapsto \omega^*\left(F_{3,(1)}(F_{3,(2)} \cap F_{4,(2)})(F_{5,(2)} \cap F_{6,(2)}) \right) \in \bbc[\Conf(3,6)].
\end{equation}
This hexapod has a boundary 3-cycle based at the third bi-colored boundary vertex $F_{3,\bullet}$. Applying the boundary 3-cycle skein relation produces a factorization of the right hand side of the above as $\omega^*\left(F_{3,(1)}F_{4,(2)}\right) \omega^*((F_{3,(2)} \cap F_{5,(2)})F_{6,(2)})$. The factor $\omega^*\left(F_{3,(1)}F_{4,(2)}\right)$ is a frozen bipod, and the factor $\omega^*((F_{3,(2)} \cap F_{5,(2)})F_{6,(2)})$ is a tripod on the white parts the boundary vertices $3,5$ and $6$. Denoting this latter tripod by $\Delta^{3,5,6}$, the dot action is
\begin{equation}\label{eq:upperindices}
(\Phi^* \circ \rho^* \circ \Psi^*) \cdot \Delta_{2,3,5} = \Delta^{3,5,6} \in \bbc[\Conf(3,6)].
\end{equation}

\end{example}

We call a fork in a (bivariant) tensor diagram \emph{black} (resp. \emph{white}) if the boundary vertices $i$ and $i+1$ it connects to are both black (resp. white). 
We call the fork \emph{even} (resp. \emph{odd}) if the boundary vertex $i$ is even (resp. \emph{odd}). Note that it can happen that a fork is both even and odd, but this happens precisely when the fork is the interior vertex of a tripod on three consecutive vertices. 

\begin{lem}\label{lem:HowToDrawII} Let $B$ be an indecomposable biweb invariant, not a frozen variable or a tripod on three consecutive vertices. Let $B' = \Phi^* \rho^* \Psi^* \cdot B$ and $B'' = \Phi^* (\rho^2)^* \Psi^* \cdot B$ be the dot actions of $B$ under the $\rho$ and $\rho^2$. By Lemma~\ref{lem:HowToDraw}, these biwebs are well-defined. Then if $B$ contains an odd black fork, then both of $B'$ and $B''$ contain an odd white fork. If $B$ has an even white fork, then both of $B',B''$ have an even black fork.  
\end{lem}

\begin{proof} To calculate $\Phi^* \rho^* \Psi^*([B])$, we take the first diagram in Figure~\ref{HowToDrawFig}, rotate it by one unit, and plug it in to the second diagram. The result is the first diagram in Figure~\ref{HowToDrawFig2}. When we plug in $B$ to Figure~\ref{HowToDrawFig2}, as argued in Lemma~\ref{lem:HowToDraw}, we obtain $B'$ by planarizing all the crossings that are introduced via the crossing removal skein relation, and then removing all boundary 3-cycles to remove factor out frozen variables. If $B$ has an odd black fork (e.g. $F_{1,(1)}F_{2,(1)})$, then the strand labeled by $F_{1,(1)}$ contributes an odd white fork. This odd black fork persists through the planarizing step and does not contribute to a boundary 3-cycle (on the contrary, even black forks $F_{6,(1)}F_{1,(1)}$ give rise to boundary 3-cycles). 

The analogous diagram for the square of the cyclic shift $\Phi^* (\rho^2)* \Psi^*([B])$ is also drawn in Figure~\ref{HowToDrawFig2}. The same style of argument for the three remaining cases (an even white fork for $\Phi^* \rho^* \Psi^*([B])$, as well as odd black forks and even white forks for $\Phi^* (\rho^2)* \Psi^*([B])$) completes the proof. 
\end{proof}


We make use of the following variant \cite{Passman} of the well-known Ping Pong Lemma. 

\begin{lem}\label{lem:PingPong}
Let $H$ be a group with nonidentity subgroups $H_1,H_2$, whose non-identity elements are denoted by $H_1^\#,H_2^\#$ respectively. 
Suppose that $H$ acts on a set $X$ having distinct nonempty subsets $X_1,X_2$, satisfying $H_1^\# X_1 \subseteq X_2, H_2^\# X_2 \subseteq X_1$. Finally, suppose $|H_2| \geq 3$. Then the group generated by $H_1,H_2$ inside $H$ is naturally isomorphic to the free product $H_1 * H_2$. 
\end{lem}

Now we establish our presentation for $\CMG(\Gr(3,9))$, and by the same methods, part (3) of Theorem~\ref{thm:GrThreeNine}.

\begin{proof}[Proof of Theorem~\ref{thm:CMG3presentation}]
We identify the cluster modular groups for $\Gr(3,9)$ and $\FG(3,6)$ via the maps $\Phi^*$ and $\Psi^*$, so we can think of the cyclic shift $\rho^*$ of vectors as an element of $\CMG(\FG(3,6))$ and so on. We elide the notational difference between a map and its pullback, writing $P$, rather than $P^*$, as an element of the cluster modular group. 

The main step is to study how the subgroup $\la \rho, P\ra$ acts on the Fock-Goncharov coordinate $\Delta_{1,3,5}$, which is a tripod on the black part of boundary vertices $1,3,5$. This cluster variable is fixed by $P^2$, so $B_3 / Z(B_3) \cong \PSL_2(\bbz) \cong \bbz / 2 \bbz \ast \bbz / 3 \bbz$ is the group that naturally acts on this orbit. We use the Ping Pong Lemma to show that this action of $\PSL_2(\bbz)$ is faithful. 

Let $X$ be the $\PSL_2(\bbz)$-orbit of $\Delta_{1,3,5}$ with respect to the dot action. Since every element of $X$ is stabilized by $P^2$ (and $P^2  = \rho^3$ commutes with $\rho$),  $X$ does not contain any tripods on three consecutive vertices. 

Next we define the Ping Pong sets $X_1,X_2 \subset X$. Let $X_1$ consist of those biwebs in $X$ that have no even white forks or odd black forks. Likewise, let $X_2$ consist of those webs in $X$ that have no even black forks or odd white forks. Notice that $\Delta_{1,3,5}$ is in both $X_1$ and $X_2$, so these two sets are nonempty. On the other hand, these two sets are distinct -- the element $\rho \cdot \Delta_{1,3,5}$, which is a hexapod 
$(F_{1,\bullet},\dots,F_{6,\bullet}) \mapsto (F_{1,(2)} \cap F_{2,(2)})(F_{3,(2)} \cap F_{4,(2)})(F_{5,(2)} \cap F_{6,(2)}),$ is in $X_1 \setminus X_2$. 

Clearly, $P(X_1) \subset X_2$ (in fact, $P(X_1) = X_2)$. It remains to show that $\rho (X_2) \subset X_1$ and $\rho^2 (X_2) \subset X_1$. Indeed, let $W \in X_2$ and let $W' = \rho \cdot W$ and $W'' = \rho^2 \cdot W$. By Lemma~\ref{lem:HowToDrawII}, if $W'$ had an odd black fork, then $\rho^2 \cdot W' = W$ must have an odd white fork, contradicting $W \in X_2$. Repeating this argument completes the proof. By Lemma~\ref{lem:PingPong}, the group $\la \rho,P \ra $ acts on $X$ as a free product  of cyclic groups of orders $2$ and $3$. 

Now let $G$ be the group with presentation given in the statement of the theorem. All of the relations in $G$ are satisfied in $\CMG$, so there is a homomorphism $G \to \CMG$. The group generated by the twist is normal, so this homomorphism descends to a well-defined map 
$\Aut(B_3) \cong G / \la \tau \ra \to \CMG / \la \tau \ra$. Here, we recall that the only outer automorphism of $B_3$ is the map inv$ \colon \sigma_i \mapsto {\sigma_i}^{-1}$. We have a short exact sequence $B_3 / Z(B_3) \injects \Aut(B_3) \surjects \la \text{inv}\ra$, and a similar short exact sequence $\la \rho,P \ra  \injects \CMG / \la \tau \ra \to \la \theta \ra$. We have shown that $B_3 / Z(B_3) \cong \la \rho,P \ra$. By the five lemma it follows that $\Aut(B_3) \cong \CMG / \la \tau \ra$. Finally, we ``unquotient'' by the center. We have a composition $G \surjects \CMG \surjects \Aut(B_3)$ whose kernel is generated by $\la \tau \ra$, and consists of $6$ elements. Thus $\ker(G \to \CMG) \subset \ker(G \to \Aut(B_3)) = \la \tau \ra$. It is easy to check on a particular cluster variable that the powers of $\tau$ are distinct elements of the cluster modular group. Thus the kernel of $G \to \CMG$ is trivial (i.e., the stated relations give a presentation).  
\end{proof}

\begin{proof}[Proof of Theorem \ref{thm:GrThreeNine} (3)] 
Let $B$ be one of the two \emph{single cycle biwebs} for~$\Conf(3,6)$. That is, $B$ is a biweb obtained by drawing a hexagon in the interior of the disk, and connecting the six vertices of the hexagon to the six boundary vertices with colors as demanded by bipartiteness. There are two such diagrams (rotations of each other). This is an indecomposable non-arborizable biweb. Notice that $B$ is stabilized by $P^2$ (it is also stabilized by reflection and the DT transformation, so the $\la \rho,P \ra$-orbit is the entire $\CMG$-orbit of $W$). Let $X_B$ be the orbit of $B$ with respect to the dot action of $\rho$ and $P$. Let $X_{1,B} \subset X_B$ and $X_{2,B} \subset X_B$ be the Ping Pong sets defined in terms of black/white  and even/odd forks as in the preceding proof. One checks that these two subsets are nonempty and disjoint. It follows that $\PSL_2(\bbz)$ acts freely on $X_B$, i.e. the orbit is infinite (in fact, the elements in the orbit are labeled by $\PSL_2(\bbz)$). Translating to~$\Gr(3,9)$ via the quasi-isomorphism, the same is true for the single cycle web in Figure~\ref{fig:SCW}. 
\end{proof}

\subsection{Proofs for $\Gr(4,8)$}
We prove Theorem~\ref{thm:GrFourEight}, and with quite a bit more work, prove Theorem \ref{thm:CMG4presentation}.

\begin{proof}[Proof of Theorem~\ref{thm:GrFourEight}] 
The argument is the same as the proof of parts (1) and (2) of Theorem~\ref{thm:GrThreeNine}. On the fundamental domain, every cluster variable is a web that is already in tree form. One checks that the renormalized Artin generator $\tilde{\sigma_1}^*$ (cf.~Remark~\ref{rmk:renormd}) preserves the set of $\SL_4$ web invariants. As before, this  is not a priori obvious, since evaluating $\tilde{\sigma_1}^*$ by ``reattaching'' (i.e., the $\SL_4$ analogue of Figures~\ref{HowToDrawFig} and \ref{HowToDrawFig2}) introduces crossings. However, all such crossings can be planarized by a combination of the  $\SL_4$ crossing removal and degeneracy relations. On the other hand, evaluating $\tilde{\sigma_1}^*$ only involves reattaching certain boundary edges to the boundary in a new way via certain trees, and thus sends tree invariants to tree invariants. Since every cluster variable on the fundamental domain is a web and a tree, the same is true for every cluster variable (up to a monomial in the frozens). The statement about Pl\"uckers and octapods is true because it holds on a fundamental domain, as noted in the proof to Lemma~\ref{lem:computercheck}.   
\end{proof}

The proof of Theorem~\ref{thm:CMG4presentation} is more ``hands on'' than the proof of Theorem~\ref{thm:CMG3presentation} for several reasons. First, there is no basis of $\SL_4$ web invariants, which makes equality testing hard. Second, there are proportionality relations amongst Artin generators that are difficult to see at the level of tensor diagrams (cf.~Remark~\ref{rmk:exotic}). Third, letting $\Delta \in B_4$ denote the half-twist braid, the key step in our proof of Theorem~\ref{thm:CMG4presentation} is to show that $\CMG^+/\la \Delta^2\ra = \MCG(S^2,4)$. After this, we ``unquotient'' by the center as was done in the proof of Theorem \ref{thm:CMG3presentation}. However, unlike in that proof, in the current situation, to establish $\CMG^+/\la \Delta^2\ra = \MCG(S^2,4)$ it is not enough to consider the orbit of a single cluster variable (which will have a nontrivial stabilizer), but rather to consider several cluster variables at a time. 

We let $\Delta \in B_4$ denote the half-twist, and consider also the braids $\aa = \sigma_2 \sigma_1, \bb = \sigma_2 \sigma_1 \sigma_2$. The next lemma is well known: 
\begin{lem}
The elements $\aa,\bb$ satisfy $\aa^3=\bb^2 =1$ in the quotient $B_4 / \la \Delta^2, \sigma_1 \sigma_2 \sigma_3^2 \sigma_2 \sigma_1\ra = \MCG(S^2,4)$. They generate a free product $\bbz / 2 \bbz \ast \bbz / 3 \bbz  \subset \MCG(S^2,4)$, of index four in $\MCG(S^2,4)$, with the braids $\{1,\sigma_3,\sigma_3 \sigma_2, \sigma_3 \sigma_2 \sigma_1\}$ serving as right coset representatives.  
\end{lem}

\begin{proof}
The proof is geometric. For any $k$, there is an inclusion  $B_k / Z(B_k) \subset \MCG(S^2,k+1)$ of the braid group modulo its center into the mapping class group of a sphere with $k+1$ marked points, described as follows. The braid group $B_k$ can be thought of as the mapping class group of a closed disk $D$ with $k$ punctures. Let $D'$ be a once-punctured disk. Gluing $D$ and $D'$ along their boundary produces a sphere $S^2$ with $k+1$ punctures. We view the $k+1$st puncture as sitting at the North Pole. The inclusion of spaces $D \subset S^2$ induces a homomorphism $B_k \to \MCG(S^2,k+1)$ of mapping class groups. As a mapping class, the full twist $\Delta^2 \in B_k$ corresponds to a Dehn twist along the boundary of the disk $D$, and this mapping class represents a trivial element of $\MCG(S^2,k+1)$. The homomorphism factors $B_k/Z(B_k) \to \MCG(S^2,k+1)$ to the quotient by $\Delta^2$, and this latter homomorphism is injective (as follows from the \emph{Birman exact sequence} \cite{FarbMargalit}). The subgroup $B_k / Z(B_k) \subset \MCG(S^2,k+1)$ is identified with the normal subgroup of mapping classes that fix the North Pole. In our case of interest, $k=3$, we get a copy of $B_3 / Z(B_3) = \PSL_2(\bbz) \subset \MCG(S^2,4)$, generated by the elements $\aa = \sigma_2\sigma_1$ and $\bb = \sigma_1 \sigma_2 \sigma_1$. There is a surjection $\MCG(S^2,4) \to \mfs_4$ (sending $\sigma_i$ to the transposition $(i, i+1)$), and $\PSL_2(\bbz) \subset \MCG(S^2,4)$ is the inverse image of $\mfs_3 \subset \mfs_4$. It follows that $\PSL_2(\bbz)$ has index $4$ in $\MCG(S^2,4)$, and that the four braids above serve as coset representatives.  
\end{proof}


\begin{lem} Let $w \in \bbz / 2 \bbz \ast \bbz / 3 \bbz $ be a (reduced) word in the generators $\aa,\bb$. If $w$ defines a trivial element of the cluster modular group, then $w$ is conjugate to a power of $\bb \aa \bb \aa$.
\end{lem}

\begin{proof}
We prove this by considering the $\la \aa,\bb \ra$-orbit of the Pl\"ucker coordinate $\Delta_{2468}$.  We let $f_\aa, f_\bb \in \Aut(\tGr^\circ(4,8))$ be the automorphisms corresponding to the braids $\aa$ and $\bb$, and let $f_w$ be the corresponding composition of $f_\aa$ and $f_\bb$. Since both $f_\aa$ and $f_\bb$ commute with the cyclic shift $\rho^4$, this is also true of $f_w$. 

If $w$ is trivial in the cluster modular group, then it fixes each cluster variable up to a frozen Laurent monomial. In particular, $f_w^*(\Delta_{2468}) = M \Delta_{2468} $ for a frozen Laurent monomial $M$. This Laurent monomial is necessarily $\rho^4$ invariant, e.g. if $\Delta_{1234}$ appears in the numerator of $M$, then $\Delta_{5678}$ also appears in the numerator of $M$, and so on. It follows that the homogeneous degree of the monomial $M$ with respect to the $\bbz^8$-grading on $\bbc[\tGr^\circ(4,8)]$ is an integer multiple of the all-ones vector $(1,\dots,1) \in \bbz^8$. If $e_1,\dots,e_8$ are standard basis vectors for $\bbz^8$, then 
$f_w^*(\Delta_{2468})$ has degree $e_2+e_4+e_6+e_8+m(1,\dots,1)$ for some $m \in \bbz$. Recall that $\sigma_1^*$ acts by permuting $e_1 \leftrightarrow e_2$ and $e_5 \leftrightarrow e_6$, while $\sigma_2^*$ acts by permuting $e_2 \leftrightarrow e_3$ and $e_6 \leftrightarrow e_7$ (cf.~the formula for $\sigma_1(\vec{v} \cdot \mathbf{t})$ before the statement of Theorem~\ref{thm:braidgroupacts}). Thus $\sigma_1^*$ and $\sigma_2^*$, and therefore $f_\aa^*$ and $f_\bb^*$, preserve the all-ones vector. One calculates the following action of $f_\aa^*$ and $f_\bb^*$ on $e_2+e_4+e_6+e_8 \in \bbz^8 / (1,\dots,1)$: 
\begin{equation}\label{eq:S3onZ8}
\xymatrix{ e_2+e_4+e_6+e_8 \ar@{->}^{f_\aa^*}[dr] \ar@(ul,ur)^{f_\bb^*} &   \\
e_1 + e_4 + e_5 + e_8 \ar@{->}^{f_\aa^*}[u] & e_3 + e_4 + e_7 + e_8 \ar@{<->}^{f_\bb^*}[l] \ar@<-1ex>_{f_\aa^*}[l]
}.
\end{equation}
That is, modulo the all ones vector, the group $\PSL_2(\bbz)$ permutes the three vectors in \eqref{eq:S3onZ8}. We get a homomorphism $\PSL_2(\bbz) \twoheadrightarrow \mfs_3$ in which $\aa$ acts by a 3-cycle and $\bb$ acts by a transposition. The symmetric group $\mfs_3$ has a dihedral presentation $\la \aa,\bb \colon \aa^3 = \bb^2 = 1, \: \bb \aa \bb = \aa^{-1}. \ra$, so ker$(\PSL_2(\bbz) \twoheadrightarrow \mfs_3)$ is the normal subgroup generated by $\bb \aa \bb \aa.$. If $w$ acts trivially in the cluster modular group then it lies in this kernel, i.e. $w$ is conjugate to a power of  $\bb \aa \bb \aa$. 
\end{proof}

\begin{lem} The quasi-automorphism determined by the braid $\bb \aa \bb \aa$ has infinite order in $\CMG(\tGr(4,8))$. 
\end{lem}

\begin{figure}
\begin{tikzpicture}
\node at (-6,3) {\begin{tikzpicture}[scale=.5]
\xymatrixrowsep{.2in}
\xymatrix{
&\Delta_{1278} \ar[d]&\Delta_{1678} \ar[d]&\Delta_{5678}& \\
&\Delta_{1378} \ar[ur] \ar[d] &\Delta_{3678} \ar[l] \ar[r] \ar[d]&\Delta_{4678} \ar[d] \ar[u]& \\
\Delta_{1238} \ar[dr]&\Delta_{2378} \ar[l] \ar[ur] \ar[dr]&\Delta_{2367} \ar[l] \ar[r]&\Delta_{3467} \ar[ul] \ar[dl] \ar[r]&\Delta_{4567} \ar[ul]\\
&\Delta_{2348} \ar[u] \ar[d]&\Delta_{2347} \ar[u] \ar[l] \ar[r]&\Delta_{3457} \ar[u] \ar[dl]& \\
&\Delta_{1234}&\Delta_{2345} \ar[u]&\Delta_{3456} \ar[u]& 
}
\end{tikzpicture}};
\node at (-4.5,-4.5) {
\xymatrix{
\boxed{\Delta_{1378}} \ar[d] & \Delta_{3678} \ar[l] \ar[r] & \boxed{\Delta_{4678}} \ar[d]  \\
\Delta_{2378}  \ar[ur] \ar[dr]&&\Delta_{3467} \ar[ul] \ar[dl] \\
\boxed{\Delta_{2348}} \ar[u] &\Delta_{2347} \ar[l] \ar[r]&\boxed{\Delta_{3457}} \ar[u]  \\
}
};

\begin{scope}[yshift= -4.5cm,xshift = 1cm]
\def  \rsizeo{.3};
\def  \rsizet{1.5};
\draw [fill= lightgray] (0,0) circle [radius = \rsizeo];
\draw (\rsizeo,0) arc [radius = \rsizeo , start angle = 0, end angle = 360];
\draw (\rsizet,0) arc [radius = \rsizet , start angle = 0, end angle = 360];
\draw [thick] (0,\rsizeo)--(0,\rsizet);
\draw [thick] (0,-\rsizeo)--(0,-\rsizet);
\draw [thick] (-90: \rsizeo) to [out = -40, in = 270] (0:.6*\rsizet) to [out = 90, in = -40] (90: \rsizet);
\draw [thick] (90: \rsizeo) to [out = 150, in = 90] (180:.6*\rsizet) to [out = 270, in = 150] (-90: \rsizet);
\node at (100:.8*\rsizet) {$a$};
\node at (-80:.75*\rsizet) {$b$};
\node at (50:.85*\rsizet) {$c$};
\node at (-155:.85*\rsizet) {$d$};
\node at (0:1.15*\rsizet) {$e$};
\node at (180:1.15*\rsizet) {$f$};
\node at (0:1.5*\rsizeo) {$g$};
\node at (180:1.5*\rsizeo) {$h$};
\draw [fill= black] (0,-\rsizet) circle [radius = .07];
\draw [fill= black] (0,-\rsizeo) circle [radius = .07];
\draw [fill= black] (0,\rsizet) circle [radius = .07];
\draw [fill= black] (0,\rsizeo) circle [radius = .07];
\node at (3,2) {$a : \: \Delta_{2378}$};
\node at (3,1.5) {$b : \: \Delta_{3467}$};
\node at (3,1) {$c : \: \Delta_{3678}$};
\node at (3,.5) {$d : \: \Delta_{2347}$};
\node at (3,0) {$e : \: \Delta_{4678}$};
\node at (3,-.5) {$f : \: \Delta_{2348}$};
\node at (3,-1) {$g : \: \Delta_{1378}$};
\node at (3,-1.5) {$h : \: \Delta_{3457} $};
\end{scope}
\end{tikzpicture}
\caption{The extended quiver for a seed in $\Gr(4,8)$. Below it, an initial seed for the subalgebra $\mca'$ obtained by freezing certain variables and setting certain frozen variables equal to $1$. To the right, we show the triangulation of the annulus that corresponds to this quiver. The element $w = \bb \aa \bb \aa$ acts on this cluster subalgebra by twisting the outer boundary one unit clockwise.\label{fig:annulus}} 
\end{figure}
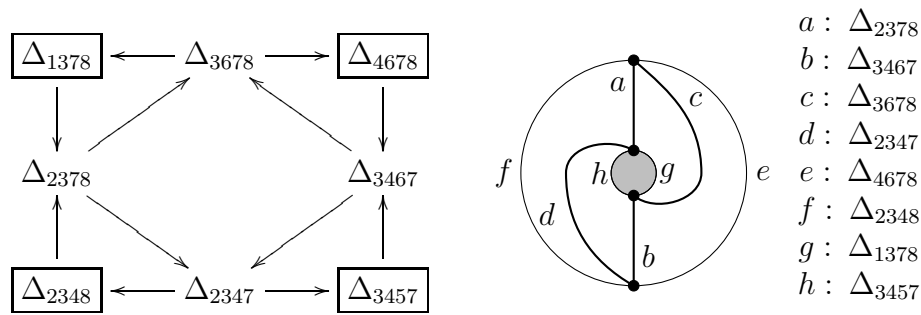

\begin{proof}
Using the braid relations one finds that $\bb \aa \bb \aa = \aa^3 \sigma_1 \sigma_2^2 \sigma_1 \in B_4$, and thus $\bb \aa \bb \aa = \sigma_1 \sigma_2^2 \sigma_1 \in \MCG(S^2,4)$. So it suffices to show that the quasi-automorphism $(\sigma_1 \sigma_2^2 \sigma_1)^*$ has infinite order in the cluster modular group. 

By direct calculations with the $\sigma_i^*$, each of the Pl\"ucker coordinates $\Delta_{1378},\Delta_{2367},\Delta_{3457}$ is fixed by $(\sigma_1 \sigma_2^2 \sigma_1)^*$ up to frozen variables. Likewise, the Pl\"ucker coordinates $\Delta_{4678} \leftrightarrow \Delta_{2348}$ are swapped by $(\sigma_1 \sigma_2^2 \sigma_1)^*$ up to frozen variables. For an example calculation, one sees that $\sigma_1^*(\Delta_{1378}) = \Delta_{2378}$. Next, the renormalized Artin generator satisfies  $\tilde{\sigma_2^*}(\Delta_{2378}) = \omega^*(v_2v_3(v_6v_7 \cap v_8v_1v_2) v_8) = \omega^*(v_1v_2v_3v_8) \omega^*(v_2v_6v_7v_8) \propto \Delta_{2678}$ by a similar argument as in the proof of Proposition~\ref{prop:upandown}. Next, $\sigma_2^*(\Delta_{2678}) = \Delta_{3678}$,  and finally $\tilde{\sigma_1}^2(\Delta_{3678}) = \omega^*(v_3 (v_5v_6 \cap v_7 v_8v_1)v_7v_8) = \omega^*(v_5v_6v_7v_8)\omega^*(v_1v_3v_7v_8) \propto \Delta_{1378}$ as claimed. The other calculations have a similar flavor. 

We have described the action of $(\sigma_1 \sigma_2^2 \sigma_1)^*$ on the Pl\"ucker coordinates
\begin{equation}
\Delta_{1378},\Delta_{2348},\Delta_{2367},\Delta_{4678},\Delta_{3457}. \label{eq:fivefixed}
\end{equation}
These are weakly separated, and can be extended to a cluster $\mcc$ for $\bbc[\tGr(4,8)]$ by adding $\Delta_{2347},\Delta_{2378},\Delta_{3678},\Delta_{3467}$.
The extended quiver $\tQ(\mcc)$ is in Figure \ref{fig:annulus}. We let $X \in \bbc[\tGr(4,8)]$ be the element $X((v_1,\dots,v_8)) = \omega^*(v_7 v_8 v_1 (v_4  v_6 v_7 \cap v_2 v_3))$, and $Y \in \bbc[\tGr(4,8)]$ be the element $Y((v_1,\dots,v_8)) = \omega^*(v_3 v_4 v_5 (v_2  v_3 v_8 \cap v_6 v_7))$. By another direct calculation, one sees that $(\sigma_1 \sigma_2^2 \sigma_1)^*$ acts on the four extra Pl\"ucker coordinates in $\mcc$ by
\begin{align}
\Delta_{2347} \mapsto \Delta_{2378} \hspace{1cm} &\Delta_{3678} \mapsto \Delta_{3467}  \label{eq:octapod1} \\
\Delta_{2378} \mapsto X \hspace{1cm}  &\Delta_{3467} \mapsto Y \label{eq:octapod2}
\end{align}
The exotic-looking cluster variables $X$ and $Y$ are both a single mutation away from the cluster in Figure~\ref{fig:annulus}. Specifically, $\mu_\mcc(\Delta_{2347}) = X$ and $\mu_\mcc(\Delta_{3678}) = Y$, which one sees by checking the corresponding exchange relation \eqref{eq:exchangereln}. 

Now we freeze the five variables in \eqref{eq:fivefixed}, obtaining a cluster subalgebra of $\bbc[\tGr(4,8)]$ with $8+5 = 13$ frozen variables and clusters with $4$ mutable variables. We set all 8 of the original frozen variables equal to $1$, and also set $\Delta_{2367} =1$. We are left with a cluster algebra $\mca'$ with four frozen variables. Our choice of initial seed for this cluster subalgebra is pictured in Figure \ref{fig:annulus}. This cluster algebra $\mca'$ is the cluster algebra associated to an annulus with $2$ points on each boundary component, and with four frozen variables given by the boundary arcs. Our initial seed corresponds to the triangulation of the annulus pictured in Figure \ref{fig:annulus}. From our description of how $f_{\sigma_1 \sigma_2^2 \sigma_1}$ acts on the four frozen variables, as well as the formulas \eqref{eq:octapod1} and \eqref{eq:octapod2}, we see that $\sigma_1 \sigma_2^2 \sigma_1$ acts on the annulus by rotating the outer boundary one unit clockwise. Rotating the outer boundary has infinite order when thought of as a map on clusters in $\mca'$. It follows that $(\sigma_1 \sigma_2^2 \sigma_1)^*$ has infinite order in $\CMG(\Gr(4,8))$. 
\end{proof}

\begin{proof}[Proof of Theorem \ref{thm:CMG4presentation}]
Let $G$ be the group with the presentation in the theorem statement. We have a surjection $G \to \CMG$ by Lemma~\ref{lem:computercheck}. It descends to the quotient by $\Delta^2$. By the two previous lemmas, we have an injection $\PSL_2(\bbz) \hookrightarrow \CMG / \la \Delta^2 \ra$, where $\PSL_2(\bbz) = \la \aa,\bb\ra$. We claim that this extends to an injection $\MCG(S^2,4) \hookrightarrow \CMG / \la \Delta^2 \ra$. Let $K$ be the kernel of this homomorphism. Conjugation by $\PSL_2(\bbz)$ acts on the four cosets of $\PSL_2(\bbz) \subset \MCG(S^2,4)$, fixing the identity coset. The other three cosets are transitively permuted by conjugation (they are permuted by $\aa$). Thus, it suffices to prove that the coset $\PSL_2(\bbz)\sigma_3 \sigma_2 \sigma_1$ intersects $K$ trivially.  

Remember that $\sigma_3 \sigma_2 \sigma_1 \propto \rho$ on $\Gr(4,8)$. Suppose $w = w' \rho \in K$ with $w' \in \la \aa,\bb \ra$. Then $f_w^*(\Delta_{1357}) \propto \Delta_{1357}$, so $f_{w'}^*(\Delta_{1357}) \propto \Delta_{2468}$. Thus, $\Delta_{1357}$ would be in the $\PSL_2(\bbz)$-orbit of $\Delta_{2468}$, but this is not possible from the calculation with the $\bbz^8$-grading \eqref{eq:S3onZ8}. One now argues that $G \to \CMG$ is an isomorphism by ``unquotienting'' the map  
$\MCG(S^2,4) \hookrightarrow \CMG / \la \Delta^2 \ra$ in a similar way as was done at the end of the proof of \ref{thm:CMG3presentation}.
\end{proof}

\end{document}